\documentclass[11pt]{article}
\usepackage[T1]{fontenc}
\usepackage[margin=1in]{geometry}
\usepackage{lmodern,amsmath,amssymb,amsthm,mathtools,bm,mathrsfs}
\usepackage{aliascnt,enumitem,graphicx,booktabs,tabularx,array,microtype}
\usepackage{xcolor,pgfplots}
\pgfplotsset{compat=1.18}
\definecolor{cA}{HTML}{0072B2}\definecolor{cB}{HTML}{D55E00}
\definecolor{cC}{HTML}{009E73}\definecolor{cD}{HTML}{CC79A7}
\usepackage{hyperref}
\hypersetup{hidelinks,hypertexnames=false}
\usepackage[nameinlink,capitalize]{cleveref}
\usepackage[bottom]{footmisc}
\renewcommand{\footnotesize}{\scriptsize}
\usepackage{caption}
\usepackage[compact]{titlesec}
\allowdisplaybreaks[2]
\AtBeginDocument{%
  \setlength{\abovedisplayskip}{7pt plus 2pt minus 3pt}%
  \setlength{\belowdisplayskip}{7pt plus 2pt minus 3pt}%
  \setlength{\abovedisplayshortskip}{4pt plus 2pt minus 2pt}%
  \setlength{\belowdisplayshortskip}{4pt plus 2pt minus 2pt}%
}

\titlespacing*{\section}
  {0pt}{1.55ex plus 0.35ex minus 0.25ex}{0.75ex plus 0.15ex}

\titlespacing*{\subsection}
  {0pt}{1.25ex plus 0.30ex minus 0.20ex}{0.55ex plus 0.12ex}

\titlespacing*{\subsubsection}
  {0pt}{1.00ex plus 0.25ex minus 0.15ex}{0.45ex plus 0.10ex}

\titleformat{\section}
  {\normalfont\large\bfseries}{\thesection.}{0.5em}{}

\titleformat{\subsection}
  {\normalfont\normalsize\bfseries}{\thesubsection.}{0.5em}{}

\titleformat{\subsubsection}
  {\normalfont\normalsize\itshape}{\thesubsubsection.}{0.5em}{}

\setlist{
  topsep=2pt plus 1pt minus 1pt,
  itemsep=1pt plus 0.5pt,
  parsep=0pt,
  partopsep=0pt,
  leftmargin=*
}

\makeatletter
\def\thm@space@setup{%
  \thm@preskip=6pt plus 2pt minus 2pt
  \thm@postskip=6pt plus 2pt minus 2pt
}
\makeatother

\makeatletter
\renewcommand{\maketitle}{%
  \begin{center}
    {\Large\bfseries \@title\par}
    \vskip 0.45em
    {\normalsize
      \lineskip .5em
      \begin{tabular}[t]{c}
        \@author
      \end{tabular}\par}
    \vskip 0.35em
    {\small \@date\par}
  \end{center}
  \par\vskip 0.35em
  \@thanks
  \vskip 0.7em
}
\makeatother

\renewenvironment{abstract}
  {\small\begin{center}\bfseries Abstract\end{center}\vspace{-0.6em}\small}
  {\par\vspace{0.6em}}

\newtheorem{theorem}{Theorem}[section]

\newaliascnt{lemma}{theorem}
\newtheorem{lemma}[lemma]{Lemma}
\aliascntresetthe{lemma}

\newaliascnt{proposition}{theorem}
\newtheorem{proposition}[proposition]{Proposition}
\aliascntresetthe{proposition}

\newaliascnt{corollary}{theorem}
\newtheorem{corollary}[corollary]{Corollary}
\aliascntresetthe{corollary}

\theoremstyle{definition}

\newaliascnt{definition}{theorem}

\aliascntresetthe{definition}

\newaliascnt{assumption}{theorem}

\aliascntresetthe{assumption}

\theoremstyle{remark}

\newaliascnt{remark}{theorem}
\newtheorem{remark}[remark]{Remark}
\aliascntresetthe{remark}

\crefname{theorem}{theorem}{theorems}
\Crefname{theorem}{Theorem}{Theorems}

\crefname{lemma}{lemma}{lemmas}
\Crefname{lemma}{Lemma}{Lemmas}

\crefname{proposition}{proposition}{propositions}
\Crefname{proposition}{Proposition}{Propositions}

\crefname{corollary}{corollary}{corollaries}
\Crefname{corollary}{Corollary}{Corollaries}

\crefname{definition}{definition}{definitions}
\Crefname{definition}{Definition}{Definitions}

\crefname{assumption}{assumption}{assumptions}
\Crefname{assumption}{Assumption}{Assumptions}

\crefname{remark}{remark}{remarks}
\Crefname{remark}{Remark}{Remarks}

\newcounter{prooflink}
\crefname{prooflink}{proof}{proofs}
\Crefname{prooflink}{Proof}{Proofs}
\newcommand{\prooflink}[1]{\par\noindent{\small\hyperref[proof:#1]{Proof in the appendix.}}}

\DeclareMathOperator{\Tr}{Tr}
\DeclareMathOperator{\diag}{diag}
\DeclareMathOperator{\spec}{spec}

\newcommand{\R}{\mathbb R}\newcommand{\C}{\mathbb C}
\newcommand{\N}{\mathbb N}\newcommand{\Z}{\mathbb Z}
\newcommand{\E}{\mathcal E}\newcommand{\cJ}{\mathcal J}
\newcommand{\dd}{\,\mathrm d}
\newcommand{\ind}{\mathbf 1}

\newcommand{\sech}{\operatorname{sech}}
\newcommand{\ip}[2]{\langle #1,#2\rangle}
\newcommand{\nrm}[1]{\lVert #1\rVert}
\newcommand{\fm}{\mathfrak m}\newcommand{\fr}{\mathfrak r}
\newcommand{\fn}{\mathfrak n}\newcommand{\cL}{\mathcal L}
\newcommand{\cB}{\mathcal B}\newcommand{\cK}{\mathcal K}
\newcommand{\cF}{\mathcal F}\newcommand{\cG}{\mathcal G}
\newcommand{\sK}{\mathsf K}
\newcommand{\hD}{\widehat D}\newcommand{\hd}{\hat{\mathfrak d}}
\newcommand{\fD}{\mathfrak D}\newcommand{\fE}{\mathfrak E}
\newcommand{\cO}{\mathcal O}\newcommand{\rI}{\mathrm i}

\renewcommand{\Re}{\operatorname{Re}}\renewcommand{\Im}{\operatorname{Im}}
\numberwithin{equation}{section}

\setlist[enumerate]{label=\textup{(\alph*)},leftmargin=2.3em}
\title{Resolvent-Diagonal Energy of Isospectral Jacobi Matrices}
\author{Dmitry Pasechnyuk-Vilensky\thanks{Mohamed bin Zayed University of Artificial Intelligence, Abu Dhabi, UAE. \texttt{dmivilensky1@gmail.com}}, Martin Tak\'{a}\v{c}\thanks{Mohamed bin Zayed University of Artificial Intelligence, Abu Dhabi, UAE. \texttt{Takac.MT@gmail.com}}}
\date{}
\begin{document}
\maketitle
\begin{abstract}
We study the minimum squared resolvent diagonal among Jacobi matrices with prescribed eigenvalues. For the zeros of $P_n^{(\alpha,\beta)}$, we construct a canonical realization and bound its asymptotic cost against the optimum. Its rescaled inverse converges to the Green operator of a reflection-symmetric singular Schr\"odinger operator with Bessel-zero spectrum. A Green-coordinate identity places this construction and all finite realizations in the same density framework. Mass-constrained spectral problems then give computable lower bounds. The bounds meet at $\alpha=1/2$, with optimal energy $2/15$, and have first-order contact with derivative $-13/75$. We classify the constrained extremal profiles and the three endpoint geometries of the canonical density. Two polynomial families with the same limiting edge spectrum have different limiting optimal energies. For Nesterov momentum $(s-1)/(s+1)$, the optimal leading endpoint Hessian-drift gain among nonnegative-link Jacobi realizations is $\sqrt2$ times the original gain. An explicit isospectral flow strictly improves the standard one-sided realization for $\alpha>-1/2$, $\alpha\ne1/2$.
\end{abstract}

\section*{Introduction}

A polynomial with simple real zeros determines the spectrum of a Jacobi matrix. Its positive spectral weights determine the matrix itself, and therefore the distribution of the diagonal of its resolvent. We seek a realization that distributes this fixed trace with minimum squared diagonal. For the zeros of a Jacobi polynomial the asymptotic problem is
\[
 e(\alpha,\beta)=\liminf_{n\to\infty}\ \inf_{J\in\cJ(P_n^{(\alpha,\beta)})}
 n^{-3}\sum_{i=0}^{n-1}[(I-J)^{-1}]_{ii}^2,
 \qquad \alpha,\beta>-1,
\]
where $\cJ(p)$ consists of real symmetric tridiagonal matrices with positive off-diagonal entries and eigenvalues equal to the zeros of $p$. We construct a canonical realization and assess its cost by lower bounds valid throughout the fibre.

This question also arises in accelerated iterations. A terminal residual determines the output on a fixed quadratic, while spectral weights specify the intermediate recurrence. The response to a Hessian perturbation at step $i$ factors through the two principal blocks on either side of that vertex. At the spectral endpoint this cofactor is a resolvent diagonal entry. The energy therefore measures an endpoint channel of the terminal first variation. Parts~I--II \cite{PT1,PT2} develop this connection for the Chebyshev minimax residual; \cref{prop:drift-energy} gives the identity for general Jacobi realizations.

The canonical choice is the persymmetric matrix, whose weights are proportional to $|p'(\zeta_k)|^{-1}$. Our main construction follows these finite weights to a singular Schr\"odinger operator with spectrum $\{j_{\alpha,k}^2\}$. We identify its boundary branch, its midpoint conditions and the limiting Green diagonal, obtaining an asymptotic energy $c_*(\alpha)$ independent of $\beta$. Uniform control of both spectral edges allows the diagonal and its squared integral to pass to the limit. The Green coordinate then represents the same operator by a density on the line, preserving the full spectrum, the energy and both scalar normalizations.

The lower construction starts from the packing of Part~II. Every finite realization gives a density with the same cost and asymptotic spectral data. The sharp second-moment inequality already gives the exact answer at $U$, the index $\alpha=1/2$. Away from this point, retaining the mass, the largest eigenvalue, the regularized determinant or the root mass yields constrained variational problems. We classify their extremal profiles and obtain
\[
 \ell(\alpha)\le e(\alpha,\beta)\le c_*(\alpha),\qquad
 e(1/2,\beta)=2/15,\qquad
 \partial_\alpha e(1/2,\beta)=-13/75.
\]
The bounds have first-order contact at $U$. Their numerical comparison measures the quality of the canonical choice across the Jacobi family.

Two consequences explain the information carried by these constructions. At momentum $(s-1)/(s+1)$, the exact $U$ energy determines the optimal leading endpoint drift gain among all nonnegative-link Jacobi realizations of the Nesterov residual, including arbitrary restarts and both horizon parities. The original Nesterov recurrence has the same terminal regularization output and a gain smaller by a factor $\sqrt2$. More generally, the standard Jacobi matrix realizes the one-sided Bessel geometry, and an explicit flow of its spectral weights gives a strict asymptotic improvement for $\alpha>-1/2$, $\alpha\ne1/2$. A second example gives two polynomial families with identical limiting edge spectra and different limiting optimal energies. These results show how the finite realization enters both algorithmic sensitivity and asymptotic optimization.

Finite spectral reconstruction and persymmetric uniqueness are classical \cite{Hoch67,Hald76,dBG,dBS86,GTVZ17}. We use Gel'fand--Levitan reconstruction and its discrete Cholesky form \cite{GeLe51,Mar86,DBK05} to identify the canonical limit and control its energy. The constrained classifications build on sharp one-dimensional interpolation and eigenvalue inequalities.

The operator method of moments connects polynomial recurrences, finite operator approximations and acceleration of linear iterations \cite[Chapters~III--IV]{Vor58}; its Lanczos moment-matching formulation is reviewed in \cite[\S3]{Str09}. Here the eigenvalues are fixed while the spectral weights vary to control the resolvent-diagonal energy.

The density framework below also has a classical mechanical reading: under $\xi=\tanh x$ the kernel $\sqrt{m(x)}e^{-|x-y|}\sqrt{m(y)}$ becomes the Green operator of an inhomogeneous string on $[-1,1]$ with fixed ends. This is the Liouville transformation that Beals, Sattinger and Szmigielski used for the multipeakons of the Camassa--Holm equation \cite[\S3]{BSS00} (\cref{rem:string}). The dictionary runs as follows:
\begin{itemize}
 \item finite realizations become beaded strings in the sense of Gantmacher and Krein \cite{GK02} (\cref{rem:beads});
 \item the canonical limit becomes a symmetric string with the Bessel-zero spectrum;
 \item for $\alpha>0$, the standard realization becomes a power-law string with the same spectrum (\cref{rem:power-string});
 \item the resolvent flow of \cref{ss:flow} becomes the multipeakon flow (\cref{rem:peakons}).
\end{itemize}
In this language the energy problem is a weighted isoperimetric problem for strings with prescribed spectrum, in the tradition of Krein \cite{Kr51a}.

\paragraph{Formal verification.}
The mathematical results are accompanied by formal proofs in Lean. The package is available at
\href{https://github.com/dmivilensky/PolynomialRealizationsProgram/tree/main/jacobi}{\texttt{PolynomialRealizationsProgram/jacobi}}\footnote{\url{https://github.com/dmivilensky/PolynomialRealizationsProgram/tree/main/jacobi}}.

\paragraph{Notation.}
Integrals display their domain and differential; $\dd s$ denotes Lebesgue measure, and other measures are specified explicitly. The spaces $L^p(\Omega,\dd s)$ and $H^1(\Omega)$ use Lebesgue measure; $H^1_0(\Omega)$ is the closure of $C_c^\infty(\Omega)$ in $H^1(\Omega)$. Function spaces are real unless complex parameters are under consideration. For functions, $\nrm{\cdot}_p$ denotes the $L^p$ norm on the indicated domain; for operators, $\nrm{\cdot}$ and $\nrm{\cdot}_1$ denote the operator and trace norms. Finite-dimensional vectors use the Euclidean norm, and sequence norms carry the subscript $\ell^p$. The energy inner product on $\mathcal H=H^1(\R)$ is specified in \cref{app:foundation}.

\section{Finite realizations and spectral data}

\subsection{The fibre and its energy}

Throughout, $f\sim g$ means $f/g\to1$ in the stated limit. Positive constants $c,C$ may depend on the fixed polynomial parameters and change between estimates.

Let $p$ be monic of degree $n$, with zeros $1>\zeta_1>\cdots>\zeta_n>-1$. For any scalar normalization of a polynomial, $\cJ(p)$ denotes the fibre of its monic normalization. Given $J\in\cJ(p)$, number its vertices by $0,\ldots,n-1$, put $b_i=J_{ii}$, and write the entry joining $i-1$ and $i$ as $a_i>0$. The notation $J_{a:b}$ denotes the principal block on vertices $a,\ldots,b$. Set
\begin{equation}
 d_i(J)=[(I-J)^{-1}]_{ii},\qquad
 \E(J)=\sum_i d_i(J)^2,\qquad
 \sum_i d_i(J)=\frac{p'(1)}{p(1)}.
\end{equation}
For $p=P_n^{(\alpha,\beta)}$ the fibre depends on the zeros, and its trace is
$n(n+\alpha+\beta+1)/(2(\alpha+1))$. Cauchy--Schwarz gives
\[
 n^{-3}\E(J)\ge\left[\frac{n+\alpha+\beta+1}{2n(\alpha+1)}\right]^2.
\]
This fixes the scale used in the definition of $e(\alpha,\beta)$.

Let $\varphi_k$ be the unit eigenvector at $\zeta_k$, chosen positive at vertex zero. The weights $w_k=\varphi_k(0)^2$ are positive and sum to one. Conversely, orthonormalizing $1,x,\ldots,x^{n-1}$ for $\sum_k w_k\delta_{\zeta_k}$ gives polynomials $\psi_i$ and their Jacobi recurrence \cite{dBG}. Hence
\begin{equation}\label{eq:eigvec}
 \varphi_k(i)=\sqrt{w_k}\,\psi_i(\zeta_k),\qquad
 d_i(J)=\sum_k\frac{w_k\psi_i(\zeta_k)^2}{1-\zeta_k}.
\end{equation}
The open probability simplex therefore parametrizes the fibre analytically.

Let $\mathsf R$ be the reversal permutation matrix and let $w^\vee$ be the weights of $\mathsf RJ\mathsf R$. A matrix is persymmetric when $\mathsf RJ\mathsf R=J$. The endpoint cofactor and its residue at $\zeta_k$ give
\[
 \varphi_k(0)\varphi_k(n-1)p'(\zeta_k)=\prod_{i=1}^{n-1}a_i,
 \qquad w_k^\vee=\frac{[w_kp'(\zeta_k)^2]^{-1}}
 {\sum_\ell[w_\ell p'(\zeta_\ell)^2]^{-1}}.
\]
Solving $w=w^\vee$ gives
\begin{equation}\label{lem:persym}
 w_k^*=\frac{|p'(\zeta_k)|^{-1}}{\sum_\ell|p'(\zeta_\ell)|^{-1}}.
\end{equation}
These are the weights of the unique persymmetric matrix $J_p^*$ \cite{dBS86,GTVZ17}. Its eigenvectors have reflection sign $(-1)^{k-1}$ and its resolvent diagonal is symmetric.

\begin{proposition}[Criticality]
The persymmetric matrix $J_p^*$ is a critical point of $\E$ on $\cJ(p)$.
\end{proposition}
\begin{proof}
Use the analytic coordinates $\xi_k=\log(w_k|p'(\zeta_k)|)$ modulo constants. Reversal sends $\xi$ to $-\xi$ and preserves $\E$. The resulting even analytic function has zero differential at $\xi=0$, which corresponds to \eqref{lem:persym}.
\end{proof}

Reciprocal-derivative weights also enter the construction of spin chains with perfect state transfer \cite{VZ12}. For a fixed spectrum, Ueoka et al.\ prove that the persymmetric matrix uniquely maximizes the endpoint transfer ceiling, and establish a universal normalized endpoint transfer profile and quantitative near-persymmetry \cite[Theorems~4.4, 5.1 and~8.2]{Ue26}.

\subsection{The first variation of a realization}

Fix an interval $\mathcal I=[\lambda_-,\lambda_+]$, $0\le\lambda_-<\lambda_+$, and put
$\Delta=\lambda_+-\lambda_-$, $z_0=(\lambda_++\lambda_-)/\Delta$, $z(\lambda)=z_0-2\lambda/\Delta$.
Fix $J\in\cJ(p)$ and let $P_i(z)=\det(zI-J_{0:i-1})$, $P_0=1$, so $P_n=p$. Interlacing gives $P_i(z_0)>0$. Let $H_i$ be self-adjoint operators on one finite-dimensional real Hilbert space. Given an initial error $e_0$ and $e_{-1}=e_0$, the normalized recurrence defines a realization $\mathcal R$:
\begin{equation}\label{eq:jacobi-iteration}
 \begin{gathered}
 e_{i+1}=(I-\omega_{i+1}H_i)e_i+\delta_i(e_i-e_{i-1}),\\
 \omega_{i+1}=\frac{2P_i(z_0)}{\Delta P_{i+1}(z_0)},\qquad
 \delta_i=\frac{a_i^2P_{i-1}(z_0)}{P_{i+1}(z_0)}\ (1\le i<n),\qquad \delta_0=0.
 \end{gathered}
\end{equation}
For a fixed $H$, its terminal error is $p(z(H))e_0/p(z_0)$. Each step uses one application of $H_i$.

Let $H$ be a fixed self-adjoint operator on this space, with spectrum in $\mathcal I$ and spectral projections $\Pi_\lambda$. Write $\mathcal T_{\mathcal R}(H_0,\ldots,H_{n-1})e_0=e_n$. For self-adjoint perturbations $Z_i$, its first variation at $H_i=H+\varepsilon Z_i$ has the spectral form
\[
 \left.\partial_\varepsilon\mathcal T_{\mathcal R}(H+\varepsilon Z_0,\ldots,H+\varepsilon Z_{n-1})\right|_{\varepsilon=0}
 =\sum_{i=0}^{n-1}\sum_{\lambda,\lambda'\in\spec H}K_i(\lambda,\lambda')\Pi_\lambda Z_i\Pi_{\lambda'}.
\]
The sharp spectral-pair gain is
\begin{equation}\label{eq:gain-definition}
 A_n(\mathcal R;\mathcal I)=\sup_{\lambda,\lambda'\in \mathcal I}
       \left(\sum_{i=0}^{n-1}|K_i(\lambda,\lambda')|^2\right)^{1/2}.
\end{equation}
For a fixed spectral pair it is the operator norm from the stepwise $\ell_2$ amplitudes to the terminal first variation.

\begin{proposition}[Energy and endpoint drift]\label{prop:drift-energy}
Write $Q_{n-i-1}^{(i+1)}(z)=\det(zI-J_{i+1:n-1})$, with empty determinant one. The realization \eqref{eq:jacobi-iteration} has
\begin{equation}\label{eq:general-drift-kernel}
 K_i(\lambda,\lambda')=-\frac{2P_i(z(\lambda'))Q_{n-i-1}^{(i+1)}(z(\lambda))}{\Delta p(z_0)},
 \qquad
 K_i(\lambda_-,\lambda_-)=-\frac{2p(1)}{\Delta p(z_0)}d_i(J).
\end{equation}
Consequently $A_n(\mathcal R;\mathcal I)\ge2p(1)\sqrt{\E(J)}/[\Delta p(z_0)]$.
For $\mathcal I=[0,1]$ the endpoint gain is exactly $2\sqrt{\E(J)}$.
\end{proposition}
\begin{proof}
Set $v_i=P_i(z_0)e_i$. The recurrence becomes
$v_{i+1}=(z(H_i)-b_i)v_i-a_i^2v_{i-1}$, with the last term zero at $i=0$. An insertion at step $i$ is $-2Z_iv_i/\Delta$ and propagates to the terminal vertex by the right-tail polynomial. This gives the first identity. The diagonal cofactor formula
$P_i(1)Q_{n-i-1}^{(i+1)}(1)/p(1)=[(I-J)^{-1}]_{ii}$ gives the second. Taking the $\ell_2$ norm proves the gain assertions. This is the Jacobi kernel factorization of Part~I \cite{PT1} in the present normalization.
\end{proof}
Write $T_n$ and $U_n$ for the Chebyshev polynomials of the first and second kinds. For the Chebyshev minimax residual, $p(1)/p(z_0)=\epsilon_n^*=1/T_n(z_0)$, so this recovers the gain bound motivating the energy problem. Induced-gain design for momentum methods \cite{Gur25} and gain-margin optimization \cite{WCJG26} likewise quantify the propagation of oracle errors; here \eqref{eq:general-drift-kernel} expresses that propagation in finite inverse-spectral coordinates.

\subsection{The limiting spectral data}

For $\alpha,\beta>-1$ write $J_n^*=J_{P_n^{(\alpha,\beta)}}^*$, with eigenvectors $\varphi_{n,k}^*$ and weights $w_{n,k}^*$ as in \eqref{eq:eigvec}. Let $j_{\alpha,k}$ be the ordered positive zeros of the Bessel function $J_\alpha$ of the first kind, and put
\begin{equation}\label{eq:edgedata}
 \sigma_k=j_{\alpha,k}^2,\quad \Lambda_k=\frac2{\sigma_k},\quad
 M=\sum_k\Lambda_k=\frac1{2(\alpha+1)},\quad
 T_2=\sum_k\Lambda_k^2=\frac1{4(\alpha+1)^2(\alpha+2)}.
\end{equation}
The Bessel product gives the two sums by comparison of its first coefficients \cite[10.21.15]{DLMF}. Write $\Lambda^{(\alpha)}=(\Lambda_k)_{k\ge1}$, suppressing the superscript when the index is fixed. Define $\nu=\alpha/2+1/4$ and $\gamma=\nu+1/2>1/4$. The labels $T,P,U$ refer to the indices $-1/2,0,1/2$, respectively.

\begin{lemma}[Spectral convergence]\label{lem:spec}
For the decreasing zeros $\zeta_{n,k}$ of $P_n^{(\alpha,\beta)}$, pad
$r_{n,k}=[n^2(1-\zeta_{n,k})]^{-1}$ by zero for $k>n$. Then
$\|r_n-\Lambda\|_{\ell^1}\to0$.
\end{lemma}
\begin{proof}
Mehler--Heine \cite[18.11.5]{DLMF} gives $2n^2(1-\zeta_{n,k})\to\sigma_k$ for each fixed $k$. The trace formula gives $\sum_k r_{n,k}\to M=\sum_k\Lambda_k$. Given $K$, the tail sum is the total sum minus its first $K$ terms. Letting first $n\to\infty$ and then $K\to\infty$ proves the claim by positivity.
\end{proof}
We next determine how the canonical weights distribute these levels over the vertices. This requires uniform estimates for both ends $\zeta\approx1$ and $\zeta\approx-1$; \cref{prop:zeros} supplies them.

\section{The canonical continuum limit}

\subsection{The limiting operator and energy}

The canonical limit identifies the upper construction in the comparison of costs. Its profile is a Green diagonal on $(0,1)$, with endpoint exponent $\gamma=\alpha/2+3/4$ and reflection about the midpoint.

\begin{theorem}[Canonical limit]\label{thm:A}
For every $\alpha>-1$ the construction in \cref{sec:contpersym} gives a positive self-adjoint realization $H_\alpha=-\partial_t^2+V_\alpha$ on $L^2((0,1),\dd t)$, with continuous Green kernel $G_\alpha$. Its potential is locally continuous and satisfies $V_\alpha(1-t)=V_\alpha(t)$. Its simple eigenvalues are $\sigma_k=j_{\alpha,k}^2$ and its normalized eigenfunctions have parity $(-1)^{k-1}$ about $1/2$. Set
\begin{equation}
 \phi_\alpha(t)=2G_\alpha(t,t),\qquad c_*(\alpha)=\int_0^1\phi_\alpha(t)^2\dd t.
\end{equation}
For every $\beta>-1$,
\[
 \sup_{n,i}\frac{d_i(J_n^*)}{n}<\infty,\qquad
 \frac{d_{\lfloor nt\rfloor}(J_n^*)}{n}\longrightarrow\phi_\alpha(t)\quad(0<t<1),\qquad
 n^{-3}\E(J_n^*)\longrightarrow c_*(\alpha).
\]
For $0<t,s<1$, $n^{-1}(I-J_n^*)^{-1}_{\lfloor nt\rfloor,\lfloor ns\rfloor}\to2G_\alpha(t,s)$.
With $(\mathcal U_nv)(t)=\sqrt n\,v_{\lfloor nt\rfloor}$, $0\le t<1$,
\[
 \mathcal U_n\bigl(2n^2(I-J_n^*)\bigr)^{-1}\mathcal U_n^*
 \longrightarrow H_\alpha^{-1}\quad\text{in trace norm}.
\]
In particular $e(\alpha,\beta)\le c_*(\alpha)$.
\end{theorem}
\prooflink{thm:A}

At $\alpha=1/2$, $H_\alpha$ is the Dirichlet Laplacian and
$\phi_{1/2}(t)=2t(1-t)$, $c_*(1/2)=2/15$.
At $\alpha=-1/2$ it gives the sine--Jacobi constant of Parts~I--II. These two explicit reference points will reappear in the lower-bound comparison. The resolvent limit, its shifts and the corresponding left-edge limits are proved in \cref{app:gram,app:shift}.

\subsection{Reconstruction from spectral data}\label{sec:contpersym}

We reconstruct the limit by perturbing an explicit reference spectrum. Its quadratic part matches the Bessel levels through their constant term, so the remaining displacements are summable. The reference potential and eigenvalues are
\begin{equation}
 q_\alpha(t)=\frac{\gamma(\gamma-1)\pi^2}{\sin^2\pi t}+c_0,
 \qquad \tau_k=\pi^2(k-1+\gamma)^2+c_0,
 \qquad c_0=\tfrac14-\alpha^2.
\end{equation}
Let $f(t,z)$ solve $-f''+q_\alpha f=zf$ with
$f(t,z)=t^\gamma(1+O(t^2))$ as $t\downarrow0$.
The order of this remainder specifies the branch for every $\gamma>1/4$.
Its Gegenbauer eigenfunctions and their norms are given in \cref{lem:f}. McMahon's expansion gives $\sigma_k-\tau_k=O(k^{-2})$. Define
\begin{equation}\label{eq:norming-main}
 c_k^0=\left(\int_0^1f(t,\tau_k)^2\dd t\right)^{-1},\qquad
 c_k=c_k^0\prod_{\ell\ne k}\frac{\tau_\ell-\tau_k}{\sigma_\ell-\sigma_k}.
\end{equation}
The product converges absolutely and supplies the continuum form of the finite reciprocal-derivative weights \eqref{lem:persym}.

On $0\le s,t\le1/2$ form
\begin{equation}\label{eq:F}
 F_\alpha(t,s)=\sum_{k\ge1}
 \bigl[c_k f(t,\sigma_k)f(s,\sigma_k)-c_k^0f(t,\tau_k)f(s,\tau_k)\bigr].
\end{equation}
Let $\cF_{\alpha,t}$ be the integral operator with kernel $F_\alpha$ on $L^2((0,t),\dd s)$, $0<t\le1/2$. Write $F=F_\alpha$ and $K=K_\alpha$. The transformed waves are determined by
\begin{equation}\label{eq:GL}
 \begin{aligned}
 K(t,s)+F(t,s)+\int_0^tK(t,r)F(r,s)\dd r&=0,\qquad 0\le s\le t\le1/2,\\
 u(t,z)&=f(t,z)+\int_0^tK(t,s)f(s,z)\dd s.
 \end{aligned}
\end{equation}
\Cref{prop:reconstruction} proves a uniform positive lower bound for $I+\cF_{\alpha,t}$, continuity of $K$, and the differential equation
\[
 -u''+V_\alpha u=zu,\qquad
 V_\alpha=q_\alpha+2\frac{\dd}{\dd t}K(t,t),\qquad
 u(t,z)=t^\gamma(1+O(t)).
\]
The odd levels $\sigma_k$ satisfy $u_t(1/2,\sigma_k)=0$ and the even ones satisfy $u(1/2,\sigma_k)=0$. Thus
\[
 \upsilon_k(t)=\sqrt{c_k}\,u(t,\sigma_k)\quad(0\le t\le1/2),\qquad
 \upsilon_k(1-t)=(-1)^{k-1}\upsilon_k(t)
\]
form a complete orthonormal basis. The spectral operator $H_\alpha\upsilon_k=\sigma_k\upsilon_k$ and the uniformly convergent series
$G_\alpha(t,s)=\sum_k\upsilon_k(t)\upsilon_k(s)/\sigma_k$ define the operator and kernel of \cref{thm:A}.

The reference data are the classical P\"oschl--Teller formulas \cite{PTel33,DLMF}. Kostenko--Sakhnovich--Teschl \cite[Theorem~2.8]{KST10} treat perturbed spherical operators with exponent $\ell+1$, $\ell\ge-1/2$; Guliyev \cite[Theorem~4.8]{Gul23} proves symmetric inverse uniqueness in a class of integer inverse-square indices. Here the normalized branch runs through all $\gamma>1/4$, including the nonprincipal branch when $\gamma<1/2$. For the explicit data \eqref{eq:norming-main}, \cref{app:closure} proves existence and identifies the two midpoint conditions. Their reflection decomposition is the half-interval structure studied in \cite[\S4]{FS25}.

\subsection{The finite reconstruction}

Let $J_n^0$ be the explicit comparison chain of \cref{app:aw}, with corrected cosine nodes $\eta_{n,k}$, canonical weights $w_{n,k}^0$, and orthonormal polynomials $\psi_{n,i}^0$. Its nodes approximate the Jacobi zeros to second order at both ends of the interval. The chain belongs to the $q$-Racah/para-Racah families \cite{SZ96,LVZ18}; its trigonometric form gives positivity when the Askey--Wilson base $q$ lies on the unit circle and uniform control as both degree and lattice size grow. Beaulac et al.\ \cite[Proposition~1]{BCL26a} obtain P\"oschl--Teller limits from a quadratic para-Racah bi-lattice, with convergence of each fixed spectral mode away from the singular points. The algebraic families are developed further in \cite{BCL26b}. Here the comparison chain leads to the reference operator; the Gram transformation carries it to the Bessel-zero spectrum and supplies the uniform bounds needed for the energy.

The target canonical measure is positive on $n$ distinct nodes, so its Gram form is positive definite on polynomials of degree less than $n$. Form the Gram matrix
\[
 \delta_{ij}+(B_n)_{ij}=\sum_{k=1}^n w_{n,k}^*
       \psi_{n,i}^0(\zeta_{n,k})\psi_{n,j}^0(\zeta_{n,k}),\qquad t_i=i/n.
\]
Uniformly for $i,j\le n/2$, the two edges give
\[
 n(B_n)_{ij}=F_\alpha(t_i,t_j)+(-1)^{i+j}F_\beta(t_i,t_j)+o(1).
\]
Here $F_\beta$ is \eqref{eq:F} with $\alpha$ replaced by $\beta$. Pairing adjacent coordinates separates their sums and differences into the two continuous positive forms. This gives a lower bound for every principal Gram section. Cholesky orthogonalization then converges to the transformation \eqref{eq:GL}; the alternating terms cancel in each fixed right-edge wave. In particular the normalized waves converge to $\upsilon_k$, while the left index enters the finite endpoint weights through an explicit scalar normalization.

This is the discrete GL--Cholesky reconstruction of Natterer \cite{Nat89}, as presented in \cite[\S2]{DBK05}. A related continuum-limit analysis of discrete inverse spectral reconstruction on optimal grids is given in \cite{BDK05c}; that work assumes a regular, strictly positive coefficient.

In the present scaling, \cref{prop:finite-limit} proves the two-edge convergence, the uniform bound $n|\varphi_{n,k}^*(i)|^2\le C$ for the required half spectrum, and the limit of each fixed wave. The diagonal identity \eqref{eq:eigvec} then splits into finitely many convergent modes and a tail controlled by \cref{lem:spec}. It yields the uniform diagonal bound, its pointwise limit and, by dominated convergence, the energy limit. The same waves give trace-norm convergence of the embedded inverses.

\subsection{The geometry of the profile}\label{sec:trichotomy}

The Green diagonal also determines a density with the full limiting spectrum. A density here is a nonnegative function $m\in L^1(\R,\dd x)$. For a positive compact operator $A$, write $\lambda(A)=(\lambda_k(A))_{k\ge1}$ in nonincreasing order with multiplicities, padded by zeros for finite rank. Let $\sK_m$ be the integral operator on $L^2(\R,\dd x)$ with kernel
\[
 \sK_m(x,y)=\sqrt{m(x)}e^{-|x-y|}\sqrt{m(y)}.
\]
Define its mass, root mass and energy by
\begin{equation}\label{eq:density-functionals}
 \fm(m)=\int_\R m(x)\dd x,\qquad
 \fr(m)=\int_\R\sqrt{m(x)}\dd x,\qquad
 E(m)=\int_\R m(x)^{3/2}\dd x.
\end{equation}
The last two functionals may take the value $+\infty$.

\begin{proposition}[The Green coordinate]\label{prop:green-coordinate}
The kernel $G_\alpha$ is strictly positive on $(0,1)^2$. Define
\begin{equation}\label{eq:xpic}
 \frac{\dd x}{\dd t}=\phi_\alpha(t)^{-1},\qquad x(1/2)=0,\qquad
 m_\alpha(x(t))=\phi_\alpha(t)^2,
\end{equation}
and extend $m_\alpha$ by zero beyond $X_\alpha=x((0,1))$. Then
\begin{equation}\label{eq:green-density-kernel}
 2G_\alpha(t,s)=\sqrt{\phi_\alpha(t)\phi_\alpha(s)}\,
                  e^{-|x(t)-x(s)|}.
\end{equation}
The nonzero eigenvalues of $\sK_{m_\alpha}$, with multiplicities, are $\Lambda_k=2/\sigma_k$. Moreover
\[
 E(m_\alpha)=c_*(\alpha),\qquad \int_\R m_\alpha(x)\dd x=M,\qquad
 \int_\R\sqrt{m_\alpha(x)}\dd x=1.
\]
\end{proposition}
\begin{proof}
The finite resolvents $R=(I-J)^{-1}$ are positive entrywise and obey
$R_{ij}R_{j\ell}=R_{i\ell}R_{jj}$ for $i\le j\le\ell$. The latter identity is the classical single-pair structure of inverses of nonsingular Jacobi matrices \cite[\S3, 4$^\circ$]{GK37}; see also \cite[\S2.1]{Meu92}. Both properties follow directly from principal cofactor formulas. The kernel limit in \cref{prop:finite-limit} gives these properties for $G_\alpha$. For $f\in C_c^\infty(0,1)$, the spectral domain in \cref{prop:reconstruction} gives $H_\alpha f=-f''+V_\alpha f$. Testing $H_\alpha^{-1}H_\alpha f=f$ against its continuous kernel yields $(-\partial_s^2+V_\alpha(s))G_\alpha(t,s)=\delta_t$ distributionally. Thus the kernel is $C^2$ off the diagonal and its first derivative has jump $-1$. A zero diagonal entry would, by positive definiteness and Cauchy--Schwarz, make its kernel row zero and contradict this jump. Thus the diagonal is positive. Continuity and the displayed multiplicative identity then give positivity between arbitrary points by subdivision into neighbouring pairs.

Write $t\wedge s=\min(t,s)$ and $t\vee s=\max(t,s)$. The multiplicative identity gives positive solutions $h_\pm$ of $-h''+V_\alpha h=0$ with factorization
\[
 G_\alpha(t,s)=\frac{h_-(t\wedge s)h_+(t\vee s)}{W},\qquad
 W=h_-'h_+-h_-h_+'>0.
\]
The derivative jump fixes the constant Wronskian. Hence
$x=\tfrac12\log(h_-/h_+)+\text{constant}$ has derivative $1/\phi_\alpha$, which proves \eqref{eq:green-density-kernel}. The unitary map
\[
 (\mathcal Vf)(t)=\phi_\alpha(t)^{-1/2}f(x(t)),\qquad
 \mathcal V:L^2(X_\alpha,\dd x)\longrightarrow L^2((0,1),\dd t),
\]
conjugates the restriction of $\sK_{m_\alpha}$ to $2H_\alpha^{-1}$ by this kernel identity. On $L^2(\R\setminus X_\alpha,\dd x)$ the operator is zero, so its nonzero spectrum is $\Lambda$. Finally the substitution $\dd x=\dd t/\phi_\alpha$ gives the three integrals, using $\int_0^1\phi_\alpha(t)\dd t=2\Tr H_\alpha^{-1}=M$.
\end{proof}

The kernel identity \eqref{eq:green-density-kernel} is the Davies--Harrell representation through the Green diagonal \cite{DH87}, in the form of \cite[Theorem~2.1(3)]{HW13}. The corresponding Green coordinate underlies the packing of Part~II \cite{PT2}. For the canonical limit, it identifies the full spectral data, both scalar normalizations and the $L^{3/2}$ cost with the limiting finite energy. The endpoint exponent determines the geometry. For $\nu<0$ put $c_\alpha=\sum_k\Lambda_kc_k$; \cref{thm:E} proves that this constant is finite and positive.
\begin{center}
\begin{tabular}{lll}
\toprule
Parameter & $\phi_\alpha(t)$ as $t\downarrow0$ & Ends of $m_\alpha$\\
\midrule
$\nu<0$ & $\sim c_\alpha t^{2\gamma}$ & compact support\\
$\nu=0$ & $\sim2t\log(1/t)$ & doubly exponential\\
$\nu>0$ & $\sim t/\nu$ & exponential, rate $2/\nu$\\
\bottomrule
\end{tabular}
\end{center}
\Cref{thm:E} proves positivity and these asymptotics with their precise constants. \Cref{fig:profiles}(a) shows the three regimes. Integrating $1/\phi$ gives a finite endpoint distance when $\nu<0$ and an infinite distance in the other two cases. For $U$, $t=(1+\tanh x)/2$ and $m_U=\tfrac14\sech^4x$, where $\sech x=1/\cosh x$.

\begin{remark}[String form of a density]\label{rem:string}
Let $m\ge0$ be integrable on $\R$. Put $\xi=\tanh x$, so that $\dd x=\dd\xi/(1-\xi^2)$, and
\[
 \rho(\xi)=m(x)\cosh^4x=\frac{m(\operatorname{artanh}\xi)}{(1-\xi^2)^2},
 \qquad
 (\mathcal Uu)(\xi)=\frac{u(x)}{\cosh x}.
\]
Integration by parts gives
\[
 \int_\R(u'w'+uw)\dd x=\int_{-1}^1(\mathcal Uu)'(\mathcal Uw)'\dd\xi,
 \qquad
 \int_\R m\,uw\dd x=\int_{-1}^1\rho\,(\mathcal Uu)(\mathcal Uw)\dd\xi.
\]
The boundary term $[\xi f^2/(1-\xi^2)]_{-1}^{1}$ vanishes for $f\in H^1_0(-1,1)$.
 
Hence $\mathcal U$ is a unitary map from the space $\mathcal H=H^1(\R)$ of \cref{app:foundation} onto $H^1_0(-1,1)$ with the Dirichlet norm. It sends the feature vector $v_x$ to $\sqrt2\cosh x\,g(\tanh x,\cdot)$, where
\[
 g(\xi,\eta)=\tfrac12(1+\xi\wedge\eta)(1-\xi\vee\eta)
\]
is the influence function of a string of unit tension fixed at $\pm1$. By \eqref{eq:features} and \eqref{eq:app-BS}, $\cK_mu=\lambda u$ if and only if $f=\mathcal Uu$ solves
\[
 -f''=z\rho f\ \text{ on }(-1,1),\qquad f(\pm1)=0,\qquad z=2/\lambda,
\]
in the weak sense. Thus $\lambda(\sK_m)=(2/z_k)_{k\ge1}$, where the $z_k$ are the Dirichlet eigenvalues of the string with density $\rho$. This is the Liouville transformation of \cite[(3.1)--(3.5)]{BSS00}, written there for the Lax operator of the Camassa--Holm equation; see also \cite[(36)]{EKT16}.
 
In string form the three integrals of \cref{prop:green-coordinate} read
\[
 \begin{aligned}
 \fm(m)&=\int_\R m(x)\dd x=\int_{-1}^1(1-\xi^2)\rho(\xi)\dd\xi=2\sum_k\frac1{z_k},\\
 \fr(m)&=\int_\R\sqrt{m(x)}\dd x=\int_{-1}^1\sqrt{\rho(\xi)}\dd\xi,\\
 E(m)&=\int_{-1}^1(1-\xi^2)^2\rho(\xi)^{3/2}\dd\xi.
 \end{aligned}
\]
The first is twice the trace of the string's Green operator.
 
The weights are those of the hyperbolic length $\dd\xi/(1-\xi^2)$. A translation $x\mapsto x-x_0$ acts on strings by the automorphism $\xi\mapsto(\xi-\tanh x_0)/(1-\xi\tanh x_0)$ of $(-1,1)$ and multiplies the density by $(\cosh x_0-\xi\sinh x_0)^{-4}$.
 
The root mass controls the density of the spectrum. For a string with $\int_{-1}^1\rho(\xi)\dd\xi<\infty$, Krein's theorem \cite{Kr52} (see \cite[Theorem~5.1]{KP13}) gives
\[
 \sqrt{z_k}\sim\frac{\pi k}{\displaystyle\int_{-1}^1\sqrt{\rho(\xi)}\dd\xi},
\]
and the leading term does not depend on the boundary conditions.
 
By reflection symmetry, $\rho_\alpha$ is even: a symmetric string with Dirichlet spectrum $\{\sigma_k\}$. The finite persymmetric uniqueness \cite{Hoch67,dBS86} thus has its continuum counterpart in Krein's determination of the density of a symmetric string from its spectrum \cite{Kr51b}.
 
Normalize $h_+(t)=h_-(1-t)$ in the proof of \cref{prop:green-coordinate}. The composite map $t\mapsto\xi$ is then the classical Liouville transformation of $H_\alpha$ by its positive null solution $Y=(h_-+h_+)/\sqrt{2W}$:
\[
 \frac{\dd\xi}{\dd t}=Y^{-2},\qquad \rho_\alpha(\xi(t))=Y(t)^4 .
\]
At $U$ we have $h_-=t$, $h_+=1-t$ and $W=1$. Hence $\xi=2t-1$ and $\rho_{1/2}\equiv\frac14$: the canonical $U$ density is the homogeneous string, with $z_k=\pi^2k^2$.
 
\Cref{thm:E} gives the three endpoint geometries in string form.
\begin{itemize}
 \item For $\nu<0$, $\rho_\alpha$ is supported on $[-\tanh x_*,\tanh x_*]$, so the string carries massless end segments.
 \item For $\nu=0$, $\rho_\alpha$ vanishes at $\pm1$ faster than any power of $1-\xi^2$.
 \item For $\nu>0$, $\rho_\alpha(\xi)=(1-\xi^2)^{1/\nu-2+o(1)}$ as $\xi\to\pm1$. The exponent vanishes only at $\nu=1/2$, the index $U$. The density tends to zero at the ends for $\nu<1/2$, blows up for $\nu>1/2$, and is not integrable for $\nu>1$.
\end{itemize}
For $\nu<1$ the density $\rho_\alpha$ is integrable, and the asymptotics $j_{\alpha,k}\sim\pi k$ agree with $\int_\R\sqrt{m_\alpha(x)}\dd x=1$.
\end{remark}

The density formulation also accommodates arbitrary finite realizations. We now use it to compare $c_*(\alpha)$ with lower bounds valid over the entire fibre.
\section{Variational lower bounds}\label{sec:lower}

\subsection{Packing all realizations}

The upper construction has produced one density with spectrum $\Lambda$. To test every realization, we first pack finite matrices into that same class. A second-moment inequality then gives the exact answer at $U$. The subsequent refinements retain the mass and successively more spectral or geometric information; their active constraints explain the changes of extremal profile.

Here symmetric decreasing means even and nonincreasing on $[0,\infty)$. For nonnegative integrable $m$, $\sK_m$ is positive trace class with trace $\int_\R m(x)\dd x$. Packing preserves the cost $E(m)$ and the two scalar normalizations.

\begin{proposition}[Resolvent packing]\label{prop:packing}
Let $J$ be a symmetric tridiagonal $n\times n$ matrix with nonnegative off-diagonal entries and $\nrm J<1$. There are intervals $I_i$ of lengths $1/d_i$, with disjoint interiors, such that $m_J=\sum_i n^{-2}d_i^2\ind_{I_i}$ satisfies, with $M_J=n^{-2}\Tr(I-J)^{-1}$,
\begin{equation}\label{eq:packing-main}
 \begin{gathered}
 \int_\R m_J(x)\dd x=M_J,\qquad \int_\R\sqrt{m_J(x)}\dd x=1,\qquad
 E(m_J)=n^{-3}\E(J),\\
 \nrm{\lambda(\sK_{m_J})-\lambda(n^{-2}(I-J)^{-1})}_{\ell^1}
 \le4\sqrt{M_J/n}.
 \end{gathered}
\end{equation}
\end{proposition}
\prooflink{prop:packing}
The Green-coordinate construction is from Part~II \cite[Theorem~2.3]{PT2}. The proof in \cref{app:packing} gives its feature-operator representation, checks the mass normalization and treats reducible matrices. The three exact integrals identify the cost and normalizations; the spectral error transfers continuous lower bounds to the finite fibres. \Cref{lem:spectral-compactness} supplies the compactness needed for this transfer.

\begin{remark}[Finite realizations as beaded strings]\label{rem:beads}
Let $J$ have positive links and $\nrm J<1$. The single-pair identity $R_{ij}R_{jl}=R_{il}R_{jj}$ ($i\le j\le l$) from the proof of \cref{prop:green-coordinate} gives points $X_0<\cdots<X_{n-1}$, unique up to a common translation, with
\[
 R_{ij}=\sqrt{d_id_j}\,e^{-|X_i-X_j|},\qquad
 X_{i+1}-X_i=\tfrac12\log\frac{d_id_{i+1}}{R_{i,i+1}^2}>0 .
\]
Put $\mu=n^{-2}\sum_id_i\delta_{X_i}$. For a finite positive measure $\mu$, define the operator on $L^2(\R,\mu)$ by
\[
 (\sK_\mu f)(x)=\int_\R e^{-|x-y|}f(y)\dd\mu(y).
\]
If $\dd\mu(x)=m(x)\dd x$, multiplication by $\sqrt m$ identifies this operator with the restriction of $\sK_m$ to $L^2(\{m>0\},\dd x)$. For the atomic measure above, its matrix in the normalized indicator basis is $n^{-2}R$, the Gram matrix $\Gamma$ in the proof of \cref{prop:packing}. The spectral and weighted-mass identities of \cref{rem:string} extend to finite measures; the root mass and energy in \eqref{eq:density-functionals} are defined for densities.
 
Put $\xi_i=\tanh X_i$ and $w_i=n^{-2}d_i\cosh^2X_i$. Let $\ell_0,\ldots,\ell_n$ be the lengths of the segments into which the points $\xi_i$ divide $[-1,1]$. Let $\mathsf W=\diag(w_i)$, and let $\mathsf S$ be the stiffness matrix of the string, with $\mathsf S_{ii}=\ell_i^{-1}+\ell_{i+1}^{-1}$ and $\mathsf S_{i,i+1}=-\ell_{i+1}^{-1}$. Then
\[
 2n^2(I-J)=\mathsf W^{-1/2}\mathsf S\,\mathsf W^{-1/2},
 \qquad
 n^{-2}d_i=2w_i\,g(\xi_i,\xi_i)=w_i(1-\xi_i^2).
\]
Indeed $\mathsf S^{-1}=[g(\xi_i,\xi_j)]$ and $2\sqrt{w_iw_j}\,g(\xi_i,\xi_j)=n^{-2}\sqrt{d_id_j}\,e^{-|X_i-X_j|}$.
 
Every $J\in\cJ(p)$ is therefore the normalized matrix of small oscillations of a string carrying $n$ beads, fixed at $\pm1$, with squared frequencies $2n^2(1-\zeta_k)$. The common translation of the $X_i$ is the automorphism freedom of \cref{rem:string}.
 
This is the setting of Gantmacher and Krein \cite{GK02}. Their single-pair matrices \cite[\S3, 4$^\circ$]{GK37} are, up to diagonal congruence, exactly the influence matrices $[g(\xi_i,\xi_j)]$. Krein solved the inverse problem for the beaded string by Stieltjes continued fractions \cite[Supplement~II]{GK02}. The same discrete strings encode multipeakons; in particular $n^{-2}d_i=w_i(1-\xi_i^2)$ is the peakon mass $m_j=g_j(1-y_j^2)$ of \cite[(3.7)--(3.8)]{BSS00}.
 
In this language,
\[
 n^{-4}\E(J)=\sum_iw_i^2(1-\xi_i^2)^2
\]
is the sum of the squared mass-weighted self-influence coefficients of the beads. Packing spreads each bead over an interval $I_i$ of optical length $\int_{I_i}\sqrt{m_J(x)}\dd x=1/n$.
 
The passage from discrete to continuous strings is classically controlled by the continuity of Krein's correspondence between strings and spectral functions \cite{Kas75}; see \cite[\S1]{Kot13}. Here \cref{lem:spectral-compactness} plays this role, and it also keeps the energy in view.
\end{remark}

\subsection{The second moment and the exact point}\label{sec:second}

The simplest spectral constraint is the second moment
$G_2(m)=\Tr\sK_m^2=\iint_{\R^2}m(x)e^{-2|x-y|}m(y)\dd x\dd y$.
Duality for the sharp one-dimensional Gagliardo--Nirenberg inequality \cite[\S3.1]{DELL14}, with the normalization checked in \cref{app:moment}, gives
\begin{equation}\label{eq:sharp-main}
 G_2(m)\le C_{\rm GN}E(m)^{4/3},\qquad C_{\rm GN}=(5/18)^{1/3}.
\end{equation}
It holds for every nonnegative $m\in L^{3/2}(\R,\dd x)$, with $G_2$ understood through its double integral when the mass is infinite; the nonzero equality cases are precisely $a\sech^4(x-x_0)$ \cite{DELL14}. Hence $G_2\ge T_2$ implies
\[
 E(m)\ge\ell_Y=(T_2/C_{\rm GN})^{3/4}
 =\frac{2\sqrt3}{5}\left(\frac{5T_2}{2}\right)^{3/4}.
\]
At $U$, the canonical $m_U=\tfrac14\sech^4x$ is such an extremizer and
\begin{equation}
 \int_\R m_U(x)\dd x=\tfrac13,\qquad \int_\R\sqrt{m_U(x)}\dd x=1,\qquad E(m_U)=\tfrac2{15},
 \qquad \spec_{\ne0}\sK_{m_U}=\{2/(\pi^2k^2):k\ge1\}.
\end{equation}
The full spectrum follows from \cref{prop:green-coordinate} for the Dirichlet Laplacian. The sharp moment bound and the canonical upper value now give $e(1/2,\beta)=2/15$.

In the string form of \cref{rem:string}, $m_U$ is the homogeneous string $\rho\equiv\frac14$ on $(-1,1)$. The equality cases $a\sech^4(x-x_0)$ of \eqref{eq:sharp-main} are the strings $a(\cosh x_0-\xi\sinh x_0)^{-4}$, the images of a homogeneous string under the automorphisms of $(-1,1)$. The four identities above become
\[
 \tfrac14\int_{-1}^1(1-\xi^2)\dd\xi=\tfrac13,\qquad
 \tfrac12\int_{-1}^1 1\dd\xi=1,\qquad
 \tfrac18\int_{-1}^1(1-\xi^2)^2\dd\xi=\tfrac2{15},\qquad
 z_k=\pi^2k^2 .
\]
Exactness at $U$ thus says that the uniform string is both the canonical limit and a free second-moment extremal.

The free extremizer uses mass $\sqrt{5T_2/2}$. Comparing this with the available $M$ gives the first refinement:

\[
 L_2=\inf\{E(m):0\le m\in L^1(\R,\dd x),\ \fm(m)\le M,\ G_2(m)\ge T_2\}.
\]
Thus $L_2=\ell_Y$ when $5T_2\le2M^2$, and $L_2>\ell_Y$ otherwise. For Jacobi data $T_2/M^2=1/(\alpha+2)$, so the threshold is $\alpha=1/2$.

One ODE explains the constrained regime. Write $r_+=\max(r,0)$ and let $*$ denote convolution on $\R$. For a symmetric decreasing minimizer put $W=e^{-2|\cdot|}*m$. The normal optimality conditions proved in \cref{app:foundation} give, for $\kappa>0,h\ge0$,
\begin{equation}\label{eq:moment-mechanism}
 -W''+4W=4m,\qquad \sqrt m=\kappa(W-h)_+,\qquad
 (W')^2=4W^2-\frac{8\kappa^2}{3}(W-h)_+^3.
\end{equation}
The last identity follows by multiplying the ODE by $W'$ and using the decay at infinity. For $h=0$ it gives the free profile; for $h>0$ the clipping produces compact support. Eliminating the scale leaves one shape parameter. Its strict monotonicity, unique determination and elliptic quadratures are established in \cref{thm:elliptic}. Thus the change of branch reflects an active mass constraint.

\subsection{Retaining more spectral information}

The second moment summarizes the spectrum by one sum. Testing its largest level gives a complementary bound:

\[
 L_\infty=\inf\{E(m):0\le m\in L^1(\R,\dd x),\ \fm(m)\le M,\ \lambda_1(\sK_m)\ge\Lambda_1\}.
\]
The free Keller value is $\frac\pi4(3\Lambda_1/2)^{3/2}$, with mass $3\Lambda_1/2$ \cite{Kel,LT}. It is the answer when this mass is at most $M$; the complementary regime has a uniquely parametrized compact branch. For its positive ground state $u$, normalized by $u(0)=1$, the optimality condition has the form $\sqrt m=a(u^2-v_e)_+$, with $a>0$ and $0\le v_e<1$, and the equation $-u''+u=(2/\Lambda_1)mu$ gives a first integral. \Cref{thm:keller} gives the full formula and its mass threshold.

The regularized determinant $\det_2$ \cite{Sim} combines all trace moments of orders $p\ge2$ into one constraint. For $s>0$ define
\begin{equation}
 \hD_s(m)=-\log\det{}_2(I-s\sK_m)
   =\sum_{p\ge2}\frac{s^p}{p}\Tr\sK_m^p
 \quad\text{if }s\lambda_1(\sK_m)<1,
\end{equation}
and set it to $+\infty$ otherwise. For $0<s<1/\Lambda_1$, let $\hd_s=\sum_k[-\log(1-s\Lambda_k)-s\Lambda_k]$ and
\[
 B_s=\inf\{E(m):0\le m\in L^1(\R,\dd x),\ \fm(m)\le M,\ \hD_s(m)\ge\hd_s\}.
\]
Small $s$ emphasizes the second moment; near the first pole the top level becomes more influential. Let $g(x)=(-\partial_x^2+1-2sm)^{-1}(x,x)$. Below the first pole, the gradient is $s(2g-1)$, so the optimality condition is $m=A(g-g_e)_+^2$, with $A>0$ and $g_e\ge1/2$. The Gel'fand--Dikii identity \cite{GD} closes an ODE for $g$ after this substitution. Its first integral gives the profiles and quadratures in \cref{app:fredholm}. There is at most one candidate in each of the free and mass-active branches; at least one admissible branch exists, and the least admissible energy is exactly $B_s$ (\cref{thm:fred}).

The Birman--Schwinger equation \cite{Bir61,Sch61} is $-u''+u=(2/\lambda)mu$: fixing $\lambda$ fixes the negative energy $-1$ of $-\partial_x^2-2m/\lambda$. This places the problems alongside optimization of finite collections of negative Schr\"odinger eigenvalues \cite{FGL21,FGL25} and the one-dimensional eigenvalue-sum inequality of \cite{RS26}. The free interpolation and Keller values supply our starting profiles. The classifications in \cref{app:moment,app:fredholm} determine their deformation under a mass cap and the constrained determinant value at this fixed energy.

\begin{remark}[Isoperimetric problems for strings]\label{rem:iso}
By \cref{rem:string}, $L_\infty$ asks for the string on $(-1,1)$ that minimizes $\int_{-1}^1(1-\xi^2)^2\rho(\xi)^{3/2}\dd\xi$, subject to a fundamental Dirichlet eigenvalue $z_1\le2/\Lambda_1$ and a weighted mass $\int_{-1}^1(1-\xi^2)\rho(\xi)\dd\xi\le M$.
 
Extremal problems for the eigenvalues of a string of given length and total mass, with or without a bound on the density, go back to Krein \cite{Kr51a}. For constraints $\int_{-1}^1\rho(\xi)^\gamma\dd\xi=1$ on a fixed interval see Egorov and Kondrat'ev \cite{EK96}. The present problems differ in their hyperbolic weights, which make all functionals invariant under the automorphisms of $(-1,1)$, that is, under translations on the line.
 
The determinant constraint has a similar reading. Let $\varphi_\rho(\xi,z)$ solve $-\varphi''=z\rho\varphi$ with $\varphi(-1)=0$ and $\varphi'(-1)=1$. If $\int_{-1}^1\rho(\xi)\dd\xi<\infty$ (in particular for every compactly supported $m$, such as a packed density), then
\[
 \det(I-s\sK_m)=\prod_k\Bigl(1-\frac{2s}{z_k}\Bigr)=\tfrac12\varphi_\rho(1,2s).
\]
To see this, note four facts about $\varphi_\rho(1,\cdot)$: it is entire of order at most one; its zeros are the simple eigenvalues $z_k$; $\varphi_\rho(1,0)=2$; and its logarithmic derivative at $0$ equals $-\frac12\int_{-1}^1(1-\xi^2)\rho(\xi)\dd\xi=-\sum_k1/z_k$. Hadamard's theorem then gives the product. Thus $B_s$ prescribes the characteristic function of the string at the single point $z=2s$ below its spectrum. For Jacobi data this is the Bessel product used in \cref{app:fredholm}.
\end{remark}

\subsection{Retaining the root mass}

For $\alpha>1/2$, the free second-moment extremizer uses less than $M$. The mass cap permits diffuse additions carrying the remaining mass at negligible cost. Packing supplies the stronger information $\int_\R m(x)\dd x=M$ and $\int_\R\sqrt{m(x)}\dd x\le1$. After symmetric decreasing rearrangement, the root cap gives $m(x)\le1/(4x^2)$, hence uniform control of the mass in the tails. It preserves exact mass under the limiting passage of \cref{lem:spectral-compactness}. This is the refinement that uses the spatial normalization of the finite chain.

We retain these two constraints through their multipliers. For $b>0$ and $\eta\ge-34/5$ put
\begin{equation}\label{eq:root-main}
 \begin{aligned}
 V_{b,\eta}(T_2)&=\inf_{m\ge0,\ G_2(m)\ge T_2}
 \left\{E(m)+\frac{\eta b}{4}\fm(m)+b^2\fr(m)\right\},\\
 \ell_H(\Lambda;b,\eta)&=V_{b,\eta}(T_2)-\frac{\eta b}{4}M-b^2.
 \end{aligned}
\end{equation}
For a competitor with exact mass and root mass at most one, its value in braces is at most $E(m)+(\eta b/4)M+b^2$. Taking the infimum proves $E(m)\ge\ell_H$. The parameter $\eta b/4$ is the multiplier of the exact mass.

The inner problem is solved by global pointwise minimization of a cubic in $\sqrt m$. A needle variation, corrected to preserve $G_2$, selects this global branch. The minimizing density has finite support and jumps from $\sqrt m=b$ to zero at its boundary. \Cref{app:root} proves this principle and reduces the answer to four quadratures and one strictly monotone scalar equation. The polynomial certificate at the end of \cref{app:root} proves monotonicity throughout the stated parameter range.

\subsection{The resulting lower bound}

The second moment gives exactness at $U$; the other constructions refine the comparison away from that point. Their inputs and outputs are:
\begin{center}
\begin{tabularx}{\linewidth}{lXll}
\toprule
Value & Retained constraints & Type of answer & Formula\\
\midrule
$L_2$ & mass cap, second moment & exact infimum & \cref{app:moment}\\
$L_\infty$ & mass cap, top level & exact infimum & \cref{app:moment}\\
$B_s$ & mass cap, regularized determinant & exact infimum & \cref{app:fredholm}\\
$\ell_H$ & exact mass, root mass, second moment & dual lower bound & \cref{app:root}\\
\bottomrule
\end{tabularx}
\end{center}
Here the data may be any nonincreasing positive summable sequence $\Lambda$: set $M=\sum_k\Lambda_k$, $T_2=\sum_k\Lambda_k^2$, $\Lambda_1=\max_k\Lambda_k$. In particular $T_2<M^2$ and $\Lambda_1<M$. Define
\begin{equation}
 \ell(\Lambda)=\max\left\{L_2,L_\infty,
 \sup_{0<s<1/\Lambda_1}B_s,
 \sup_{b>0,\,\eta\ge-34/5}\ell_H(\Lambda;b,\eta)\right\}.
\end{equation}
For Jacobi data write $\ell(\alpha)=\ell(\Lambda^{(\alpha)})$.

\begin{theorem}[Lower bound for general edge data]\label{thm:transfer-main}
Let $p_n$ have degree $n$ and simple zeros $1>\zeta_{n,1}>\cdots>\zeta_{n,n}>-1$, and suppose $r_{n,k}=[n^2(1-\zeta_{n,k})]^{-1}$, padded by zeros, converges in $\ell^1$ to $\Lambda$. Then
\[
 \liminf_{n\to\infty}\inf_{J\in\cJ(p_n)}n^{-3}\E(J)\ge\ell(\Lambda).
\]
\end{theorem}
\begin{proof}
An infinite lower limit satisfies the claim. For a finite lower limit, take an asymptotically minimizing subsequence and pack it by \eqref{eq:packing-main}. Rearrangement increases the retained spectral functionals and preserves energy and both scalar integrals. The compactness and semicontinuity in \cref{lem:spectral-compactness} pass their constraints to the limit; with bounded root mass they also preserve exact mass. Applying the four bounds and then taking the parameter suprema proves the claim. \Cref{app:packing} records the full transfer, including joint moment constraints and the root-constrained Fredholm class.
\end{proof}
The canonical construction and the transfer now bound the same finite-fibre optimum from above and below.

\section{Comparison and first-order contact}\label{sec:comparison}

\begin{figure}[t]
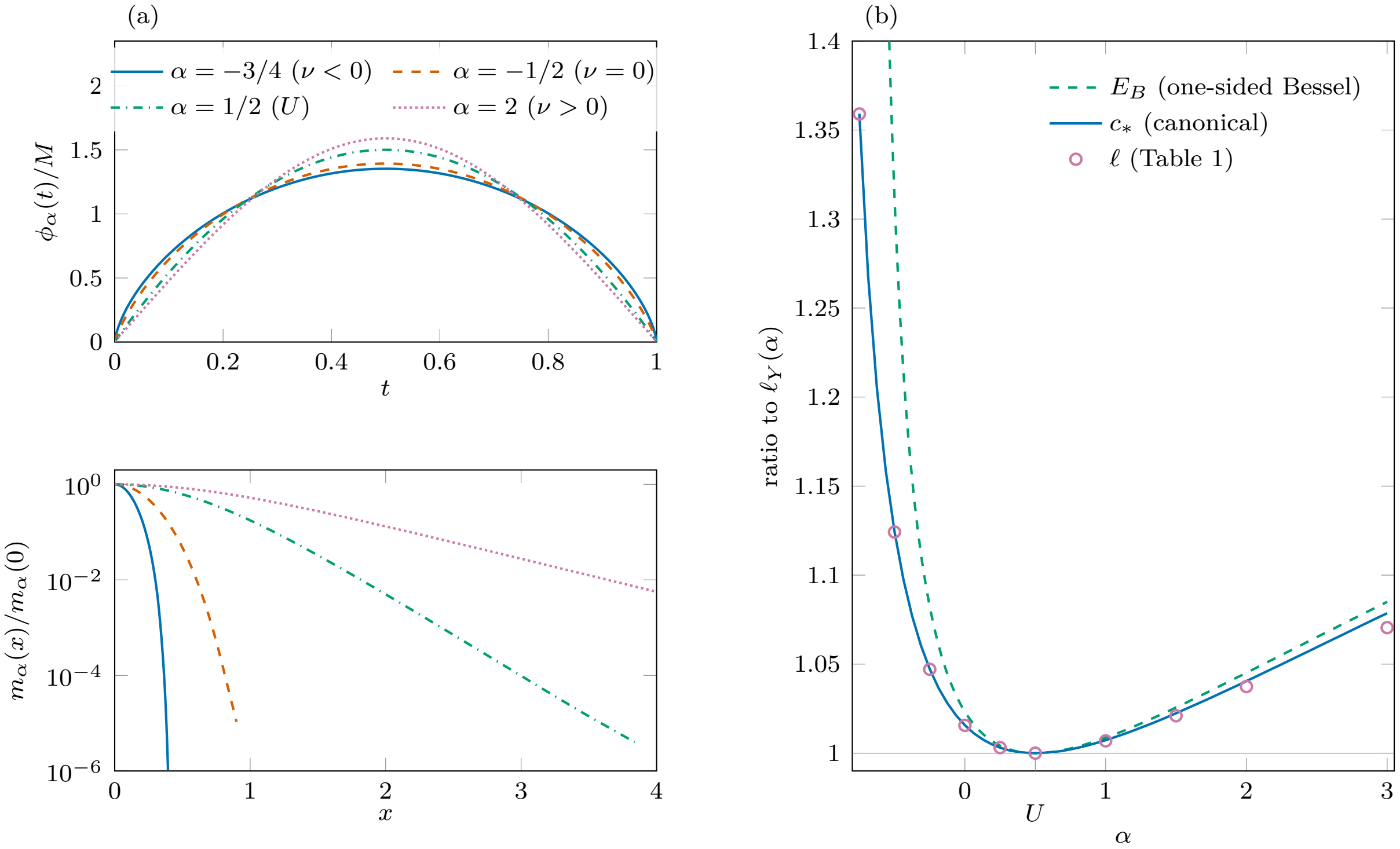

\centering

\caption{(a) Canonical profiles $\phi_\alpha=2G_\alpha(t,t)$, normalized by the mass $M=1/(2(\alpha+1))$ (top), and the corresponding Green-coordinate densities $m_\alpha(x)=\phi_\alpha(t(x))^2$ on $x\ge0$ (bottom, logarithmic scale). The three endpoint regimes of \cref{sec:trichotomy} appear as a compact support with a vertical edge ($\nu<0$), a doubly exponential tail ($\nu=0$) and exponential tails of rate $2/\nu$ ($\nu>0$; $m_U=\tfrac14\sech^4x$). Curves are computed from $d_i(J_n^*)/n$ with $n=2000$, $\beta=\alpha$, where $dx=1/d_i$ and $m=(d_i/n)^2$ on the $i$th cell, as in \cref{prop:packing}; for $\alpha=-3/4$ the last cells are completed by the endpoint asymptotics of \cref{thm:E}, and the support ends at $x_*\approx0.42$. (b) The canonical cost $c_*$, the lower bound $\ell$ and the one-sided Bessel cost $E_B$ of \cref{sec:bessel-geometry}, divided by the explicit second-moment bound $\ell_Y(\alpha)=\frac{2\sqrt3}{5}(5T_2/2)^{3/4}$. All three curves touch $1$ with zero slope at $U$ (\cref{thm:C} and \cref{sec:bessel-geometry}); $E_B=+\infty$ for $\alpha\le-3/4$. The $c_*$ curve is the Richardson extrapolation of $n^{-3}\E(J_n^*)$ from $n=500,1000,2000$ on a grid of step $1/16$; at the ten indices of Table~\ref{tab:main} it agrees with the tabulated $\widetilde c_*$ to within $2.5\cdot10^{-6}$ in relative terms.}\label{fig:profiles}
\end{figure}

\begin{theorem}[Two-sided estimate]\label{thm:B}
For all $\alpha,\beta>-1$,
\[
 \ell(\alpha)\le e(\alpha,\beta)\le c_\ast(\alpha),\qquad
 e(1/2,\beta)=\ell(1/2)=c_\ast(1/2)=2/15.
\]
\end{theorem}
\begin{proof}
Combine \cref{thm:A,thm:transfer-main} and the exact $U$ calculation in \cref{sec:second}.
\end{proof}

The canonical cost $c_\ast$ and universal lower bound $\ell$ bracket the lower limiting optimum $e$. For the exact bounds, $\delta=(c_\ast-\ell)/c_\ast$ gives an asymptotic energy ratio at most $c_\ast/\ell=1/(1-\delta)$. Table~\ref{tab:main} evaluates the comparison numerically. \Cref{fig:profiles}(b) shows $c_*$, $\ell$ and the one-sided Bessel cost $E_B$ of \cref{sec:bessel-geometry} across the family. Its Fredholm and root columns use the selected admissible parameters recorded in \cref{app:table}; their maximum gives the displayed numerical gap. The endpoint identity in \cref{prop:drift-energy} transfers energy bounds to algorithmic sensitivity, and \cref{sec:nesterov} evaluates the full spectral gain in an explicit case.

\begin{table}[t]
\centering
\caption{Numerical comparison, independent of $\beta$. The columns evaluate $L_2$, $L_\infty$, $B_s$ and $\ell_H$ at the computed shape solutions or selected parameters. Bold marks the largest displayed lower values; $\widetilde\delta$ is the relative gap from their unrounded maximum to $\widetilde c_\ast$. The $U$ value is exact. Reproducible parameters and refinement diagnostics are given in \cref{app:table}.}\label{tab:main}
\setlength{\tabcolsep}{3.5pt}
\begin{tabular}{rcccccc}
\toprule
$\alpha$ & $\widetilde L_2$ & $\widetilde L_\infty$ & $\widehat B$ & $\widehat H$ & $\widetilde c_\ast$ & $\widetilde\delta$\\
\midrule
$-.75$ & $4.47796$ & $\mathbf{4.47864}$ & $\mathbf{4.47864}$ & $4.47795$ & $4.47953$ & $2.0\,10^{-4}$\\
$-.5\ (T)$ & $1.14206$ & $\mathbf{1.14246}$ & $\mathbf{1.14246}$ & $1.14206$ & $1.14269$ & $2.1\,10^{-4}$\\
$-.25$ & $.515817$ & $.516001$ & $\mathbf{.516006}$ & $.515817$ & $.516082$ & $1.5\,10^{-4}$\\
$0\ (P)$ & $.294010$ & $.293990$ & $\mathbf{.294080}$ & $.294010$ & $.294106$ & $8.8\,10^{-5}$\\
$.25$ & $.190247$ & $.189767$ & $\mathbf{.190263}$ & $.190247$ & $.190269$ & $2.8\,10^{-5}$\\
$.5\ (U)$ & $\mathbf{2/15}$ & $.131620$ & $.133333$ & $.133333$ & $2/15$ & $0$\\
$1$ & $.0755344$ & $.0725430$ & $.0755344$ & $\mathbf{.0760573}$ & $.0760850$ & $3.6\,10^{-4}$\\
$1.5$ & $.0481470$ & $.0449825$ & $.0481470$ & $\mathbf{.0491596}$ & $.0492267$ & $1.4\,10^{-3}$\\
$2$ & $.0331363$ & $.0301295$ & $.0331363$ & $\mathbf{.0343741}$ & $.0344777$ & $3.0\,10^{-3}$\\
$3$ & $.0182059$ & $.0157136$ & $.0182059$ & $\mathbf{.0194901}$ & $.0196375$ & $7.5\,10^{-3}$\\
\bottomrule
\end{tabular}
\end{table}

The leading computed values reflect the information retained. At the negative indices Keller and the near-pole Fredholm evaluations agree to the displayed precision. Fredholm leads at the intermediate indices; above $U$ the root constraint supplies the improvement over the free second-moment problem. The exact thresholds in \cref{thm:elliptic,thm:keller} identify activation of the mass constraint within each problem; the table also locates changes of the largest computed bound.

For the Chebyshev row, the certified bound $1.14268353624$ from Part~II gives a relative gap below $7.6\cdot10^{-6}$ to $c_{\sin}=1.1426921798\ldots$. The four variational columns display their contributions across the Jacobi family. Real shifted resolvents admit the same two-sided comparison (\cref{cor:realw}).

Analytic dependence gives an exact local comparison.
\begin{theorem}[Analyticity and contact]\label{thm:C}
The function $c_\ast$ is real-analytic on $(-1,\infty)$. For all $\alpha,\beta>-1$, $\ell_Y(\alpha)\le e(\alpha,\beta)\le c_\ast(\alpha)$, and, uniformly in $\beta$ as $\alpha\to1/2$,
\[
 0\le c_\ast(\alpha)-e(\alpha,\beta)
 \le c_\ast(\alpha)-\ell_Y(\alpha)=O((\alpha-1/2)^2).
\]
In particular $e(\cdot,\beta)$ is differentiable at $1/2$ and
$\partial_\alpha e(1/2,\beta)=c_\ast'(1/2)=-13/75$.
\end{theorem}
\begin{proof}
\Cref{app:analytic} proves analytic dependence by continuing all Bessel zeros on one complex disk and then continuing the uniformly controlled GL construction. The analytic nonnegative function $c_\ast-\ell_Y$ vanishes at $1/2$; its derivative vanishes there, giving the quadratic bound. Squeezing $e$ between the two functions gives its derivative. Finally
\[
 (\log\ell_Y)'=-\tfrac34\left(\frac2{\alpha+1}+\frac1{\alpha+2}\right),
 \qquad \ell_Y(1/2)=2/15,
\]
so that $\ell_Y'(1/2)=-13/75$.
\end{proof}
The two-sided bounds and their contact apply uniformly in the left index $\beta$. The exact $U$ value now gives a quantitative comparison of two implementations of a Nesterov residual.

\section{Nesterov residuals and Hessian drift}\label{sec:nesterov}

\subsection{Terminal-equivalent block realizations}

We apply the exact energy comparison to the ordering of accelerated steps. Jacobi factors already appear in semi-iterative regularization \cite{Sti58,Bra87,EHN96}, average-case acceleration \cite{PS20} and accelerated gossip \cite{BBG20}. For the Landweber iteration \cite{Lan51} with Nesterov acceleration \cite{Nes83}, convergence rates for linear ill-posed problems were studied by Neubauer \cite{Neu17}. Kindermann \cite{Kin21} identifies a Gegenbauer residual and interprets it as a combination of semi-iteration and Landweber damping. Its quadratic transformation gives an explicit block realization. We compare the sensitivity of these two recurrences while preserving their fixed-operator output and application count. In error coordinates the iteration is
\begin{equation}\label{eq:nesterov-iteration}
 e_1=(I-H_0)e_0,\qquad
 e_{s+1}=(I-H_s)\left(e_s+\frac{s-1}{s+b}(e_s-e_{s-1})\right),\qquad b>-1.
\end{equation}
Denote this realization by $\mathcal N_b$. For a fixed $H$ with spectrum in $[0,1]$, write $e_k=r_k(H)e_0$.
Kindermann's Gegenbauer formula \cite[Theorem~1]{Kin21}, followed by the quadratic transformation in \cref{app:nesterov}, gives
\begin{equation}\label{eq:nesterov-factor}
 \begin{gathered}
 k-1=2m+\epsilon,\qquad \epsilon\in\{0,1\},\qquad k-m=m+1+\epsilon,\\
 r_k(\lambda)=(1-\lambda)^{k-m} Q_m(\lambda),\qquad
 Q_m(\lambda)=\frac{P_m^{(b/2,\epsilon-1/2)}(1-2\lambda)}{P_m^{(b/2,\epsilon-1/2)}(1)}.
 \end{gathered}
\end{equation}
For each $J\in\cJ(P_m^{(b/2,\epsilon-1/2)})$, let $\mathcal B_{k,J}$ consist of its $m$ Jacobi steps \eqref{eq:jacobi-iteration}, calibrated on $[0,1]$, followed by $k-m$ Landweber steps $e\mapsto(I-H_s)e$. Its terminal residual is $r_k$, and its application count is $m+(k-m)=k$. At $m=0$ the implementation consists of the $k$ Landweber steps.

This identity also preserves the regularized output. Let $A:X\to Y$ be bounded between Hilbert spaces, let $x^\dagger$ be the minimum-norm solution of $Ax=y$, and suppose $\|y^\delta-y\|\le\delta$. Starting from $x_0=0$, both implementations have
\[
 x_k^\delta=g_k(A^*A)A^*y^\delta,\qquad
 g_k(\lambda)=\frac{1-r_k(\lambda)}{\lambda},
\]
with the polynomial continuation at zero. Suppose $\|A^*A\|\le1$ and the exact solution has the source representation $x^\dagger=(A^*A)^\mu w$, $0<\mu\le(b+1)/4$. The a priori choice
$k(\delta)=\lceil\delta^{-1/(2\mu+1)}\rceil\sim\delta^{-1/(2\mu+1)}$ as $\delta\downarrow0$ gives the common error
$O(\delta^{2\mu/(2\mu+1)})$ \cite[Theorem~4]{Kin21}.
The source exponent threshold is the edge index of the Jacobi factor, $\nu=(b+1)/4$. Its entire hard-edge limit and ordered small zeros are given in \cref{prop:nesterovedge}.

For a realization $\mathcal R$, denote its endpoint channel by
$a_k(\mathcal R)^2=\sum_{s<k}|K_s(0,0)|^2$.
The Jacobi-stage insertions in $\mathcal B_{k,J}$ are $-2d_i(J)$, and the Landweber-stage insertions are $-1$. Thus
\begin{equation}\label{eq:nesterov-block-energy}
 a_k(\mathcal B_{k,J})^2=4\E(J)+k-m,\qquad
 \lim_{k\to\infty}k^{-3}a_k(\mathcal B_{k,J_m^*})^2=\tfrac12c_*(b/2),
\end{equation}
and applying \cref{thm:transfer-main} to the two parity subsequences gives the lower bound $\tfrac12\ell(b/2)$ for the corresponding lower limiting optimum. The factor $1/2$ is $4(m/k)^3$. Equations \eqref{eq:nesterov-factor} and \eqref{eq:nesterov-block-energy} transfer the article's energy bounds to a specified implementation of the full Nesterov residual.

\subsection{An exact comparison at \texorpdfstring{$b=1$}{b=1}}

At $b=1$ the spectral supremum in \eqref{eq:gain-definition} can be evaluated for the original method and for the canonical even-horizon blocks. Put
\[
 F_j(\lambda)=(1-\lambda)^{j/2}U_j(\sqrt{1-\lambda}),\qquad 0\le\lambda\le1.
\]
This is a polynomial, with $F_0=1$ and $\max_{[0,1]}|F_j|=F_j(0)=j+1$.

\begin{theorem}[Sharp Nesterov and block gains]\label{thm:nesterov-gains}
For \eqref{eq:nesterov-iteration} with $b=1$,
\begin{equation}\label{eq:nesterov-sharp-gain}
 K_s^{\rm N}(\lambda,\lambda')=-\frac{F_{k-s-1}(\lambda)F_s(\lambda')}{k},\qquad
 \bigl(A_k^{\rm N}\bigr)^2=\frac{(k+1)^5-(k+1)}{30k^2}.
\end{equation}
For $k=2m+2$, $m\ge1$, let $J_m^U$ have zero diagonal and off-diagonal entries $1/2$. Then $\mathcal B_{k,J_m^U}$ has the same terminal residual and
\begin{equation}\label{eq:block-sharp-gain}
 A_k(\mathcal B_{k,J_m^U};[0,1])^2
 =4\E(J_m^U)+m+2
 =\frac{k^3}{15}+\frac{k}{2}+1-\frac{16}{15k}.
\end{equation}
Moreover
\[
 \lim_{m\to\infty}\frac{\displaystyle\inf_{J\in\cJ(U_m)}A_{2m+2}(\mathcal B_{2m+2,J};[0,1])}{(2m+2)^{3/2}}
 =\frac1{\sqrt{15}},\qquad
 \frac{A_k(\mathcal B_{k,J_m^U};[0,1])}{A_k^{\rm N}}\longrightarrow\sqrt2.
\]
\end{theorem}
\prooflink{thm:nesterov-gains}

The spectral endpoint $\lambda=\lambda'=0$ attains both gains. Nesterov distributes the Landweber action through its extrapolated queries at every step, giving constant $1/\sqrt{30}$; grouping the Landweber factors after a Jacobi stage has optimal constant $1/\sqrt{15}$. The full spectral gains above give the exact even-horizon comparison. \Cref{t:full} extends the endpoint optimality to all nonnegative-link Jacobi realizations, with arbitrary restarts and both parities. The fixed-operator regularization outputs coincide by \eqref{eq:nesterov-factor}.

For general $b>-1$, the original endpoint gain satisfies
\[
 \lim_{k\to\infty}k^{-3}a_k(\mathcal N_b)^2
 =\frac{2}{3(b+3)(2b+3)};
\]
\cref{lem:nesterov-endpoint} proves this by propagating a scalar insertion through \eqref{eq:nesterov-iteration}. Adaptive stopping and parameter selection for these accelerated regularizations are studied in \cite{Jin25}.

\subsection{The Bessel geometry of the original recurrence}\label{sec:bessel-geometry}

The comparison at $b=1$ mixes two effects: the length of the horizon on which the Jacobi factor is realized, and the shape of the resulting drift profile. They can be separated for every $b>-1$, because the endpoint profile of the original recurrence is itself a Green diagonal with the same limiting spectrum.

Put $\alpha=b/2$, the right index of the Jacobi factor in \eqref{eq:nesterov-factor}. The proof of \cref{lem:nesterov-endpoint} shows that the insertion at step $s\approx tk$ is $-k f_b(t)(1+o(1))$, with $f_b(t)=t(1-t^b)/b$. Consider the one-sided Bessel operator
\begin{equation}
 B_\alpha=-\partial_t^2+\frac{\alpha^2-1/4}{t^2}\quad\text{on }L^2((0,1),\dd t),
\end{equation}
with the branch $t^{\alpha+1/2}$ at $t=0$ and the Dirichlet condition at $t=1$. For $\alpha\ne0$, its adapted null solutions are $h_-(t)=t^{\alpha+1/2}$ and $h_+(t)=t^{1/2-\alpha}-t^{\alpha+1/2}$, with Wronskian $2\alpha=b$. At $\alpha=0$, take $h_+(t)=\sqrt t\log(1/t)$ and Wronskian one. The corresponding Green kernel is the product of the adapted solutions divided by their Wronskian; in particular,
\[
 G_B(t,t)=\frac{t-t^{2\alpha+1}}{2\alpha}=f_b(t),
\]
with $f_0(t)=t\log(1/t)$ by continuity. The eigenfunctions of $B_\alpha$ are $\sqrt t\,J_\alpha(j_{\alpha,k}t)$, so its spectrum is $\{\sigma_k\}=\{j_{\alpha,k}^2\}$, the spectrum of $H_\alpha$. The endpoint profile of $\mathcal N_b$ is therefore the Green diagonal of a second continuum realization of the Bessel-zero spectrum. This is consistent with \cref{prop:nesterovedge}: $y(z)=\Gamma(\alpha+1)(2/z)^\alpha J_\alpha(z)$ satisfies $y''+(b+1)z^{-1}y'+y=0$, which is the quadratic case of the ODE $\ddot X+(r/t)\dot X+\nabla f(X)=0$ of Su, Boyd and Cand\`es with friction $r=b+1$; their generalized scheme uses the momentum $(k-1)/(k+r-1)$, which is \eqref{eq:nesterov-iteration} \cite{SBC16}.

In the normalization of the article, where a realization with diagonal profile $\phi=2G(t,t)$ has energy $\int_0^1\phi(t)^2\dd t$, the one-sided realization has
\begin{equation}\label{eq:EB}
 E_B(\alpha)=\int_0^1\bigl(2G_B(t,t)\bigr)^2\dd t=\frac{8}{3(2\alpha+3)(4\alpha+3)},\qquad \alpha>-3/4.
\end{equation}
This range contains every Nesterov index, since $b>-1$ gives $\alpha>-1/2$. In particular the constant of \cref{lem:nesterov-endpoint} is $E_B(b/2)/4$. For $\alpha>-1/2$ the kernel $G_B$ is continuous, positive and factorized as in the proof of \cref{prop:green-coordinate}, so the same Green coordinate turns $2G_B(t,t)$ into a density $m_B$ with nonzero spectrum $\Lambda^{(\alpha)}$, mass $M$, root mass one and cost $E(m_B)=E_B(\alpha)$. Thus $m_B$ is admissible in each problem of \cref{sec:lower}, and $\ell(\alpha)\le E_B(\alpha)$ for $\alpha>-1/2$.

Combining this with \eqref{eq:nesterov-block-energy},
\begin{equation}
 \frac{\lim_kk^{-3}a_k(\mathcal B_{k,J_m^*})^2}{\lim_kk^{-3}a_k(\mathcal N_b)^2}
 =2\cdot\frac{c_*(b/2)}{E_B(b/2)}.
\end{equation}
The factor $2$ comes from compression. Taking $Z_s=I$ in the spectral form shows that the endpoint insertions sum to $r_k'(0)$, with absolute total asymptotic to $k^2M/2$. The original recurrence spreads this total over $k$ steps; the block uses $m\approx k/2$ active steps. Compressing a fixed shape into half the length while preserving its integral doubles its squared $\ell_2$ norm.

The factor $c_*/E_B$ compares the shapes $2G_\alpha(t,t)$ and $2G_B(t,t)$ on the unit interval. They have the same spectrum, mass and root mass. The canonical operator $H_\alpha$ has exponent $\gamma=\alpha/2+3/4$ at both ends; $B_\alpha$ has exponent $\alpha+1/2$ at zero and a regular Dirichlet end at one.

\begin{remark}[Two strings with one spectrum]\label{rem:power-string}
For $\alpha>0$, choose the additive constant in the Green coordinate of $G_B$ so that $x=\frac12\log(h_-/h_+)$. Then $\xi=\tanh x=2t^{2\alpha}-1$, and \cref{rem:string} turns $m_B$ into the power-law string
\[
 \rho_B(\xi)=\frac1{16\alpha^2}\Bigl(\frac{1+\xi}2\Bigr)^{1/\alpha-2},
 \qquad -1<\xi<1 .
\]
Its Dirichlet problem is solved by $\sqrt{1+\xi}\,J_\alpha\bigl(c\,(1+\xi)^{1/(2\alpha)}\bigr)$ \cite[10.13.3]{DLMF}, which recovers the spectrum $\{j_{\alpha,k}^2\}$. Direct integration recovers
\[
 \int_{-1}^1(1-\xi^2)\rho_B(\xi)\dd\xi=M,\qquad \int_{-1}^1\sqrt{\rho_B(\xi)}\dd\xi=1,\qquad \int_{-1}^1(1-\xi^2)^2\rho_B(\xi)^{3/2}\dd\xi=E_B(\alpha).
\]
 
So $c_*(\alpha)$ and $E_B(\alpha)$ are the energies of two strings with the same Dirichlet spectrum: the symmetric string $\rho_\alpha$ and the power-law string $\rho_B$. The latter is regular at $\xi=1$ and has exponent $1/\alpha-2$ at $\xi=-1$. Its exponent matches the canonical exponent $1/\nu-2=4/(2\alpha+1)-2$ of \cref{rem:string} only at $\alpha=1/2$, where both strings are homogeneous.
 
A string is not determined by its spectrum. The missing data are its norming constants, the continuum counterpart of the spectral weights; Krein's correspondence between strings and spectral functions \cite{KK74} accounts for them.
\end{remark}

The two exponents agree precisely at $\alpha=1/2$, where $H_{1/2}=B_{1/2}$ is the Dirichlet Laplacian and $c_*=E_B=2/15$. The factor $\sqrt2$ in \cref{thm:nesterov-gains} therefore comes from the length of the horizon. For the other tabulated Nesterov indices, Table~\ref{tab:main} and \eqref{eq:EB} give $c_*<E_B$, so canonical blocks recover part of the compression loss. \Cref{t:improve} gives an analytic strict improvement over $E_B$ throughout $\alpha>-1/2$, $\alpha\ne1/2$.

The ratio $2c_*(b/2)/E_B(b/2)$ is stationary at $b=1$. Logarithmic differentiation of \eqref{eq:EB} gives $E_B'(1/2)/E_B(1/2)=-\tfrac24-\tfrac45=-\tfrac{13}{10}$, the logarithmic derivative of $\ell_Y$ at $1/2$ computed in the proof of \cref{thm:C}. Hence $E_B'(1/2)=-13/75$, and $E_B$ is tangent at $U$ to $\ell_Y$, to $c_*$ and to $e$. Since the lower limiting optimum of the block endpoint channel lies between $\tfrac12\ell(b/2)$ and $\tfrac12c_*(b/2)$, its ratio to the Nesterov endpoint channel lies between $2\ell_Y(b/2)/E_B(b/2)$ and $2c_*(b/2)/E_B(b/2)$, and both ends equal $2+O((b-1)^2)$. The lower end is explicit:
\[
 \rho(b)=\frac{2\ell_Y(b/2)}{E_B(b/2)}
 =\frac{3\sqrt3\,5^{3/4}}{10}\,\frac{(b+3)(2b+3)}{(b+2)^{3/2}(b+4)^{3/4}},\qquad
 (\log\rho)'(b)=-\frac{(b-1)(2b^2+11b+18)}{4(b+2)(b+3)(b+4)(2b+3)}.
\]
The quadratic factor is positive, so $\rho$ increases on $(-1,1)$, attains its maximum $\rho(1)=2$ and then decreases; $\rho(-1^+)\approx1.52$, and $\rho>1$ for $-1<b<b_+$ with $b_+\approx139.8$. On this range, $\mathcal N_b$ has a strict asymptotic endpoint-gain advantage over all nonnegative-link Jacobi realizations (\cref{t:full}). For the block class, the full $\ell$ at the tabulated indices with $-1/2\le b\le6$ gives $2\ell/E_B\ge1.93$.

For general $b$ we compare endpoint channels; the full spectral gain \eqref{eq:gain-definition} is evaluated at $b=1$. We now extend the endpoint comparison to all three-term realizations and improve the standard realization by an isospectral flow.

% Part of the Nesterov section; compile jacobi.tex.

\subsection{All three-term realizations}\label{ss:closure}

To include restarts, allow nonnegative off-diagonal entries and denote the
resulting fibre by $\overline{\cJ}(p)$. Here $p$ may have repeated zeros in
$[-1,1)$. Every matrix in this fibre satisfies $-I\le J<I$, so all leading
determinants $P_i(1)$ are positive. Thus \eqref{eq:jacobi-iteration} remains
defined, with $a_i=0$ exactly when $\delta_i=0$.

For the Nesterov residual \eqref{eq:nesterov-factor}, put
\[
 p_k(z)=(z+1)^{k-m}\widehat P_m^{(b/2,\epsilon-1/2)}(z),
 \qquad k-1=2m+\epsilon,\quad\epsilon\in\{0,1\},
\]
where the hat denotes monic normalization. Then
$p_k(1-2\lambda)/p_k(1)=r_k(\lambda)$, and
$\overline{\cJ}(p_k)$ parametrizes its nonnegative-link three-term
realizations. Positive-link Jacobi matrices have simple spectrum; hence
restarts are necessary when $k-m\ge2$. The closed fibre contains the block
matrices $J\oplus(-I_{k-m})$ underlying $\mathcal B_{k,J}$, and the matrices
corresponding to interleavings of Landweber steps and restarted Jacobi
segments. Write $a_k(J)$ for the endpoint channel.

\begin{lemma}\label{l:closure}
The drift identity of \cref{prop:drift-energy} extends to
$\overline{\cJ}(p)$; in particular $a_k(J)^2=4\E(J)$ on $[0,1]$.
The transfer theorem \cref{thm:transfer-main} also holds for closed fibres,
with zeros in $[-1,1)$ counted with multiplicity. Each finite closed fibre
has an attained energy minimum.
\end{lemma}
\begin{proof}
The drift identity uses positive leading determinants and the diagonal
cofactor formula, both valid for zero links. For transfer, approximate $J$
by $(1-\eta)J$, whose norm is less than one. In each fixed dimension the
resolvent, its spectrum and $\E$ are continuous at $J<I$. Choosing
$\eta_n\downarrow0$ sufficiently fast preserves the limiting scaled spectrum
and normalized energy. Apply \cref{prop:packing}, which already permits
zero links, and the compactness argument of \cref{thm:transfer-main}.
Finally, the closed fibre is compact, and its fixed spectrum stays below
one; continuity of $\E$ gives attainment.
\end{proof}

Retain the three amplitude-homogeneous lower bounds by writing
\[
 \ell^\circ(\alpha)=\max\left\{L_2,L_\infty,
                 \sup_{0<s<1/\Lambda_1}B_s\right\}
 \quad\text{at }\Lambda^{(\alpha)}.
\]
Thus $\ell_Y\le\ell^\circ\le\ell$. Define the optimal squared endpoint channel
$C_k(b)=\min_{J\in\overline{\cJ}(p_k)}a_k(J)^2$.

\par
\begin{samepage}
\begin{theorem}[All three-term realizations]\label{t:full}
For every $b>-1$,
\[
 \liminf_{k\to\infty}\frac{C_k(b)}{k^3}\ge\tfrac12\ell^\circ(b/2).
\]
At $b=1$ this is sharp:
\[
 \lim_{k\to\infty}\frac{C_k(1)}{k^3}=\frac1{15}
       =2\lim_{k\to\infty}\frac{a_k(\mathcal N_1)^2}{k^3}.
\]
The canonical blocks $\mathcal B_{k,J_m^*}$ asymptotically attain this
minimum in both parities.
\end{theorem}
\end{samepage}
\begin{proof}
The active roots contribute $(m/k)^2$ times the Jacobi scaled resolvent
spectrum. The $k-m$ roots at $-1$ contribute a total $(k-m)/(2k^2)$.
Since $m/k\to1/2$, \cref{lem:spec} gives the full $\ell^1$ limit
$\Lambda^{(b/2)}/4$. By \cref{l:closure},
\[
 \liminf_k k^{-3}C_k(b)\ge4\ell(\Lambda^{(b/2)}/4)
                         \ge\tfrac12\ell^\circ(b/2).
\]
For the last inequality, multiplying a quarter-spectrum density by four
multiplies its mass and top eigenvalue by four, its second moment by
sixteen and its cost by eight. The determinant constraint transforms by
$\hD_s(4m)=\hD_{4s}(m)$. This lifts competitors for each retained bound to
the original data.

At $b=1$, \cref{sec:second} gives $\ell^\circ(1/2)=2/15$. Conversely,
\[
 a_k(\mathcal B_{k,J_m^*})^2=4\E(J_m^*)+k-m,
 \qquad \E(J_m^*)\sim2m^3/15
\]
by \cref{l:closure,thm:A}. Together with \cref{lem:nesterov-endpoint},
these identities give the stated limit and comparison.
\end{proof}

Thus at $b=1$ the regrouping in \cref{thm:nesterov-gains} is asymptotically
optimal among all these three-term realizations, with arbitrary restarts
and placement of Landweber steps. The factor $\sqrt2$ is the advantage of
the extrapolated query over this entire class.

\subsection{The standard realization}\label{ss:standard}

Fix $\alpha>-1/2$ and $\beta>-1$. Let
$J_n^{\rm st}\in\cJ(P_n^{(\alpha,\beta)})$ be the matrix of the orthonormal
Jacobi recurrence for the probability measure proportional to
$(1-x)^\alpha(1+x)^\beta\dd x$ on $[-1,1]$. Its spectral weights are the
normalized Gauss--Jacobi weights. Write
\[
 \phi_\alpha^B(t)=2G_B(t,t)=2t^{2\alpha+1}\int_t^1v^{-2\alpha-1}\dd v
 =\begin{cases}t(1-t^{2\alpha})/\alpha,&\alpha\ne0,\\
                 2t\log(1/t),&\alpha=0.\end{cases}
\]
For $0<s\le t\le1$, $2G_B(s,t)=(s/t)^{\alpha+1/2}\phi_\alpha^B(t)$;
the kernel extends symmetrically and continuously to $[0,1]^2$.
Thus $\phi_\alpha^B=2f_b$ for $b=2\alpha$, and
\[
 \int_0^1\phi_\alpha^B(t)\dd t=M=\frac1{2(\alpha+1)},\qquad
 \int_0^1\phi_\alpha^B(t)^2\dd t=E_B(\alpha).
\]

\begin{theorem}[The standard one-sided limit]\label{t:std}
Set $R^{\rm st}=(I-J_n^{\rm st})^{-1}$ and extend its entries by zero at
the terminal index $n$. For each $\delta>0$,
\[
 n^{-1}R^{\rm st}_{\lfloor sn\rfloor,\lfloor tn\rfloor}
       \longrightarrow2G_B(s,t)
 \quad\text{uniformly for }\delta\le s,t\le1.
\]
Moreover $\sup_{n\ge1,\,i<n}d_i(J_n^{\rm st})/n<\infty$ and
$n^{-3}\E(J_n^{\rm st})\to E_B(\alpha)$. In particular,
$e(\alpha,\beta)\le E_B(\alpha)$.
\end{theorem}
\begin{proof}
Let $u_i=q_i(1)>0$ for the orthonormal Jacobi polynomials, including the
next polynomial $q_n$, and set $c_i=a_{i+1}u_iu_{i+1}$. The discrete
Green formula, or the first-difference factorization of $I-J_n^{\rm st}$,
gives
\begin{equation}\label{e:green}
 R^{\rm st}_{ij}=u_iu_j\sum_{l=\max(i,j)}^{n-1}c_l^{-1}.
\end{equation}
This is the standard tridiagonal inverse construction \cite{Meu92}.
The Jacobi recurrence and norm formulas \cite[\S\S18.3,18.9]{DLMF},
together with the Gamma-ratio asymptotic, yield
\begin{equation}\label{e:stdasym}
 a_i\longrightarrow\tfrac12,\qquad u_i\sim\kappa i^{\alpha+1/2}
 \quad(\kappa>0).
\end{equation}
For indices $i/n\to s>0$, $j/n\to t>0$, the tail sum in
\eqref{e:green} is a Riemann sum for $2G_B(s,t)$.
This moving-mesh limit gives uniform convergence on $[\delta,1]^2$
by compactness, including the zero terminal row.

The same power bounds give $u_i^2/c_l\le C$ for $i\le l$, hence
$d_i\le Cn$. Positivity gives $|R_{ij}|^2\le d_id_j$.
The first $\delta n$ indices therefore contribute $O(\delta)$ to the
normalized energy. On the remaining indices the diagonal converges
uniformly, so $n^{-3}\E(J_n^{\rm st})\to\int_0^1\phi_\alpha^B(t)^2\dd t$.
\end{proof}

Thus $\mathcal N_b$ spreads the standard Jacobi endpoint profile over the
full horizon.

\subsection{An isospectral descent}\label{ss:flow}

For $\alpha\ne1/2$, an explicit flow of the standard spectral weights
supplies competitors of strictly lower asymptotic cost.

\begin{proposition}[Resolvent flow]\label{p:flow}
For $J\in\cJ(p)$ of size $n\ge1$, put $R=(I-J)^{-1}$ and
$B=n^{-2}(R_--R_+)$, using strict lower and upper triangular parts. The equation
$\dot J=[B,J]$ defines a global flow on $\cJ(p)$ with constant
$s_n=n^{-2}\Tr R$ and weights
\begin{equation}\label{e:tilt}
 w_j(\tau)=
 \frac{w_j(0)e^{-2\tau/[n^2(1-\zeta_j)]}}
      {\sum_l w_l(0)e^{-2\tau/[n^2(1-\zeta_l)]}}.
\end{equation}
Along this flow,
\begin{equation}\label{e:flowder}
 \dot\E=\frac4{n^2}\sum_{i>j}(d_i-d_j)R_{ij}^2,\qquad
 |\dot\E|\le4s_n\E,\qquad |\ddot\E|\le24s_n^2\E.
\end{equation}
Consequently, for every real $\tau$,
\begin{equation}\label{e:taylor}
 \E(\tau)\le\E(0)+\tau\dot\E(0)
              +12s_n^2e^{4s_n|\tau|}\tau^2\E(0).
\end{equation}
For roots $\zeta_1>\cdots>\zeta_n$, the matrix limits at positive and
negative infinite time are respectively
$\operatorname{diag}(\zeta_n,\ldots,\zeta_1)$ and
$\operatorname{diag}(\zeta_1,\ldots,\zeta_n)$.
\end{proposition}
\begin{proof}
Reconstruct the Jacobi matrix from the positive weights \eqref{e:tilt}
using the analytic inverse-spectral chart \cite{dBG}. Differentiating the
orthogonal recurrence gives $\dot J=[B,J]$; since $RJ=JR$, the commutator
is tridiagonal and
$\dot J_{i,i+1}=n^{-2}(d_i-d_{i+1})J_{i,i+1}$.
The weight formula gives a solution for every real time.

Now $\dot R=[B,R]$ and
$\dot d_i=2n^{-2}\sum_l\operatorname{sgn}(i-l)R_{il}^2$.
Pairing indices gives the first identity in \eqref{e:flowder}.
The inverse of the positive definite matrix $I-J$ is entrywise
nonnegative. Thus, with $\nrm R\le n^2s_n$,
\[
 |\dot R_{ij}|\le2n^{-2}(R^2)_{ij},\qquad
 |\dot d_i|\le2s_nd_i,\qquad |\ddot d_i|\le8s_n^2d_i.
\]
Summing the derivatives of $\E=\sum_i d_i^2$ proves both bounds.
Gronwall's inequality and Taylor's integral remainder give \eqref{e:taylor}.

For sorting, the weight ratios in \eqref{e:tilt} successively select the
roots in the stated order. The monic orthogonal polynomial of degree $l$
minimizes the weighted squared norm. Comparing with the polynomial
vanishing at the first $l$ selected roots forces convergence to their
product. The resulting prefix traces give the ordered diagonal entries;
the conserved $\Tr J^2$ then forces all off-diagonal entries to zero.
\end{proof}

\par
\begin{samepage}
\begin{theorem}[Strict improvement]\label{t:improve}
Start the flow at $J_n^{\rm st}$. Then
\[
 \frac{\dot\E_n(0)}{n^3}\longrightarrow
 D_\alpha=\frac{32(2\alpha-1)}
 {5(\alpha+1)(2\alpha+3)(4\alpha+3)(2\alpha+5)(4\alpha+5)(6\alpha+5)}.
\]
For $\alpha\ne1/2$, the standard matrix is therefore noncritical for all
sufficiently large $n$. Put
\[
 \tau_\alpha=\frac{(\alpha+1)(1-2\alpha)}
 {5(2\alpha+5)(4\alpha+5)(6\alpha+5)},\qquad
 \Theta_\alpha=\frac{16(2\alpha-1)^2}
 {25(2\alpha+3)(4\alpha+3)(2\alpha+5)^2(4\alpha+5)^2(6\alpha+5)^2}.
\]
The flowed matrices satisfy
\[
 \limsup_n n^{-3}\E(J_n(\tau_\alpha))\le E_B(\alpha)-\Theta_\alpha.
\]
In particular, when $\alpha\ne1/2$,
\[
 e(\alpha,\beta)\le
 \limsup_n\inf_{J\in\cJ(P_n^{(\alpha,\beta)})}n^{-3}\E(J)
 \le E_B(\alpha)-\Theta_\alpha<E_B(\alpha).
\]
\end{theorem}
\end{samepage}
\begin{proof}
The uniform limit and global bound in \cref{t:std} turn
\eqref{e:flowder} into
\[
 D_\alpha=4\int_0^1t^{-2\alpha-1}\phi_\alpha^B(t)^2
             \int_0^t[\phi_\alpha^B(t)-\phi_\alpha^B(u)]u^{2\alpha+1}\dd u\dd t.
\]
Integrating the powers gives the displayed rational value. At $\alpha=0$,
the logarithmic moments give $D_0=-32/5625$ directly.

Write $C=E_B(\alpha)$ and $s=M=1/(2(\alpha+1))$.
Since $n^{-3}\E_n(0)\to C$ and $s_n\to s$, \eqref{e:taylor} gives
\[
 \limsup_n n^{-3}\E_n(\tau)
       \le C+\tau D_\alpha+12s^2e^{4s|\tau|}\tau^2C.
\]
The choice $\tau_\alpha=-D_\alpha/(48s^2C)$ satisfies
$4s|\tau_\alpha|\le1/15$, so its exponential factor is less than two.
The right side is consequently at most
\[
 C+\tau_\alpha D_\alpha+24s^2C\tau_\alpha^2
       =C-\frac{D_\alpha^2}{96s^2C}=C-\Theta_\alpha.
\]
Every flowed matrix remains in the same finite fibre.
\end{proof}

For comparison with the original recurrence, write
\[
 C_k^{\rm blk}(b)=
 \inf_{J\in\cJ(P_m^{(b/2,\epsilon-1/2)})}a_k(\mathcal B_{k,J})^2.
\]
The minimum defining $C_k$ is attained; $C_k^{\rm blk}$ is an infimum over
positive-link factors.

\begin{corollary}[Optimal endpoint comparison]\label{c:sandwich}
For $b>-1$, $b\ne1$, and $\alpha=b/2$,
\begin{align*}
 \rho(b)\le\frac{2\ell^\circ(\alpha)}{E_B(\alpha)}
 &\le\liminf_k\frac{C_k(b)}{a_k(\mathcal N_b)^2}
 \le\limsup_k\frac{C_k^{\rm blk}(b)}{a_k(\mathcal N_b)^2}\\
 &\le2-\frac{12(b-1)^2}{25(b+5)^2(2b+5)^2(3b+5)^2}<2.
\end{align*}
At $b=1$, both optimal ratio sequences converge to $2$.
\end{corollary}
\begin{proof}
The lower bounds follow from \cref{t:full,lem:nesterov-endpoint} and
$\ell_Y\le\ell^\circ$. Since $C_k\le C_k^{\rm blk}$, it remains to bound
the blocks. Use the flowed competitor of \cref{t:improve} in each parity
and $m/k\to1/2$ to obtain
\[
 \limsup_k k^{-3}C_k^{\rm blk}(b)
       \le\tfrac12(E_B(\alpha)-\Theta_\alpha).
\]
Divide by the positive Nesterov limit $E_B(\alpha)/4$ and simplify.
At $b=1$, the lower bound and canonical blocks give the common limit.
\end{proof}

On the range $-1<b<b_+$, where $\rho(b)>1$, every three-term realization
has asymptotically larger endpoint gain than $\mathcal N_b$.
The optimal squared-gain penalty is strictly below two away from $b=1$.

The standard and persymmetric matrices realize two different geometries
of the same spectrum. The former gives $B_\alpha$; the latter is critical
in every finite fibre and gives the reflection-symmetric $H_\alpha$.
For the $U_n$ polynomials they coincide, consistently with $D_{1/2}=0$.

In the coordinates $\xi_j=\log(w_j|p'(\zeta_j)|)$ modulo constants, the
flow is the straight line
\[
 \xi_j(\tau)=\xi_j(0)-\frac{2\tau}{n^2(1-\zeta_j)}.
\]
It reaches the persymmetric matrix exactly when this line passes through
the origin. With two roots the line covers the whole positive fibre, so
every realization, including the persymmetric one, is reached in finite time.

At its two ends lie the oppositely sorted Richardson schedules, with
common energy $\sum_j(1-\zeta_j)^{-2}$. For Jacobi roots,
\cref{lem:spec} gives
\[
 n^{-4}\sum_j(1-\zeta_j)^{-2}\longrightarrow
 T_2=\frac1{4(\alpha+1)^2(\alpha+2)}.
\]
Thus the cubic-energy momentum realizations lie between two schedules of
quartic energy.

\begin{remark}[Peakons]\label{rem:peakons}
Through \cref{rem:beads}, the flow of \cref{p:flow} is the multipeakon flow of the Camassa--Holm equation.
 
In the normalization of \cite[(1.1), (2.2)]{BSS00}, an $n$-peakon $u=\frac12\sum_im_ie^{-2|x-x_i|}$ with $m_i>0$ evolves by
\[
 \dot x_i=u(x_i),\qquad
 \dot m_i=m_i\sum_jm_j\operatorname{sgn}(x_i-x_j)\,e^{-2|x_i-x_j|},
\]
and its Lax operator has the kernel $e^{-|x-y|}$ of $\sK$. Start from $m_i=n^{-2}d_i$ and $x_i=X_i$ for $J\in\cJ(p)$, and put $\Gamma(t)=[\sqrt{m_im_j}\,e^{-|x_i-x_j|}](t)$. Then
\[
 \Gamma(t)=n^{-2}\bigl(I-J(t/2)\bigr)^{-1},
\]
where $J(\tau)$ is the flow of \cref{p:flow}.
 
To verify this, recall that the residues $a_j$ of the Weyl function of the associated string \cite[(4.5)--(4.6)]{BSS00} evolve by $a_j(t)=a_j(0)e^{-2t/\lambda_j}$ \cite[(4.10)]{BSS00}. Here $\lambda_j=-2/r_j$, and $r_j=[n^2(1-\zeta_j)]^{-1}$ are the constant eigenvalues of $\Gamma$. Green's identity gives $a_j|\lambda_j|=f_j'(1)^2$ for the eigenfunction $f_j$ normalized in $L^2((-1,1),\rho(\xi)\dd\xi)$. Since $f_j$ is linear on the last segment, \cref{rem:beads} shows that $a_j|\lambda_j|$ is proportional to $w_j^\vee$, with a factor independent of $j$. Hence $w_j^\vee(t)\propto w_j^\vee(0)e^{tr_j}$. By the formula for $w^\vee$ preceding \eqref{lem:persym}, this is \eqref{e:tilt} at $\tau=t/2$, and the weights determine $J$.
 
Beals, Sattinger and Szmigielski realize the same flow of spectral data as the isospectral Jacobi flow generated by $-1/\lambda$. That flow belongs to the hierarchy containing the finite Toda lattice \cite[\S5, (5.1)]{BSS01}.
 
For $\mu_t=\sum_i m_i(t)\delta_{x_i(t)}$, the invariants of \cref{p:flow} are Camassa--Holm conservation laws:
\[
 s_n=\Tr\Gamma=\sum_im_i=\int_\R 1\dd\mu_t(x),\qquad
 \Tr\Gamma^2=\int_\R\bigl(4u(x)^2+u_x(x)^2\bigr)\dd x.
\]
The sorting limits of the flow are the separation of the peakons. The matrix $R$ becomes diagonal, so the gaps $X_{i+1}-X_i$ diverge. As $t\to+\infty$ the limiting masses $r_j$ increase from left to right. This matches the fact that an isolated peakon travels with speed equal to its height $m_i/2$.
\end{remark}

\begin{figure}[t]
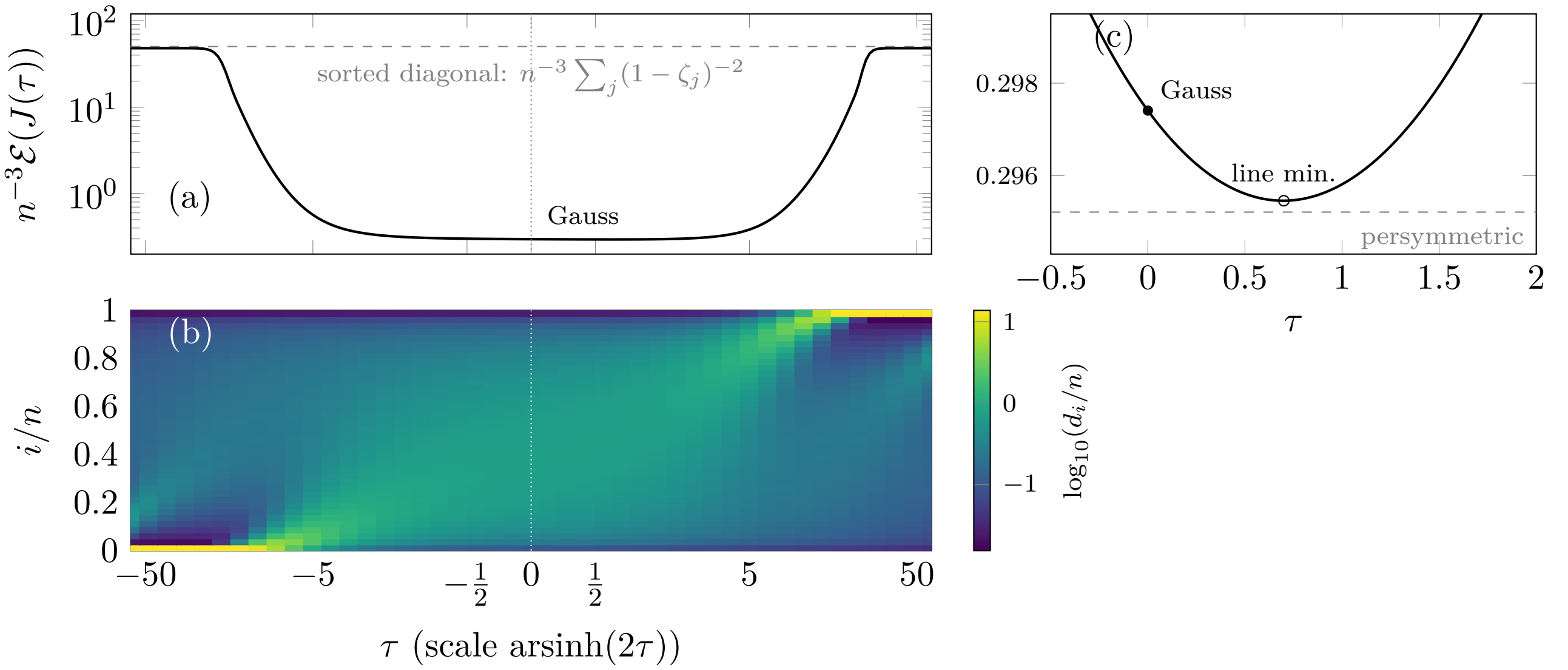

\centering

\caption{Numerical illustration of the resolvent flow through the Gauss
realization of the Legendre zeros ($\alpha=\beta=0$, $n=400$), reconstructed
from \eqref{e:tilt} by Lanczos with full reorthogonalization.
(a) Normalized energy on the scale $\operatorname{arsinh}(2\tau)$; the dashed
line is the limiting sorted-diagonal energy.
(b) The diagonal profile $d_i/n$, averaged over $40$ blocks of vertices.
At zero time it approximates $\phi_0^B$; a sorting front develops as
$|\tau|$ grows.
(c) The energy near its minimum on this orbit. The negative initial slope
is consistent with $D_0=-32/5625$. The plotted minimum near $\tau=0.70$
recovers about $89\%$ of the finite Gauss--persymmetric energy gap;
the dashed reference is the persymmetric energy.}
\label{f:toda}
\end{figure}

At $b=0$ the relative squared-gain improvement is
$2\Theta_0/E_B(0)=12/390625$.
Figure~\ref{f:toda} illustrates the energy along one flow orbit.

\section{Two families with the same limiting spectrum}

The general transfer theorem uses the limiting edge spectrum. The following two families have that same spectrum and different limiting finite-fibre optima. Let $L_n^{(a)}$ denote the Laguerre polynomial with parameter $a>-1$.

\begin{proposition}[The same spectrum, different optimal energies]
Fix $0<\theta\le1/4$, and put
\[
 p_n^{\rm U}(x)=(\theta/2)^nU_n((x-1+\theta)/\theta),\quad
 p_n^{\rm L}(x)=\frac{n!}{H_n^n}L_n^{(1/2)}(H_n(1-x)),\quad
 H_n=\frac{2n+3/2}{4\theta}.
\]
Both polynomials have simple zeros in $(-1,1)$, and both scaled resolvent spectra converge in $\ell^1$ to $\Lambda_\theta=(2/(\theta\pi^2k^2))_{k\ge1}$. Their optimal energies satisfy
\[
 \lim_n n^{-3}\inf_{\cJ(p_n^{\rm U})}\E=\frac2{15\theta^2},\qquad
 \lim_n n^{-3}\min_{\cJ(p_n^{\rm L})}\E=\frac1{9\theta^2}.
\]
\end{proposition}
\begin{proof}
The affine map $J\mapsto(1-\theta)I+\theta J$ multiplies the resolvent by $\theta^{-1}$. The $U$ case follows from its exact upper and lower limits, and its spectral limit from \cref{lem:spec}.

For the other family we use a general factorization observation. Suppose $Y=BB^\top>0$ is Jacobi and $B$ is its lower bidiagonal Cholesky factor with positive diagonal, and $\bar A_i=B_{ii}^2$, $\bar C_i=B_{i,i-1}^2$ satisfy $\bar A_i-\bar C_i=\kappa>0$, with $\bar C_0=0$. The squared row norms of $B^{-1}$ obey
\[
 r_0=1/\bar A_0,\qquad r_i=(1+\bar C_i r_{i-1})/\bar A_i;
 \qquad r_i=1/\kappa.
\]
Indeed the two terms in the row recursion for $B^{-1}$ have disjoint new and old support. Thus $(B^\top B)^{-1}$ has constant diagonal $1/\kappa$. If $H>\lambda_{\max}(Y)/2$ and $D=\diag((-1)^i)$, then $J=I-DB^\top BD/H$ is in the fibre with zeros $1-\lambda_j(Y)/H$ and has constant resolvent diagonal $H/\kappa$. The fixed trace and Cauchy--Schwarz show that its energy $nH^2/\kappa^2$ is the exact minimum.

For the Laguerre polynomial $L_n^{(a)}$, $a>-1$, its Jacobi recurrence matrix has diagonal $2i+a+1$ and off-diagonal $\sqrt{i(i+a)}$ between $i-1$ and $i$. Thus $\bar A_i=i+a+1$ and $\bar C_i=i$, so $\kappa=a+1$. Gershgorin, using $\sqrt{i(i+a)}\le i+a/2$ for $i\ge1$ and treating row zero separately, gives $\lambda_{\max}<4n+2a+2$. At $a=1/2$ and $\theta\le1/4$ this is at most $2H_n$, proving admissibility and the limit $1/(9\theta^2)$.

Let $y_{n,1}<\cdots<y_{n,n}$ be the zeros of $L_n^{(1/2)}$. Laguerre Mehler--Heine \cite[18.11.6]{DLMF} gives $4ny_{n,k}\to\pi^2k^2$ for every fixed $k$. Hence $H_n/(n^2y_{n,k})\to2/(\theta\pi^2k^2)$. The exact trace $\sum_k1/y_{n,k}=n/(3/2)$ gives convergence of their sums to $1/(3\theta)$, the sum of $\Lambda_\theta$. Positivity upgrades this to $\ell^1$, as in \cref{lem:spec}.
\end{proof}

Both families have root mass one under packing. Their different optimal energies show that the finite-fibre structure contributes to the limiting optimum in addition to the common spectral data. For Jacobi zeros, the canonical construction and the variational bounds quantify this optimum through the two-sided estimate of \cref{thm:B}.

\section*{Conclusion}

The estimate $\ell(\alpha)\le e(\alpha,\beta)\le c_*(\alpha)$ links finite
isospectral optimization to the Green geometry of its continuum limits.
It is exact at $\alpha=1/2$, with energy $2/15$ and derivative $-13/75$;
the gap between the bounds is $O((\alpha-1/2)^2)$. The two polynomial families
with identical limiting edge spectra and different optimal energies show
that the finite-fibre structure contributes additional information to the
limiting optimum.

For accelerated iterations, the same energy controls endpoint Hessian drift.
At $b=1$, the sharp comparison quantifies the $\sqrt2$ advantage of Nesterov's
extrapolated query over all nonnegative-link Jacobi realizations, including
arbitrary restarts and both horizon parities. The resolvent flow gives a
strict asymptotic improvement of the standard one-sided realization for
$\alpha>-1/2$, $\alpha\ne1/2$.

In the string form of \cref{rem:string}, the limiting problem asks for the string on $(-1,1)$ with prescribed Dirichlet spectrum that minimizes a hyperbolically weighted $L^{3/2}$ norm of its density. The canonical construction selects the symmetric string, which at $U$ is the homogeneous one. Krein's spectral theory of strings \cite{KK74} parametrizes all strings with a given spectrum by their norming constants, and so offers a natural setting for the remaining gap between $\ell$ and $c_*$.

\textbf{Use of AI tools.} During the preparation of this article, the authors used
the AI assistants Claude Opus 5.5 and ChatGPT-6.1 Sol as auxiliary tools for the
following purposes: (i) editing the manuscript and providing recommendations on
how to organize the text so that its content is better explained on first reading;
(ii) preparing the Lean package for formal verification of the proofs;
(iii) obtaining suggestions on the use of externally published, openly available
results to make the proofs more compact and to ensure proper attribution when
such results had been obtained earlier by other authors, as well as checking the
correctness of the original cited works and the results reported therein; and (iv) identifying connections with the classical spectral theory of vibrating strings and with the related literature, which informed the interpretive remarks in the text.
The authors have verified all statements in this article and take responsibility
for all results and accompanying code.

\appendix
\section{Residuals and first-variation kernels}\label{app:nesterov}

\subsection{Jacobi factors and the hard edge}

Let $C_j^{(\rho)}$ denote the Gegenbauer polynomial. For $\nu>0$, Brakhage's $\nu$-methods have residuals \cite{Bra87,Han91}
\begin{equation}
 r_k^{(\nu)}(\lambda)=\frac{C_{2k}^{(2\nu)}(\sqrt{1-\lambda})}{C_{2k}^{(2\nu)}(1)}
 =\frac{P_k^{(2\nu-1/2,-1/2)}(1-2\lambda)}{P_k^{(2\nu-1/2,-1/2)}(1)}.
\end{equation}
Thus their right index is $\alpha=2\nu-1/2$, and their edge index $\alpha/2+1/4=\nu$ is the qualification parameter of the regularization method \cite{EHN96}.
For \eqref{eq:nesterov-iteration}, Kindermann's formula is \cite[Theorem~1]{Kin21}
\begin{equation}
 r_k(\lambda)=(1-\lambda)^{(k+1)/2}
 \frac{C_{k-1}^{((b+1)/2)}(\sqrt{1-\lambda})}{C_{k-1}^{((b+1)/2)}(1)}.
\end{equation}
The quadratic Gegenbauer transformations \cite[18.7.15--16]{DLMF}, normalized at one, give \eqref{eq:nesterov-factor}. For $k-1=2m+\epsilon$, the factor of degree $m$ has indices $(b/2,\epsilon-1/2)$; when $\epsilon=1$, the square root in the transformation combines with the prefactor to give the integer exponent $k-m=m+1+\epsilon$. The Jacobi realization in \eqref{eq:jacobi-iteration} has that normalized factor as terminal polynomial, and the $k-m$ Landweber steps multiply it by $(1-\lambda)^{k-m}$.

\begin{proposition}[The Nesterov hard edge]\label{prop:nesterovedge}
For $b>-1$, locally uniformly for $z\in\C$,
\begin{equation}\label{eq:nesterov-edge}
 r_k(z^2/k^2)\longrightarrow\Gamma(b/2+1)(2/z)^{b/2}J_{b/2}(z),
\end{equation}
with its entire-series value at zero. If $\lambda_{k,j}$ is the $j$th smallest zero of $r_k$ in $(0,1)$, then $k^2\lambda_{k,j}\to j_{b/2,j}^2$.
The Jacobi block and Landweber realization of \cref{sec:nesterov} uses exactly $k$ operator applications.
\end{proposition}
\begin{proof}
For either parity put $\theta_k(z)=2\arcsin(z/k)$. Then
$1-2z^2/k^2=\cos\theta_k(z)$ and $m\theta_k(z)\to z$ uniformly on compact sets. Normalized Mehler--Heine gives \eqref{eq:nesterov-edge} for the Jacobi factor; the Landweber prefactor tends uniformly to one because $k-m=O(k)$. The two parities have the same limit. The zeros in $(0,1)$ are exactly those of the Jacobi factor. Hurwitz's theorem at each simple positive Bessel zero, together with uniform convergence on the intervening compact intervals, fixes their ordering and proves the asserted limit. Finally $m+(k-m)=k$, and every step applies the operator once.
\end{proof}

\subsection{The exact gain comparison}

\begin{proof}[Proof of \Cref{thm:nesterov-gains}]\refstepcounter{prooflink}\label{proof:thm:nesterov-gains}
For $b=1$, set $v_s=se_s$, with $v_0=0$. On a fixed spectral mode the recurrence is
\[
 v_{s+1}=2(1-\lambda)v_s-(1-\lambda)v_{s-1}.
\]
Its fundamental solution over $j$ steps is $F_j(\lambda)$. The query at step $s$, on the $\lambda'$ mode, equals $F_s(\lambda')/(s+1)$; this follows from
$r_{s+1}(\lambda')=(1-\lambda')F_s(\lambda')/(s+1)$, including $s=0$. An insertion therefore changes $v_{s+1}$ by $-F_s(\lambda')Z_s$, and propagation to $v_k$, followed by division by $k$, gives the kernel in \eqref{eq:nesterov-sharp-gain}. The bound $|U_j(t)|\le j+1$ on $[-1,1]$ gives $|F_j(\lambda)|\le j+1$, with equality for every $j$ at $\lambda=0$. Consequently
\[
 \bigl(A_k^{\rm N}\bigr)^2=\frac1{k^2}\sum_{s=0}^{k-1}(s+1)^2(k-s)^2
 =\frac{(k+1)^5-(k+1)}{30k^2}.
\]

For $k=2m+2$, the Jacobi factor in \eqref{eq:nesterov-factor} is $U_m(1-2\lambda)/(m+1)$, and its canonical matrix is $J_m^U$. Its leading and trailing polynomials are $2^{-j}U_j$. By \eqref{eq:general-drift-kernel}, the first $m$ block kernels are
\[
 K_i^{\rm B}(\lambda,\lambda')
 =-\frac{4(1-\lambda)^{k-m}}{m+1}U_{m-i-1}(1-2\lambda)U_i(1-2\lambda').
\]
For the subsequent $k-m=m+2$ Landweber steps,
\[
 K_{m+j}^{\rm B}(\lambda,\lambda')
 =-(1-\lambda)^{k-m-j-1}(1-\lambda')^j\frac{U_m(1-2\lambda')}{m+1},\qquad 0\le j<k-m.
\]
The same Chebyshev bound shows that each absolute value is maximized at $(0,0)$. The diagonal cofactor gives
\[
 d_i(J_m^U)=\frac{2(i+1)(m-i)}{m+1},\qquad
 \E(J_m^U)=\frac2{15}\left((m+1)^3-\frac1{m+1}\right).
\]
Summing the squared endpoint kernels gives \eqref{eq:block-sharp-gain}.
For every other $J\in\cJ(U_m)$, the spectral supremum is at least its endpoint channel $\sqrt{4\E(J)+m+2}$. For $k=2m+2$, \cref{thm:B} gives
\[
 \liminf_{m\to\infty}k^{-3}\inf_J A_k(\mathcal B_{k,J};[0,1])^2
 \ge\tfrac12 e(1/2,1/2)=\tfrac1{15}.
\]
The canonical block gives the reverse upper limiting bound by \eqref{eq:block-sharp-gain}. Taking square roots proves the optimal limit and the gain ratio.
\end{proof}

\begin{lemma}[The original endpoint gain for general momentum]\label{lem:nesterov-endpoint}
For $b>-1$, the Nesterov realization $\mathcal N_b$ of \eqref{eq:nesterov-iteration} satisfies
\[
 \lim_{k\to\infty}k^{-3}a_k(\mathcal N_b)^2=\frac{2}{3(b+3)(2b+3)}.
\]
\end{lemma}
\begin{proof}
At the zero spectral mode the unperturbed error and every query equal $e_0$. An insertion at step $j$ has coefficient $-w_{k,j}$, where
\[
 w_{k,j}=\sum_{r=j+1}^k\ \prod_{i=j+1}^{r-1}\frac{i-1}{i+b},\qquad
 w_{k,0}=1,\qquad 0\le w_{k,j}\le k-j.
\]
For $j/k\to t>0$, the product converges uniformly for $r/k\in[t,1]$ to $(t/(r/k))^{b+1}$: sum the logarithms and use the uniform $O(i^{-2})$ remainder. Hence $w_{k,j}/k\to f_b(t)=\int_t^1(t/u)^{b+1}\dd u$. For $b\ne0$, this is $t(1-t^b)/b$; for $b=0$, it is $t\log(1/t)$. The bound by $1-t$ controls the interval near zero in the Riemann sum. Therefore
\[
 \lim k^{-3}\sum_jw_{k,j}^2=\int_0^1f_b(t)^2\dd t
 =\frac{2}{3(b+3)(2b+3)},
\]
with the value at $b=0$ obtained by continuity or direct integration.
\end{proof}

\section{Uniform edge estimates}\label{app:edge}

This appendix strengthens the fixed-mode limit to the full range required for comparison of spectral measures. Write $\zeta_{n,k}=\cos\vartheta_{n,k}$ and set
\begin{equation}
 \gamma_\alpha=\tfrac\alpha2+\tfrac34,\quad
 \gamma_\beta=\tfrac\beta2+\tfrac34,\quad
 \bar\gamma=\tfrac12(\gamma_\alpha+\gamma_\beta),\quad
 N=n-1+2\bar\gamma,\quad h=\pi/N.
\end{equation}
Put $c_{0a}=1/4-a^2$ and
\begin{equation}
 \mathsf c_n=1-\frac{c_{0\alpha}+c_{0\beta}}{4N^2},\quad
 \mathsf d_n=\frac{c_{0\beta}-c_{0\alpha}}{4N^2},\quad
 x_k=(k-1+\gamma_\alpha)h,\quad \eta_{n,k}=\mathsf c_n\cos x_k+\mathsf d_n.
\end{equation}
Let $K_n=\lceil n/2\rceil$, $\hat d_k=\min(k,n+1-k)$ and $m_{k\ell}=\min(k+\ell,2n+2-k-\ell)$. A superscript $\vee$ swaps $\alpha,\beta$; reflection gives $\zeta_{n,n+1-k}=-\zeta_{n,k}^\vee$ and the same identity for $\eta$.

\begin{proposition}\label{prop:zeros}
For fixed $\alpha,\beta>-1$ there are $C,c>0$ such that, for all sufficiently large $n$:
\begin{enumerate}
\item for $k\le K_n$,
\[
 \left|\vartheta_{n,k}-x_k-\frac{c_{0\alpha}\cot(x_k/2)-c_{0\beta}\tan(x_k/2)}{4N^2}\right|
 \le C\left(\frac1{Nk^3}+\frac1{N^3}\right);
\]
\item $e_{n,k}:=|\zeta_{n,k}-\eta_{n,k}|\le C(N^{-2}\hat d_k^{-2}+N^{-4}\hat d_k)$;
\item $c\hat d_k/N\le\min(\vartheta_{n,k},\pi-\vartheta_{n,k})\le C\hat d_k/N$, $\vartheta_{n,k+1}-\vartheta_{n,k}\ge c/N$, and, for $k\ne\ell$,
\[
 \min\{|\zeta_{n,k}-\zeta_{n,\ell}|,|\eta_{n,k}-\eta_{n,\ell}|\}
 \ge c|k-\ell|m_{k\ell}/N^2;
\]
\item if $\mathsf c_n\cos\vartheta+\mathsf d_n$ lies between $\eta_{n,k}$ and $\zeta_{n,k}$, then $c\hat d_k/n\le\sin\vartheta\le C\hat d_k/n$.
\end{enumerate}
\end{proposition}
\begin{proof}
We give the argument for $k\le K_n$; reflection treats the other edge. The Jacobi Liouville form \cite[\S8.21]{Sz} is
\[
 u''+(N^2+c_{0\alpha}\vartheta^{-2}+r(\vartheta))u=0,
 \qquad u=\vartheta^{\alpha+1/2}g(\vartheta),
\]
where $g$ is even analytic near zero, $g(0)\ne0$, and $r\in C^1[0,3\pi/4]$ has primitive
\[
 \mathcal R(\vartheta)=c_{0\alpha}(\vartheta^{-1}-\tfrac12\cot(\vartheta/2))
                  +\tfrac12c_{0\beta}\tan(\vartheta/2).
\]
Use the Bessel basis $v_1=\sqrt\vartheta J_\alpha(N\vartheta)$, $v_2=\sqrt\vartheta Y_\alpha(N\vartheta)$, where $Y_\alpha$ is the Bessel function of the second kind; its Wronskian is $2/\pi$. Variation of constants gives $u=A_0v_1+(\pi/2)\int_0^\vartheta[v_1(\vartheta)v_2(s)-v_2(\vartheta)v_1(s)]r(s)u(s)\dd s$.

First choose $0<y_*<\min(1,j_{\alpha,1})/2$. Dividing this equation by $v_1$ on $(0,y_*/N]$ gives a Volterra kernel bounded by
\[
 Cs\,[1+(s/\vartheta)^{2\alpha}](1+|\log Ns|).
\]
Its integral is $O(N^{-2})$, uniformly in $\vartheta$, because $1+2\alpha>-1$. Thus $u/v_1=A_0(1+O(N^{-2}))$ with the sign of $A_0v_1$ throughout this interval. This estimate also initializes the phase perturbation with $\chi=O(N^{-2})$ up to $1/N$.

Write $u=Av_1+Bv_2$ with the variation-of-constants condition $A'v_1+B'v_2=0$. Set $J_\alpha=\mathsf M\cos\theta_\alpha$, $Y_\alpha=\mathsf M\sin\theta_\alpha$, and $(A,B)=\varrho(\cos\chi,\sin\chi)$. Then
\[
 u=\varrho\sqrt\vartheta\mathsf M(N\vartheta)\cos\Theta,
 \quad\Theta=\theta_\alpha(N\vartheta)-\chi,
 \quad\chi'=-\tfrac\pi2 r\vartheta\mathsf M(N\vartheta)^2\cos^2\Theta.
\]
The standard Hankel estimates \cite[10.17.9--12,10.18.17--18]{DLMF} give $x\mathsf M(x)^2=2\pi^{-1}(1+O(x^{-2}))$, with differentiated remainder $O(x^{-3})$. On $[1/N,3\pi/4]$, $\Theta'\ge cN$ and $|\Theta''|\le C[N^2(N\vartheta)^{-3}+1]$. Splitting $\cos^2\Theta=(1+\cos2\Theta)/2$ and integrating the oscillatory part by parts, with denominator $2\Theta'$, therefore gives
\begin{equation}
 \chi(\vartheta)=-\frac{\mathcal R(\vartheta)}{2N}+O(N^{-2}).
\end{equation}
All remainders are uniform on this interval, since the integrals of $(Ns)^{-2}$ and $(Ns)^{-3}$ from $1/N$ are $O(N^{-1})$.

The phase is strictly increasing from $y_*/N$. Its initial value lies in the same gap between $\pi/2+\pi\Z$ as $\theta_\alpha(y_*)$, and at $3\pi/4$ it has crossed more than $K_n$ levels. This fixes the numbering and gives the identity
\begin{equation}\label{eq:app-numer}
 \theta_\alpha(N\vartheta_{n,k})-\theta_\alpha(j_{\alpha,k})
 =\chi(\vartheta_{n,k})\qquad(k\le K_n).
\end{equation}
Using $\theta_\alpha'=1+O(k^{-2})$ between these points and McMahon's formula $j_{\alpha,k}=Nx_k+c_{0\alpha}/(2Nx_k)+O(k^{-3})$ gives (a). Finitely many small $k$ are absorbed in its $1/(Nk^3)$ error.

For (b), set $\delta=\vartheta_{n,k}-x_k$ and use $\sin x\cot(x/2)=1+\cos x$, $\sin x\tan(x/2)=1-\cos x$. The linear correction in $\cos(x_k+\delta)$ is exactly $\eta_{n,k}-\cos x_k$; the quadratic error and (a) give the claim. The adjacent-angle difference is $\pi/N+O((Nk^2)^{-1}+N^{-3})$ for large $k$; fixed small $k$ use the positive Bessel gaps in \eqref{eq:app-numer}. Reflection covers the middle and the other edge. The identity for a difference of cosines proves (c).

For (d), $1-\cos\vartheta=[(1-x)-c_{0\alpha}/(2N^2)]/\mathsf c_n$ is affine in $x$. It is comparable to $(k/N)^2$ at both ends of the segment. For fixed small $k$, positivity follows from $\sigma_k>4(\alpha+1)>c_{0\alpha}$; large $k$ follow from (c). Affinity gives the comparison throughout, and reflection finishes the proof.
\end{proof}

Frenzen--Wong \cite{FW85} give uniform zero estimates under $\alpha>-1/2$, $\alpha+\beta\ge-1$. The integrable endpoint Volterra bound above supplies the estimates throughout $\alpha,\beta>-1$.

\section{Reconstruction, convergence and geometry}

We prove the canonical limit in one route. The reference estimates below feed the half-interval reconstruction. A finite comparison chain then supplies a Gram approximation to the same transformation, with enough uniform control to pass to the energy. Throughout $\gamma=\gamma_\alpha$, $c_0=c_{0\alpha}$ and $\omega^2=z-c_0$.

\subsection{Reference waves and norms}\label{app:reference}

For analytic continuation we use the bilinear normalizations $\fn_k=\int_0^1 f(t,\tau_k)^2\dd t$. Let $B(a,b)=\Gamma(a)\Gamma(b)/\Gamma(a+b)$ be the beta function. Powers of positive real bases use their real logarithms. Define
\begin{equation}\label{eq:f}
 f(t,z)=\frac{(\sin\pi t/\pi)^{1-\gamma}}{B(\gamma,\gamma)}
 \int_0^t\left(\frac{\sin\pi s\sin\pi(t-s)}{\pi^2}\right)^{\gamma-1}
 \cos((t-2s)\omega)\dd s.
\end{equation}
\begin{lemma}[Reference data]\label{lem:f}
For $\Re\gamma>0$ this function is holomorphic in $(\gamma,\omega^2)$ and continuous down to $t=0$ on $[0,1/2]$. It solves the reference equation and has normalization $t^\gamma(1+O(t^2))$. For real $\gamma>1/4$ its reference eigenfunctions form a complete orthogonal system on $(0,1)$, with
\begin{equation}\label{eq:reference-norm}
 \begin{gathered}
 f(t,\tau_k)=\pi^{-\gamma}\sin^\gamma(\pi t)
       \frac{C_{k-1}^{(\gamma)}(\cos\pi t)}{C_{k-1}^{(\gamma)}(1)},\\
 \fn_k:=\int_0^1 f(t,\tau_k)^2\dd t
 =A(\gamma)\frac{\Gamma(k)}{(k-1+\gamma)\Gamma(k-1+2\gamma)},\qquad
 A(\gamma)=\frac{2^{1-2\gamma}\Gamma(2\gamma)^2}{\pi^{2\gamma}\Gamma(\gamma)^2}.
 \end{gathered}
\end{equation}
The norm formula continues as a nonvanishing holomorphic function for $\Re\gamma>0$, with $\fn_k\sim A(\gamma)k^{-2\gamma}$ as $k\to\infty$, uniformly on compact subsets. For $t\le1/2$, $|\Im\omega|\le b$,
\begin{equation}\label{eq:reference-bounds}
 |f(t,z)|\le C_b\min(t^{\Re\gamma},|\omega|^{-\Re\gamma}),\qquad
 |\partial_zf(t,z)|\le C_bt|\omega|^{-1}\min(t^{\Re\gamma},|\omega|^{-\Re\gamma})
\end{equation}
in the second inequality for $|\omega|\ge1$, uniformly on compact parameter sets. For real $\gamma$, also $|f(t,z)|\le f(t,c_0)e^{t|\Im\omega|}$ and $0\le f(t,c_0)\le Ct^\gamma$.
\end{lemma}
\begin{proof}
Mehler--Dirichlet \cite[14.12.1]{DLMF}, with order $1/2-\gamma$ and $\Re\gamma>0$, after $v=\pi(t-2s)$, identifies \eqref{eq:f} with
\[
 \pi^{-\gamma}\sin^\gamma(\pi t)\,
 {}_2F_1(\gamma-\omega/\pi,\gamma+\omega/\pi;
                  \gamma+1/2;\sin^2(\pi t/2)).
\]
The hypergeometric equation gives the ODE and the normalized branch. Gegenbauer orthogonality and its norm formula \cite[Tables~18.3.1,18.6.1]{DLMF} give \eqref{eq:reference-norm}; completeness follows by density of polynomials for the positive Gegenbauer weight. The identity theorem continues the formula, while the uniform gamma quotient asymptotic gives its bounds and nonvanishing.

For the uniform frequency estimates, put $\operatorname{sinc}z=\sin z/z$, continued by one at zero, and rescale $s=t\varsigma$:
\begin{equation}\label{eq:app-fresc}
 f(t,z)=\frac{t^\gamma\operatorname{sinc}(\pi t)^{1-\gamma}}{B(\gamma,\gamma)}
 \int_0^1[\varsigma(1-\varsigma)]^{\gamma-1}\Psi_t(\varsigma)
              \cos(t\omega(1-2\varsigma))\dd\varsigma,
\end{equation}
where $\Psi_t=[\operatorname{sinc}(\pi t\varsigma)\operatorname{sinc}(\pi t(1-\varsigma))]^{\gamma-1}$. Its positive base lies in $[4/\pi^2,1]$, so all $C^m$ bounds are uniform on compact parameter sets. Power-log majorants give joint holomorphy and continuity, including zero.

Partition the integral near its two endpoints. Each oscillatory piece has the form
\[
 \int_0^1\varsigma^{\gamma-1}g(\varsigma)e^{\pm2\rI\Omega\varsigma}\dd\varsigma,
\]
with $g$ smooth and zero near one. Splitting at $|\Omega|^{-1}$ and integrating by parts $m>\max\Re\gamma$ times on the rest gives $O(|\Omega|^{-\Re\gamma})$, uniformly in a bounded imaginary strip. The direct integral bound handles $|\Omega|\le1$. Differentiating in $\omega$ inserts a smooth factor of size $t$; division by $2\omega$ gives \eqref{eq:reference-bounds}. The remaining real-parameter assertions follow from the positive integrand. On compact $z,z'$ sets the same rescaled formula also gives $f(t,z)-f(t,z')=O(t^{\gamma+2}|z-z'|)$.
\end{proof}

The normalized branch is unique. If $f_0=t^\gamma(1+O(t^2))$, reduction of order gives the second solution $f_0(t)\int_t^\delta f_0(s)^{-2}\dd s$. Its independent term has order $t^{1-\gamma}$, or $t^{1/2}\log(1/t)$ when $\gamma=1/2$. A solution $O(t^{\gamma+2})$ is therefore zero, including $\gamma<1/2$. Comparing a normalized solution with $\pi^{-\gamma}\sin^\gamma\pi t$ in Green's identity gives $f_0'=\gamma t^{\gamma-1}(1+O(t^2))$ and the Wronskian of two differently energized normalized waves is $O(t^{2\gamma+1})$. General solutions are $O(t^{1-\gamma})$ for $\gamma>1/2$, $O(t^{1/2}\log(1/t))$ at equality and $O(t^\gamma)$ below it. For $\gamma\ge1/2$, every solution $O(t^\gamma)$ is a multiple of the normalized branch. The normalization $t^\gamma(1+O(t^2))$ selects the branch in all three regimes.

Two estimates used at negative energy also follow directly from \eqref{eq:app-fresc}:
\[
 0<f(t,c_0-y^2)\le Ce^{ty}\min(t^\gamma,y^{-\gamma}),\qquad
 f(t,c_0-y^2)\ge c_ae^{ty}y^{-\gamma}\quad(a\le t\le1/2, y\ge y_a).
\]
For the lower bound retain $s\le1/y$ in \eqref{eq:f}. The same upper bound with $y=|\omega|$ and $e^{t|\Im\omega|}$ holds for arbitrary complex $\omega$. Denote the integral in \eqref{eq:app-fresc}, with $\Omega=t\omega$, by $I_t(\Omega)$. This even entire function has the polynomial decay on the real axis just proved and the Laplace bound on the imaginary axis; Phragm\'en--Lindel\"of \cite{Tit} applies to $(1+\Omega)^\gamma e^{\rI\Omega}I_t(\Omega)$ in the first quadrant and to its reflected versions in the others.

We will need the midpoint constants. With $a_0=f(1/2,z)$, $b_0=f_t(1/2,z)$, $\lambda=\omega/\pi$ and $\kappa=2\pi^2A(\gamma)$, Gauss's second summation and duplication \cite[15.4.28,5.5.5]{DLMF} give
\begin{equation}
 \begin{gathered}
 a_0=\frac{\pi^{1/2-\gamma}\Gamma(\gamma+1/2)}
 {\Gamma((\gamma+1-\lambda)/2)\Gamma((\gamma+1+\lambda)/2)},\qquad
 b_0=\frac{2\pi^{3/2-\gamma}\Gamma(\gamma+1/2)}
 {\Gamma((\gamma-\lambda)/2)\Gamma((\gamma+\lambda)/2)},\\
 2a_0b_0=\kappa\Delta_\tau(z),\qquad
 \Delta_\tau(z)=\frac1{\Gamma(\gamma-\lambda)\Gamma(\gamma+\lambda)}.
 \end{gathered}
\end{equation}
The zeros of $a_0$ are the even $\tau_k$ and those of $b_0$ the odd ones. At $z=c_0-y^2$, $a_0>0$ and $b_0/a_0=y+O(y^{-1})$ by the gamma quotient expansion.

\subsection{Reconstruction and midpoint closure}\label{app:closure}

\begin{proposition}[Spectral reconstruction]\label{prop:reconstruction}
Let $\alpha>-1$, with $\gamma,q_\alpha,\tau_k,\sigma_k,c_k$ and $F$ as in \cref{sec:contpersym}. Then the series for $F$ converges uniformly on $[0,1/2]^2$, and its integral operators satisfy $I+\cF_t\ge\kappa_*I$ on $L^2((0,t),\dd s)$ for one $\kappa_*>0$, $0<t\le1/2$. The solution $K$ of \eqref{eq:GL} is continuous, and $w(t)=K(t,t)$ belongs to $C^1(0,1/2]$, with $w(0)=w(1/2)=0$.
The transformed waves solve $-u''+(q_\alpha+2w')u=zu$ in the branch $u=t^\gamma(1+O(t))$. Their midpoint conditions are Neumann at odd $\sigma_k$ and Dirichlet at even $\sigma_k$. Reflection with sign $(-1)^{k-1}$ makes $\sqrt{c_k}u(\cdot,\sigma_k)$ a complete orthonormal basis of $L^2((0,1),\dd t)$. The resulting spectral operator is the positive self-adjoint $H_\alpha$ of \cref{thm:A}, and its Green series converges uniformly on $[0,1]^2$.
\end{proposition}
\begin{proof}\refstepcounter{prooflink}\label{proof:prop:reconstruction}

We first justify the spectral normalizations in \eqref{eq:norming-main}. Let $\delta_k=\sigma_k-\tau_k=O(k^{-2})$. Both spectra have gaps bounded below by $c|\ell^2-k^2|$. The product estimate we use repeatedly is
\begin{equation}\label{eq:app-rhohyp}
 \sum_{\ell\ne k}\frac{k^{-2}+\ell^{-2}}{|\ell^2-k^2|}=O(k^{-2}).
\end{equation}
Split at $k/2$ and $2k$: the three contributions are $O(k^{-2})$, $O(k^{-3}\log k)$ and $O(k^{-3})$. Applying $\log(1+x)=x+O(x^2)$ gives absolute convergence, positivity and $\rho_k:=c_k/c_k^0=1+O(k^{-2})$.

For exact normalization put
\begin{equation}
 \Delta_\sigma(z)=\frac{\pi^\alpha}{\sqrt2}z^{-\alpha/2}J_\alpha(\sqrt z),
 \qquad
 \frac{\Delta_\tau(z)}{\Delta_\sigma(z)}
       =\prod_{\ell\ge1}\frac{\tau_\ell-z}{\sigma_\ell-z}.
\end{equation}
The entire Bessel expression is understood by its power series. The gamma and Bessel products give the second identity up to a constant. At $z=-y^2$ both characteristic functions are asymptotic to $\pi^\alpha y^{-\alpha-1/2}e^y/(2\sqrt\pi)$, whereas the logarithm of the product is $O(y^{-2})$; the constant is one. Differentiating the products at their zeros therefore gives
\begin{equation}
 \rho_k=\frac{\Delta_\tau'(\tau_k)}{\Delta_\sigma'(\sigma_k)},\quad
 \fn_k=\kappa|\Delta_\tau'(\tau_k)|,\quad
 c_k=\frac1{\kappa|\Delta_\sigma'(\sigma_k)|}
 =\mathfrak w_\alpha\frac{j_{\alpha,k}^{\alpha+1}}{|J_{\alpha+1}(j_{\alpha,k})|},
\end{equation}
where $\mathfrak w_\alpha=2^{\alpha+1}\Gamma(\gamma)^2/(\sqrt\pi\Gamma(\alpha+3/2)^2)$. These follow by $J_\alpha'(j)=-J_{\alpha+1}(j)$ and the simple gamma zeros. In signed form $c_k=(-1)^k/(\kappa\Delta_\sigma'(\sigma_k))$, and similarly for $c_k^0$.

\paragraph{Coercivity and the transformation.}
Write $F_k$ for the $k$th summand of \eqref{eq:F} and $\cF_t$ for its sum as an integral operator on $L^2((0,t),\dd s)$. The reference bounds show $\sup|F_k|\le Ck^{-2}$: the norming perturbation is $O(k^{-2})$ and the wave perturbation is $O(k^{-3})$ after normalization. For $v\in L^2((0,t),\dd s)$, extended by zero, Parseval gives
\[
 \ip{(I+\cF_t)v}{v}_{L^2((0,t),\dd s)}=\sum_k c_k|V(\sigma_k)|^2,\qquad
 V(z)=\int_0^t v(s)f(s,z)\dd s.
\]
If the form vanishes, the even entire function $V(c_0+\omega^2)$ vanishes at $\pm\sqrt{\sigma_k-c_0}$, whose counting function is at least $2r/\pi-C$. Its growth is bounded by $Ce^{t|\omega|}$. Jensen's formula \cite{Lev} would give $2R/\pi-O(\log R)\le tR+O(\log R)$, impossible for $t\le1/2$. Thus $V\equiv0$, and reference completeness gives $v=0$. Compactness of $\cF_{1/2}$ now gives $I+\cF_t\ge\kappa_*>0$ uniformly for all half-interval sections.

Consequently $K(t,\cdot)=-(I+\cF_t)^{-1}F(t,\cdot)$ is the unique solution of \eqref{eq:GL}. To obtain joint continuity, extend the sections by zero to $(0,1/2)$: their operators and right-hand sides vary continuously in Hilbert--Schmidt and $L^2$ norm. The equation then upgrades $L^2$ continuity to continuity on the closed triangle. This also proves stability of the sections and of $K$ under uniform kernel perturbations.

\paragraph{Why each parity closes.}
There is an exact splitting $\sum_{k\text{ odd}}F_k=\sum_{k\text{ even}}F_k=F/2$. To see it, integrate
\[
 f(t,z)f(s,z)(\Delta_\sigma(z)^{-1}-\Delta_\tau(z)^{-1})
\]
on $|z-c_0|=R_m^2$, $R_m=\pi(m+\gamma-1/2)$. The enclosed residues are $-\kappa\sum_{k\le m}(-1)^{k-1}F_k$. On these contours the global wave bound, reflection formula for $\Gamma$, and \eqref{eq:app-rhohyp} give respectively
\[
 |f(t,z)f(s,z)|\le CR_m^{-2\gamma}e^{(t+s)|\Im\omega|},\quad
 |\Delta_\tau|\ge cR_m^{1-2\gamma}e^{|\Im\omega|},\quad
 |\Delta_\sigma/\Delta_\tau-1|\le CR_m^{-2}.
\]
Their product is $O(R_m^{-3})$ for $t+s\le1$, so the contour integral is $O(R_m^{-1})$ and vanishes.

Let $\mathcal K$ be the Volterra operator with kernel $K$. Equation~\eqref{eq:GL}, followed by symmetry, gives
$(I+\mathcal K)(I+\cF_{1/2})(I+\mathcal K)^*=I$: its difference from the identity has zero kernel below the diagonal and is self-adjoint. Hence each system
\begin{equation}
 \{\sqrt{2c_k}\,u(\cdot,\sigma_k):k\text{ odd}\},\qquad
 \{\sqrt{2c_k}\,u(\cdot,\sigma_k):k\text{ even}\}
\end{equation}
is a Parseval frame. It is actually an orthonormal basis: before transformation each system is quadratically close to its reference parity basis, with $O(k^{-2})$ norm differences. Identify $\ell^2$ with the reference parity basis, and let $T$ synthesize the perturbed waves before transformation. Then $T-I$ is Hilbert--Schmidt, so $T$ is Fredholm of index zero. The transformed synthesis $S=(I+\mathcal K)T$ has the same index, and the frame identity gives $SS^*=I$. Thus $S$ is surjective, index zero makes it injective, and $SS^*=I$ makes it unitary. This proves both completeness and orthonormality.

\paragraph{The differential equation.}
To differentiate the transformation, first isolate the Green part of $F$. Let $\mathsf s=\sum_k\delta_k$ and $\varpi_k=\sqrt{\tau_k-c_0}$. Refining the same three-range sum gives
\[
 \rho_k=1+\mathsf s\varpi_k^{-2}+O(k^{-3}\log k),\qquad
 F=\mathsf s\mathfrak G_{c_0}+F^\flat,
\]
where $\mathfrak G_{c_0}(t,s)=\chi(t\wedge s)\chi(1-(t\vee s))/\mathcal W_0$, $\chi=f(\cdot,c_0)$, $\mathcal W_0=\kappa/\Gamma(\gamma)^2>0$. Green's identity and reference completeness identify this with its spectral series. Away from $t=0$, the differentiated remainder terms are $O(k^{-2}\log k)$ and therefore summable. Indeed the ODE and a difference-quotient estimate give $|f_t|\le C_\delta k^{1-\gamma}$ and $|f_{tz}|\le C_\delta k^{-\gamma}$.

It follows by difference quotients in GL that $K_t$ is continuous on $\delta\le t\le1/2$, $0\le s\le t$, and solves
\[
 (I+\cF_t)K_t(t,\cdot)=-\dot F(t,\cdot)-K(t,t)F(t,\cdot).
\]
Here $\dot F$ is the one-sided derivative from the appropriate side of the Green kink. Differentiating in $s$ under the integral is valid by its bounded piecewise-continuous derivative. Thus $w(t)=K(t,t)$ is $C^1$ on $(0,1/2]$, $w'=K_t(t,t)+K_s(t,t)$, and $w(0+)=0$.

We justify the differential equation through finite spectral sums.
Put $F^{(N)}=\sum_{k\le N}F_k$. Uniform kernel convergence gives section inverses bounded by $2/\kappa_*$ for all large $N$ and $K^{(N)}\to K$ uniformly. Write $w^{(N)}(t)=K^{(N)}(t,t)$. For a finite sum, write $F^{(N)}(t,s)=h(t)^\top C h(s)$, where the entries of $h$ are distinct reference waves and $C$ is diagonal. Then
\[
 K^{(N)}(t,s)=h(s)^\top a(t),\qquad
 a(t)=-\left(I+C\int_0^t h(r)h(r)^\top\dd r\right)^{-1}Ch(t).
\]
The matrix is invertible by section coercivity. Apply the reference equation in the two variables to the finite GL equation, first integrating over $[\varepsilon,t]$. Combining terms into Wronskians before sending $\varepsilon\downarrow0$ gives
\[
 \Theta(t,s)+\int_0^t\Theta(t,r)F^{(N)}(r,s)\dd r=0,\qquad
 \Theta=(\partial_t^2-q_\alpha(t)-2(w^{(N)})'-\partial_s^2+q_\alpha(s))K^{(N)}.
\]
For each fixed $N$, the finite formula gives $u^{(N)}-f=O(t^{3\gamma+1})$ and its derivative $O(t^{3\gamma})$. The reference Wronskians are $O(t^{2\gamma+1})$, and the transformed ones are $O(t^{2\gamma+1}+t^{4\gamma})$. Thus every boundary term in these Green identities vanishes at zero. Coercivity gives $\Theta=0$. A second Green identity in the transformation for $u^{(N)}$ yields
\[
 -(u^{(N)})''+(q_\alpha+2(w^{(N)})')u^{(N)}=zu^{(N)},\qquad
 [(u^{(N)})'-2w^{(N)}u^{(N)}]'=(q_\alpha-z)u^{(N)}-2w^{(N)}(u^{(N)})'.
\]

For the limit, work on the fixed space $L^2((0,1/2),\dd s)$ and write $\hat g(z)=\int_0^{1/2}g(s)f(s,z)\dd s$. The quadratic-form identity above gives $\sum_k c_k|\hat g(\sigma_k)|^2\le C\|g\|_2^2$. The reference derivative estimates imply, for $t\ge\delta$,
\[
 |\langle\partial_tF_k(t,\cdot),g\rangle|
 \le C_\delta\left(k^{-1}\sqrt{c_k}|\hat g(\sigma_k)|+k^{-2}\|g\|_2\right).
\]
For clarity, the leading differentiated term is $(\rho_k-1)h_k'(t)h_k(s)$, where $h_k=f(\cdot,\tau_k)/\sqrt{\fn_k}$. The displacement $\sqrt{c_k}[f(\cdot,\sigma_k)-f(\cdot,\tau_k)]$ is $O(k^{-3})$ and its $t$ derivative $O(k^{-2})$ away from zero, so the other terms have $L^2$ norm $O(k^{-2})$. The Bessel bounds for both reference and perturbed waves and Cauchy--Schwarz therefore give
\[
 \sup_{t\in[\delta,1/2]}
 \left\|\sum_{k>N}\partial_tF_k(t,\cdot)\right\|_{L^2((0,1/2),\dd s)}
 \le C_\delta\left[\left(\sum_{k>N}k^{-2}\right)^{1/2}
                         +\sum_{k>N}k^{-2}\right]\longrightarrow0.
\]
Uniform convergence of $F^{(N)}$ identifies this strong row limit with $\dot F$. Restricting the rows to $(0,t)$, the differentiated GL equations, uniform convergence of their kernels and the common inverse bound then give $(K^{(N)})_t\to K_t$ in these $L^2$ sections, uniformly in $t$. Hence
\[
 (u^{(N)})'=f_t+w^{(N)}f+\int_0^t(K^{(N)})_t(t,s)f(s,z)\dd s
\]
converges uniformly on $[\delta,1/2]$ to $u'$. Passing to the integrated quasi-derivative identity gives
\[
 [u'-2wu]_{t_0}^t=\int_{t_0}^t\bigl((q_\alpha(r)-z)u(r,z)-2w(r)u_t(r,z)\bigr)\dd r.
\]
Since $w\in C^1$ locally, $u\in C^2$ and $-u''+(q_\alpha+2w')u=zu$. The integral transformation gives $u=t^\gamma(1+O(t))$. The same limit in the finite-sum Green identity yields
\begin{equation}\label{eq:app-CD}
 (z-z')\int_0^x u(t,z)u(t,z')\dd t
 =u(x,z)u_t(x,z')-u_t(x,z)u(x,z').
\end{equation}
In particular $a(z)=u(1/2,z)$ and $b(z)=u_t(1/2,z)$ are entire of order at most $1/2$, with bound $C(1+|\omega|)e^{|\Im\omega|/2}$.

\paragraph{Identification of the midpoint conditions.}
Orthogonality within either parity in \eqref{eq:app-CD} gives one fixed real boundary angle for that class: $S=a\sin\theta-b\cos\theta$ vanishes on its levels. The endpoint vector $(a(z),b(z))$ is nonzero by uniqueness of the regular Cauchy problem at $1/2$. If $S(z)=0$ at another point, \eqref{eq:app-CD} makes $u(\cdot,z)$ orthogonal to that entire parity basis, forcing $u(\cdot,z)=0$, contrary to its endpoint normalization. Thus these are all zeros of $S$. Differentiating \eqref{eq:app-CD} at a level gives the nonzero squared half-interval norm; since $(a,b)\ne0$, each zero is simple.

At $z=c_0-y^2$, put $\mathcal I_y=\int_0^{1/2}f(s,z)/f(1/2,z)\dd s$. The Laplace asymptotic in \eqref{eq:f} gives
\[
 \frac{f(1/2-v/y,c_0-y^2)}{f(1/2,c_0-y^2)}\longrightarrow e^{-v}
\]
locally uniformly for $v\ge0$. On $1/4\le s\le1/2$ the negative-energy bounds dominate this ratio by $Ce^{-v}$, while $s<1/4$ contributes $O(e^{-y/4})$. Hence $\mathcal I_y\sim y^{-1}$ as $y\to\infty$, and $(b_0/a_0)\mathcal I_y\to1$. Thus
\[
 a/a_0=1+w(1/2)\mathcal I_y+o(\mathcal I_y),\qquad
 b/a_0=b_0/a_0+w(1/2)+O(y^{-1}),\qquad b/a=y+o(1).
\]
Hadamard products and $\delta_k=O(k^{-2})$ give $S_{\rm odd}/b_0=C_1(1+O(y^{-2}))$, $S_{\rm even}/a_0=C_2(1+O(y^{-2}))$. Their ratio is $(C_1/C_2)y+O(y^{-1})$. Since $b/a=y+o(1)$, a boundary form with a nonzero $b$ coefficient has order $ya$, while a pure $a$ form has order $a$. The ratio thus forces the even form to be $a$ and the odd form to be proportional to $b-ha$. Comparing its constant term with $Cy+O(y^{-1})$ gives $h=0$. Odd levels are therefore Neumann and even levels Dirichlet. The normalized product for the even form gives $a/a_0=1+O(y^{-2})$; its earlier expansion, with $\mathcal I_y\sim y^{-1}$, yields $w(1/2)=0$. Normalizing the product at negative infinity then gives
\[
 2u(1/2,z)u_t(1/2,z)=\kappa\Delta_\sigma(z).
\]

Reflect $V=q_\alpha+2w'$ and put $\upsilon_k=\sqrt{c_k}u(\cdot,\sigma_k)$ on the first half, extended with sign $(-1)^{k-1}$. The conditions just proved make these waves $C^2$ across the midpoint. The two half-bases give a full orthonormal basis, and $\int_0^1(V(t)-q_\alpha(t))\dd t=4[w(1/2)-w(0+)]=0$ as an improper integral. Define
\[
 D(H_\alpha)=\left\{v\in L^2((0,1),\dd t):\sum_k\sigma_k^2|\ip{v}{\upsilon_k}_{L^2((0,1),\dd t)}|^2<\infty\right\},\qquad
 H_\alpha v=\sum_k\sigma_k\ip{v}{\upsilon_k}_{L^2((0,1),\dd t)}\upsilon_k.
\]
This is positive self-adjoint. For $v\in C_c^\infty(0,1)$, integration by parts gives $\sigma_k\ip{v}{\upsilon_k}_{L^2((0,1),\dd t)}=\langle-v''+V_\alpha v,\upsilon_k\rangle$. Parseval therefore puts $v$ in $D(H_\alpha)$ with $H_\alpha v=-v''+V_\alpha v$. Testing finite spectral sums then gives this differential action distributionally on all of $D(H_\alpha)$. Since $\sup_{k,t}|\upsilon_k(t)|<\infty$ and $\sum_k\sigma_k^{-1}<\infty$, its Green series converges uniformly on $[0,1]^2$. This proves the operator assertions. At $\alpha=1/2$, $\gamma=1$, $c_0=0$ and $\sigma_k=\tau_k=\pi^2k^2$. Thus $F=K=0$ and $H_{1/2}$ is the Dirichlet Laplacian, giving the explicit $U$ profile in the main text.
\end{proof}

\subsection{A comparison chain with uniform control}\label{app:aw}

Use the two-edge notation of \cref{app:edge}. Define
\begin{equation}\label{eq:chain-coefficients}
 \begin{aligned}
 A_m&=\frac{\sin((m+\gamma_\alpha+1/2)h)\sin((m+2\bar\gamma)h)}
 {\sin((m+\bar\gamma)h)\sin((m+\bar\gamma+1/2)h)},\\
 C_m&=\frac{\sin(mh)\sin((m+\gamma_\beta-1/2)h)}
 {\sin((m+\bar\gamma-1/2)h)\sin((m+\bar\gamma)h)},\\
 a_m&=\tfrac12\sqrt{A_{m-1}C_m},\qquad
 b_m=\cos(\gamma_\alpha h)-\tfrac12(A_m+C_m),
 \end{aligned}
\end{equation}
where $A_{n-1}=C_0=0$; the quotient formulas define $A_m$ for $0\le m<n-1$ and $C_m$ for $1\le m<n$. Let $J_n^\natural$ have these entries, $J_n^0=\mathsf c_nJ_n^\natural+\mathsf d_nI$, and $\psi_m=\psi_m^\natural$, $\psi_{n,m}^0(x)=\psi_m((x-\mathsf d_n)/\mathsf c_n)$ denote their orthonormal polynomials.

All relevant sine arguments lie in $(0,\pi)$; thus the couplings are positive. Reflection gives $A_{n-1-m}=C_m$, so the chains are persymmetric. Their spectra are $\{\cos x_k\}$ and $\{\eta_{n,k}\}$. For the spectral identity, let $D$ have $n-1$ rows, with $D_{m,m}=\sqrt{A_m/2}$, $D_{m,m+1}=-\sqrt{C_{m+1}/2}$. The sine addition formulas give
\begin{equation}\label{eq:darboux}
 D^\top D=y_1I-J_n^\natural(\alpha,\beta),\qquad
 DD^\top=y_1I-J_{n-1}^\natural(\alpha+2,\beta),\qquad
 y_1=\cos(\gamma_\alpha h).
\end{equation}
The second chain has the same $N,h$, while its right exponent is $\gamma_\alpha+1$. The nonzero spectra of the two Gram matrices coincide. Induction removes exactly the first cosine node and proves the full spectrum, including exceptional algebraic parameters. Persymmetric uniqueness then gives $J_n^\natural(\beta,\alpha)=-\diag((-1)^m)J_n^\natural(\alpha,\beta)\diag((-1)^m)$, and the corresponding reflected weights and polynomials.

The chain is the persymmetric singular truncation of the Askey--Wilson recurrence with parameters $(a,as,c,cs)$, $a=e^{\rI\gamma_\alpha h}$, $c=-e^{\rI\gamma_\beta h}$, $s=e^{\rI h}$ and $q=s^2$ \cite{AW85,KLS,LVZ18}. Cancel common factors before evaluating the recurrence at the truncation; the resulting coefficients are exactly \eqref{eq:chain-coefficients}. The $q$-difference identity \cite[18.28.6\_1--2]{DLMF} continues in the same way wherever its evaluated denominators are nonzero. Symmetric parameters give the $q$-ultraspherical chain; when $2\gamma_\alpha\in\N$, $q$ is a root of unity, as in \cite{SZ96,VZ24}. These identifications supply algebraic context; \eqref{eq:darboux} proves our spectrum throughout the parameter range.

\paragraph{Products and the top wave.}
At $y_1$, the recurrence has a constant solution after leading-coefficient rescaling. Consequently
\begin{equation}
 \psi_m(y_1)^2=\frac{\prod_{j<m}A_j}{\prod_{j=1}^m C_j}
 =\frac{\sin((m+\bar\gamma)h)}{\sin(\bar\gamma h)}
   \frac{p_m(\gamma_\alpha+1/2)}{p_m(\gamma_\beta+1/2)}p_m(2\bar\gamma),
 \quad p_r(a)=\prod_{\ell=0}^{r-1}\frac{\sin((\ell+a)h)}{\sin((\ell+1)h)}.
\end{equation}
For fixed $a>0$, uniformly over $0\le r\le n/2+3$ as $n\to\infty$,
\[
 p_r(a)\sim\frac{\Gamma(r+a)}{\Gamma(a)\Gamma(r+1)}
                    \operatorname{sinc}(rh)^{a-1},\qquad
 c(r+1)^{a-1}\le p_r(a)\le C(r+1)^{a-1}.
\]
Indeed, separating the gamma quotient leaves the exponential of
$\sum_{\ell<r}[g((\ell+a)h)-g((\ell+1)h)]$, where $g(u)=\log\operatorname{sinc}u$. Taylor's formula and a Riemann sum give $(a-1)g(rh)+O(h)$, uniformly in this range. The gamma quotient gives the power bounds and
\[
 n^{1-a}p_{r_n}(a)\longrightarrow\Gamma(a)^{-1}(\sin\pi t/\pi)^{a-1}
 \quad(r_n/n\to t\in(0,1/2]),
\]
uniformly away from zero. The same proof gives the power bounds with step $h/2$ and $r\le n+3$.

Write $\Delta\psi_m=\psi_{m+1}-\psi_m$. It follows that
$c(m+1)^\gamma\le\psi_m(y_1)\le C(m+1)^\gamma$,
$|\Delta\psi_m(y_1)|\le C(m+1)^{\gamma-1}$, and
\[
 n^{-\gamma}\psi_m(y_1)\longrightarrow U(m/n),\quad
 U(t)=\kappa_U\sin^\gamma\pi t,\quad \kappa_U>0,\quad
 \kappa_U^2=\frac{\Gamma(\gamma_\beta+1/2)}
 {\bar\gamma\pi^{2\gamma}\Gamma(\gamma+1/2)\Gamma(2\bar\gamma)},
\]
uniformly for $m\le n/2$. The small range $m\le\delta n$ is bounded by $C(\delta+1/n)^\gamma$. Taylor expansion of $\log\sin((m+p)h)$ gives, on $1\le m\le n/2+3$, $a_m\ge c$ and $|a_m-1/2|+|b_m|\le Cm^{-2}$; the $m^{-1}$ terms cancel. Also $|b_0|\le C$.

\paragraph{One Volterra equation and its envelope.}
Reduction of order in the three-term recurrence gives
\begin{equation}\label{eq:volt}
 \begin{gathered}
 \psi_m(y)=\psi_m(y_1)+(y-y_1)\sum_{j<m}G(m,j)\psi_j(y),\\
 G(m,j)=\psi_m(y_1)\psi_j(y_1)
       \sum_{\ell=j}^{m-1}\frac1{a_{\ell+1}\psi_\ell(y_1)\psi_{\ell+1}(y_1)}.
 \end{gathered}
\end{equation}
Set $G(m,m)=0$ and $(\cG v)_m=\sum_{j<m}G(m,j)v_j$. The kernel is nonnegative and
\begin{equation}\label{eq:Gbd}
 G(m,j)\le C_G(m+1)^\gamma(j+1)^\gamma
                 \sum_{\ell=j}^{m-1}(\ell+1)^{-2\gamma}.
\end{equation}
Let $u^{(k)}=\cG^k\psi(y_1)$; then $u_m^{(k)}=0$ for $m<k$. Interchanging the sums in \eqref{eq:Gbd} uses
\[
 \sum_{\ell=k}^{m-1}(\ell+1)^{-2\gamma}
       \sum_{j=k}^{\ell}(j+1)^{2\gamma+2k}
 \le\frac{3(2\gamma+3)(m+1)^{2k+2}}{(2k+1)(2k+2)}.
\]
This follows by comparison of each power sum with its integral and endpoint term. Induction with $B=3C_G(2\gamma+3)$ gives
$u_m^{(k)}\le CB^k(m+1)^{\gamma+2k}/(2k)!$.
For first differences use
\[
 G(m+1,j)-G(m,j)=\frac{\Delta\psi_m(y_1)}{\psi_m(y_1)}G(m,j)
            +\frac{\psi_j(y_1)}{a_{m+1}\psi_m(y_1)}.
\]
The same summation gives the additional factor $2k/(m+1)$ for $k\ge1$. Thus, for $\varepsilon=|y-y_1|(m+1)^2$, both $|\psi_m(y)|$ and $\varepsilon^{-1}|\psi_m(y)-\psi_m(y_1)|$ are bounded by $C(m+1)^\gamma\cosh\sqrt{B\varepsilon}$; the difference of their first differences is bounded by $C\varepsilon(m+1)^{\gamma-1}\cosh\sqrt{B\varepsilon}$.

For $|\Re\theta|\le4\pi/5$, $|\Im\theta|\le b/N$, this implies the uniform envelope
\begin{equation}\label{eq:envelope}
 |\psi_m(\cos\theta)|\le C_b\min((m+1)^\gamma,|\theta|^{-\gamma}),
 \qquad m\le n/2+1.
\end{equation}
For completeness, if $|\theta|=O(N^{-1})$ use the preceding Volterra bound. Otherwise start at $m_0=\lfloor|\theta|^{-1}\rfloor$ and set $v_m^\pm=\psi_{m+1}-e^{\pm\rI\theta}\psi_m$ and $|v_m|=|v_m^+|+|v_m^-|$. Initially $|v_{m_0}^+|+|v_{m_0}^-|\le C|\theta|^{1-\gamma}$. The recurrence and the coefficient estimate give
\[
 |v_{m+1}|\le e^{|\Im\theta|}
 \left(1+\frac{C}{|\theta|(m+2)^2}\right)|v_m|.
\]
Its product is uniformly bounded up to $n/2$. Since $|\sin\theta|\ge c|\theta|$, recovering $\psi_m$ proves \eqref{eq:envelope}. Cauchy's formula on radius $1/N$ gives $|\partial_\theta\psi_m(\cos\theta)|\le C_bN\min((m+1)^\gamma,|\theta|^{-\gamma})$ for $|\Re\theta|\le3\pi/4$. At zero the original Volterra bound gives $|\psi_m(1)|\le C(m+1)^\gamma$.

\paragraph{The limit of the comparison waves.}
The continuum version of \eqref{eq:volt} is
\[
 V=U-(z-\gamma^2\pi^2)\int_0^t\mathsf G(t,s)V(s)\dd s,
 \qquad \mathsf G(t,s)=U(t)U(s)\int_s^t U(r)^{-2}\dd r.
\]
It has a unique solution in the literal class $V=O(t^\gamma)$. Iteration uses
\[
 \int_0^t s^{2\gamma+2k}\int_s^t r^{-2\gamma}\dd r\dd s
 =\frac{t^{2k+2}}{(2k+2)(2\gamma+2k+1)}
\]
and gives a factorially convergent series and uniqueness. Differentiating on $(0,1/2]$ identifies the solution as $\kappa_U\pi^\gamma f(t,z+c_0)$: the residual at the ground energy is $O(t^{\gamma+2})$, hence zero by the reference branch argument. This is valid also for $\gamma<1/2$.

For $z_n\to z$, put $y_n=1-z_n/(2N^2)$, $z_{1,n}=2N^2(1-y_1)$, $V_{n,m}=n^{-\gamma}\psi_m(y_n)$ and $U_{n,m}=n^{-\gamma}\psi_m(y_1)$. Then \eqref{eq:volt} reads
\[
 V_{n,m}=U_{n,m}-\frac{n^2}{2N^2}(z_n-z_{1,n})\frac1n\sum_{j<m}\frac{G(m,j)}n V_{n,j}.
\]
Here $z_{1,n}\to\gamma^2\pi^2$ and $G(m,j)/n\to2\mathsf G(t,s)$ for $m/n\to t$, $j/n\to s>0$. Equation~\eqref{eq:Gbd} gives the uniform weighted bound
\[
 \frac{G(m,j)}n\left(\frac{j+1}{m+1}\right)^\gamma\le C.
\]
In the norm $\|V\|_{\gamma,n}=\max_m|V_m|/((m+1)/n)^\gamma$, this makes the discarded $j<\delta n$ part $O(\delta)$. On $s\ge\delta$ the products converge uniformly, so their Riemann sums have an error $\epsilon_{n,\delta}\to0$. Subtract the two Volterra equations and iterate the positive bound above. The discrete Gronwall series gives
\[
 \max_{\delta n\le m\le n/2}
 \frac{|V_{n,m}-V(m/n)|}{((m+1)/n)^\gamma}
 \le C(\epsilon_{n,\delta}+\delta+o(1)),
\]
with $C$ independent of $\delta$. Near zero, both waves are $O((\delta+n^{-1})^\gamma)$. First sending $n\to\infty$, then $\delta\downarrow0$, proves
\begin{equation}\label{eq:comparison-limit}
 \sup_{m\le n/2}\left|n^{-\gamma}\psi_m(y_n)-\kappa_U\pi^\gamma
                  f(m/n,z+c_0)\right|\longrightarrow0.
\end{equation}
For the corrected variable $x$, its natural argument satisfies $2N^2(1-(x-\mathsf d_n)/\mathsf c_n)=2N^2(1-x)-c_0+o(1)$ whenever $2N^2(1-x)$ is bounded. Thus the limiting argument of $f$ is the corrected energy. Since $2N^2(1-\eta_{n,k})\to\tau_k$, symmetry and $1/w_{n,k}^0=\sum_m\psi_m(\cos x_k)^2$ give
\[
 n^{1+2\gamma}w_{n,k}^0\longrightarrow
     (\kappa_U^2\pi^{2\gamma}\fn_k)^{-1}.
\]

We also need all-node bounds. The reciprocal-derivative rule, telescoping cosine differences, gives, with $Z_n>0$ chosen so that the weights sum to one,
\[
 w_{n,k}^0=Z_n^{-1}\sin x_k\prod_{j=0}^{k-2}
 \frac{\sin((j+2\gamma_\alpha)h/2)\cos((j+2\bar\gamma)h/2)}
      {\sin((j+1)h/2)\cos((j+1+\gamma_\alpha-\gamma_\beta)h/2)}.
\]
The sine product bounds and reflection imply
\[
 c\,n^{-1}(\hat d_k/n)^{2\gamma_{(k)}}\le w_{n,k}^0
       \le Cn^{-1}(\hat d_k/n)^{2\gamma_{(k)}},\quad
 \sqrt{nw_{n,k}^0}|\psi_{n,i}^0(x)|\le C,\quad
 \sqrt{nw_{n,k}^0}|\partial_x\psi_{n,i}^0(x)|\le Cn^2/\hat d_k,
\]
where $\gamma_{(k)}=\gamma_\alpha$ for $k\le K_n$ and $\gamma_\beta$ otherwise, $i\le n/2$, and $x$ lies between $\eta_{n,k}$ and $\zeta_{n,k}$. The last two estimates follow from \eqref{eq:envelope}, its Cauchy derivative and \cref{prop:zeros}(d). These are the inputs for the target Gram matrix.

\subsection{From Gram matrices to the Green diagonal}\label{app:gram}

\begin{proposition}[Finite waves and resolvent limit]\label{prop:finite-limit}
Fix $\alpha,\beta>-1$. Choose the canonical eigenvectors positive at vertex zero and let $K_n=\lceil n/2\rceil$. With the normalized waves of \cref{prop:reconstruction},
\[
 \sup_{n\ge1}\ \max_{\substack{0\le i<n\\1\le k\le K_n}}n|\varphi_{n,k}^*(i)|^2<\infty,\qquad
 \sqrt n\,\varphi_{n,k}^*(i_n)\longrightarrow\upsilon_k(t)
 \quad\text{if }k\text{ is fixed and }i_n/n\to t\in(0,1).
\]
For each fixed $k$, $n^{1+2\gamma}w_{n,k}^*\to c_k/(\kappa_U^2\pi^{2\gamma})$.
These limits imply the diagonal, energy and trace-norm limits of \cref{thm:A}.
\end{proposition}
\begin{proof}\refstepcounter{prooflink}\label{proof:thm:A}
Write $t_i=i/n$ and $I_n=\{0,\ldots,\lfloor n/2\rfloor\}$. Both measures are canonical. Therefore their weight ratio is
\[
 \frac{w_{n,k}^\ast}{w_{n,k}^0}=\frac{R_{n,k}}{S_n},\qquad
 R_{n,k}=\prod_{\ell\ne k}\frac{|\eta_{n,k}-\eta_{n,\ell}|}{|\zeta_{n,k}-\zeta_{n,\ell}|},
 \qquad S_n=\sum_k w_{n,k}^0R_{n,k}.
\]
The logarithmic difference bound and \cref{prop:zeros} give $|\log R_{n,k}|\le C(\hat d_k^{-2}+n^{-2}\log n)$: apply the three-range product estimate \eqref{eq:app-rhohyp} using the two-edge gaps. Averaging with the weight bound proves
\begin{equation}
 |S_n-1|\le C(n^{-1-2\gamma_\alpha}+n^{-1-2\gamma_\beta}+n^{-2})\log n=o(n^{-1}).
\end{equation}
Thus $|w_{n,k}^\ast/w_{n,k}^0-1|\le C\hat d_k^{-2}+\epsilon_n$ with $n\epsilon_n\to0$. For fixed $k$ the ratio tends to $\rho_k$, and at the opposite edge to $\rho_k^\vee$. To justify the infinite product limit, cut at $L$: the remaining right-edge logarithms are $O_k(L^{-1})+O(n^{-2}\log n)$, and the opposite-edge part is $O_k(n^{-1})$.

\paragraph{The two-edge Gram limit.}
Let $I+B_n$ be the Gram matrix of $\psi_{n,i}^0$ in the target measure. Its summands are
\[
 T_k=w_{n,k}^\ast\psi_{n,i}^0(\zeta_{n,k})\psi_{n,j}^0(\zeta_{n,k})
       -w_{n,k}^0\psi_{n,i}^0(\eta_{n,k})\psi_{n,j}^0(\eta_{n,k}).
\]
Splitting the weight and node changes, the bounds just proved and the mean value theorem give $n|T_k|\le C\hat d_k^{-2}+\epsilon_n'$ with $n\epsilon_n'\to0$. For fixed right-edge $k$, \eqref{eq:comparison-limit} gives $nT_k\to F_{\alpha,k}(t_i,t_j)$, the $k$th summand of $F_\alpha$, uniformly in $i,j\in I_n$; the normalization $\kappa_U$ cancels. The left-edge summands are $(-1)^{i+j}$ times those with swapped parameters. Summable domination therefore gives
\begin{equation}\label{eq:gram-limit}
 \sup_{i,j\in I_n}|n(B_n)_{ij}-F_\alpha(t_i,t_j)-(-1)^{i+j}F_\beta(t_i,t_j)|\to0,
 \qquad |(B_n)_{ij}|\le C/n.
\end{equation}
For $i\in I_n$, let $B_{n,<i}$ be the restriction of $B_n$ to indices $0,\ldots,i-1$. The restriction of $I+B_n$ to $I_n$, and hence every principal section within $I_n$, has lower bound $\kappa_J/2$, where $\kappa_J=\min(\kappa_*(\alpha),\kappa_*(\beta))$. For a vector $v$ supported on $I_n$, pair adjacent coordinates and use their sum $P_\ell=v_{2\ell}+v_{2\ell+1}$ and difference $Q_\ell=v_{2\ell}-v_{2\ell+1}$. Up to $o(1)\nrm v^2$, their forms are Riemann discretizations of $I+\cF_{\alpha,1/2}$ and $I+\cF_{\beta,1/2}$, with the identity normalized by $\sum v_i^2=\tfrac12\sum(P_\ell^2+Q_\ell^2)$. Realize the pairs on cells of width $2/(n+4)$ inside $[0,1/2]$, padding an unpaired coordinate by zero. Give the sum and difference step functions heights $\sqrt{n+4}P_\ell/2$ and $\sqrt{n+4}Q_\ell/2$; their squared $L^2$ norms add to $\|v\|_2^2$. The cell displacement is $O(n^{-1})$. Uniform continuity makes each matrix-entry error $o(n^{-1})$, whose quadratic form is $o(1)\|v\|_2^2$ because $n^{-1}(\sum|v_i|)^2\le\|v\|_2^2$. Continuous coercivity therefore gives the claimed discrete lower bound. Extending a section vector by zero proves the same bound for every principal section.

\paragraph{Orthogonalization.}
For a positive Gram matrix with Cholesky factorization $\mathsf G=\mathsf L\mathsf L^\top$, $\mathsf L$ lower triangular with positive diagonal, orthonormalization is multiplication by $\mathsf L^{-1}$, the discrete GL reconstruction \cite[\S2]{DBK05}. For a fixed row $i$, suppress $n$ and write $B_<=B_{n,<i}$; subscripts $i,<$ and $<,i$ denote its preceding row and column. Orthogonality and normalization give
\begin{equation}\label{eq:cholesky}
 \psi_{n,i}^\ast=q_{n,i}^{-1/2}\left(\psi_{n,i}^0+\sum_{j<i}c_{n,ij}\psi_{n,j}^0\right),\quad
 c_{i,<}=-B_{i,<}(I+B_<)^{-1},\quad
 q_i=1+B_{ii}-B_{i,<}(I+B_<)^{-1}B_{<,i}.
\end{equation}
Coercivity of the principal sections and \eqref{eq:gram-limit} yield $q_i=1+O(n^{-1})$, $|c_{ij}|\le C/n$, and
\begin{equation}\label{eq:cholesky-limit}
 \sup_{j<i\le n/2}|nc_{n,ij}-K_\alpha(t_i,t_j)-(-1)^{i+j}K_\beta(t_i,t_j)|\to0.
\end{equation}
Here is the stability argument. The row equation is the discrete GL equation $nc_{ij}+nB_{ij}+n^{-1}\sum_{\ell<i}(nc_{i\ell})(nB_{\ell j})=0$. Substitute the candidate in \eqref{eq:cholesky-limit}. Same-edge products are uniform Riemann sums for GL. Mixed-edge products vanish by pairing adjacent terms of uniformly equicontinuous kernels; their error is bounded by a modulus of continuity at $1/n$ plus $O(n^{-1})$. Let $r_{ij}$ be this residual and $e_{ij}$ the difference between $nc_{ij}$ and the candidate. Put $\epsilon_n=\sup|r_{ij}|\to0$ and $\|e_i\|_{2,n}^2=n^{-1}\sum_{j<i}|e_{ij}|^2$. Subtracting the two row equations gives $e_i(I+B_<)=-r_i$, hence
\[
 \|e_i\|_{2,n}\le2\epsilon_n/\kappa_J,\qquad
 |e_{ij}|\le\epsilon_n+\frac Cn\sum_{\ell<i}|e_{i\ell}|
 \le\epsilon_n+C\|e_i\|_{2,n}\longrightarrow0.
\]
This proves the uniform limit. The same row equation first gives $\|c_i\|_2=O(n^{-1/2})$, then $\sum_j|c_{ij}|=O(1)$ and $|c_{ij}|=O(n^{-1})$; the Schur complement gives the stated estimate for $q_i$.

\paragraph{Waves, diagonal and energy.}
The all-node estimates and \eqref{eq:cholesky} now give
\[
 n|\varphi_{n,k}^\ast(i)|^2\le C\quad(k\le K_n,\ 0\le i<n),\qquad
 \sqrt n\,\varphi_{n,k}^\ast(i_n)\to\upsilon_k(t)\quad(i_n/n\to t\in(0,1)),
\]
with the sign fixed by the positive endpoint convention. For the first estimate, $\sum_{j<i}|c_{ij}|\le C$; reflect to the second half. For the limit, insert \eqref{eq:comparison-limit} and \eqref{eq:cholesky-limit} into \eqref{eq:cholesky}; the $K_\beta$ term is an alternating Riemann sum and vanishes. The remaining expression is $\kappa_U\pi^\gamma u(t,\sigma_k)$, whose weight normalization gives $\sqrt{c_k}u=\upsilon_k$. In particular
\[
 n^{1+2\gamma}w_{n,k}^\ast\to\frac{c_k}{\kappa_U^2\pi^{2\gamma}}.
\]
This records precisely where the residual dependence on $\beta$ survives before normalization.

Set $\Phi_k=\upsilon_k^2$. The exact diagonal representation is
\[
 d_i(J_n^\ast)/n=\sum_{k\le n}r_{n,k}\,n\varphi_{n,k}^\ast(i)^2.
\]
For $k>K_n$, the nodes are eventually at most $1/2$, so their entire contribution is at most $2/n$. The other summands have the uniform bound above. Fixed-mode convergence plus the $\ell^1$ tail therefore gives $d_{i_n}/n\to\sum_k\Lambda_k\Phi_k(t)=2G_\alpha(t,t)$ and a uniform diagonal bound. Dominated convergence of the step functions $(d_{\lfloor nt\rfloor}/n)^2$ proves the energy limit in \cref{thm:A}. The same argument gives $\int_0^1\Phi_k(t)\dd t=1$ and $\int_0^1\phi_\alpha(t)\dd t=M$.

The signed wave limits give $n^{-1}(I-J_n^\ast)^{-1}_{\lfloor nt\rfloor,\lfloor ns\rfloor}\to2G_\alpha(t,s)$. For the operator limit, $e_{n,k}=\mathcal U_n\varphi_{n,k}^\ast\to\upsilon_k$ in $L^2$ by bounded convergence, so their rank-one projections converge in trace norm. Cut the spectral expansion at $K$; its remaining trace norm is bounded by $\tfrac12\sum_{k>K}r_{n,k}+\sum_{k>K}\sigma_k^{-1}$. First let $n\to\infty$, then $K\to\infty$. This proves the finite-limit assertions of \cref{thm:A}.
\end{proof}

\subsection{Endpoint geometry}\label{app:geometry}

Put $u_k(t)=u(t,\sigma_k)$. The common estimates used in all three regimes are
\[
 \begin{aligned}
 |u_k(t)|&\le C\min(t^\gamma,k^{-\gamma}),\\
 |u_k(t)-f(t,\sigma_k)|&\le Ct\min(t^\gamma,k^{-\gamma}),\\
 |f(t,\sigma_k)-f(t,\tau_k)|&\le Ctk^{-3}\min(t^\gamma,k^{-\gamma}).
 \end{aligned}
\]
They follow from bounded $K$, the reference estimates and $\sigma_k-\tau_k=O(k^{-2})$. Positivity of $\phi$ and its density representation are established in \cref{prop:green-coordinate}.

When $\tau_1>0$, the reference inverse has kernel
\[
 \mathfrak G_0(t,s)=\frac{f(t\wedge s,0)f(1-(t\vee s),0)}{\mathcal W},
 \qquad \mathcal W=2f(1/2,0)f_t(1/2,0)\ne0.
\]
To check positivity, put $U_0=(\sin\pi t/\pi)^\gamma$ and $f_0=f(\cdot,0)$. The ground-state identity
\[
 [U_0^2(f_0/U_0)']'=\tau_1 U_0 f_0,\qquad
 f_0/U_0\to1,\quad U_0^2(f_0/U_0)'\to0\quad(t\downarrow0)
\]
shows that $f_0$ stays positive up to $1$ and $f_0'(1/2)>0$. Thus $\mathcal W>0$. Variation of constants shows that this kernel acts on $f(\cdot,\tau_k)$ by $1/\tau_k$: the two variation-of-constants remainders combine to $O(t^{\gamma+2})$ at zero, which selects the prescribed branch throughout $\gamma>1/4$. The positive continuous kernel therefore has the uniformly convergent Mercer expansion $\sum_k f(t,\tau_k)f(s,\tau_k)/(\tau_k\fn_k)$ \cite{RSN}. Reduction of order at zero gives
\[
 \mathfrak G_0(t,t)=f(t,0)^2\int_t^\delta f(s,0)^{-2}\dd s+Cf(t,0)^2.
\]
For $\gamma>1/2$ this is $t/(2\gamma-1)(1+o(1))$; for $\gamma=1/2$ it is $t\log(1/t)(1+o(1))$. In these two regimes $\tau_1>0$ holds automatically. The range $\gamma<1/2$ is treated below directly by its positive spectral series.

\begin{proposition}[Precise endpoint behavior]\label{thm:E}
The three asymptotics in \cref{sec:trichotomy} hold, with
\[
 c_\alpha=\sum_{k\ge1}\Lambda_k c_k
 =2\mathfrak w_\alpha\sum_{k\ge1}\frac{j_{\alpha,k}^{\alpha-1}}{|J_{\alpha+1}(j_{\alpha,k})|}\in(0,\infty)
 \quad(\alpha<-1/2).
\]
In the coordinate \eqref{eq:xpic}, $m_\alpha$ is even, has the three integrals stated there, and satisfies
\[
 \begin{cases}
 \displaystyle\lim_{|x|\to\infty}\frac{\log m_\alpha(x)}{|x|}=-2/\nu,&\nu>0,\\[2mm]
 \displaystyle\lim_{|x|\to\infty}\frac{\log\log(1/m_\alpha(x))}{|x|}=2,&\nu=0,\\[2mm]
 m_\alpha(x)=C_\alpha(x_*-|x|)^{4\gamma/(1-2\gamma)}(1+o(1)),&\nu<0,
 \end{cases}
\]
where the last support is $[-x_*,x_*]$ and
$C_\alpha=c_\alpha^2[c_\alpha(1-2\gamma)]^{4\gamma/(1-2\gamma)}$.
\end{proposition}
\begin{proof}
If $\gamma<1/2$, divide the positive spectral series by $t^{2\gamma}$. Its terms are bounded by $Ck^{2\gamma-2}$ and $u_k(t)/t^\gamma\to1$ for each fixed $k$, giving $c_\alpha$ by dominated convergence. If $\gamma=1/2$, $F=K=0$ and the reference Green asymptotic gives the answer. If $\gamma>1/2$, subtract twice the reference Green diagonal. The displayed wave bounds and $|\Lambda_k\rho_k-2/\tau_k|\le Ck^{-4}$ give
\[
 |\phi_\alpha(t)-2\mathfrak G_0(t,t)|
 \le C\sum_k k^{-4}\min((kt)^{2\gamma},1)
       +Ct\sum_k k^{-2}\min((kt)^{2\gamma},1)=o(t),
\]
by splitting at $k=1/t$. Hence $\phi\sim t/\nu$.

Changing variables in \eqref{eq:xpic} gives the three integrals. Integrating $1/\phi$ at zero gives respectively $|x|\sim\nu\log(1/t)$, $|x|\sim\frac12\log\log(1/t)$, and $x+x_*\sim t^{1-2\gamma}/[c_\alpha(1-2\gamma)]$. Substituting $m=\phi^2$ proves every claimed limit and the compact-end coefficient.
\end{proof}

The equality cases of \eqref{eq:sharp-main} add a useful comparison: for $\alpha\ge1/2$ every $L_2$ minimizer is a translate of $a_\alpha\sech^4x$, $a_\alpha=(45T_2/32)^{1/2}$. Its mass is $\sqrt{5T_2/2}$ and its tail exponent is four. The canonical exponent $2/\nu$ equals four at $\alpha=1/2$. For smaller $\alpha$ its mass is too large, and the strict inequality follows from \cref{thm:elliptic}. This separates the geometry of the canonical density from that of a relaxation optimizer.

\subsection{Shifted resolvents}\label{app:shift}

For $z\notin\spec J$ write $d_i(J;z)=[(zI-J)^{-1}]_{ii}$. Put $\cO_\alpha=\C\setminus(-\infty,-\sigma_1]$, $z_n(w)=1+w/(2n^2)$ and
\[
 \Lambda_k(w)=\frac2{\sigma_k+w},\qquad
 \phi_\alpha(t;w)=\sum_k\Lambda_k(w)\upsilon_k(t)^2
                  =2(H_\alpha+w)^{-1}(t,t).
\]
\begin{theorem}
On every compact $W\subset\cO_\alpha$ the profile series converges absolutely uniformly for $t\in[0,1]$, is holomorphic in $w$, and equals $\phi_\alpha$ at zero. Uniformly for $w\in W$, all sufficiently large $n$ satisfy
\[
 \sup_i n^{-1}|d_i(J_n^\ast;z_n(w))|\le C_W,
 \qquad n^{-1}d_{i_n}(J_n^\ast;z_n(w))\to\phi_\alpha(t;w)
 \quad(i_n/n\to t\in(0,1)),
\]
and
\[
 \begin{aligned}
 n^{-3}\sum_i d_i(J_n^\ast;z_n(w))^2&\to\int_0^1\phi_\alpha(t;w)^2\dd t,\\
 n^{-3}\sum_i |d_i(J_n^\ast;z_n(w))|^2&\to\int_0^1|\phi_\alpha(t;w)|^2\dd t.
 \end{aligned}
\]
At the left edge $d_i(J_n^\ast(\alpha,\beta);-z)=-d_i(J_n^\ast(\beta,\alpha);z)$, giving the corresponding limits $-\phi_\beta$ on $\cO_\beta$.
\end{theorem}
\begin{proof}
Write $s_{n,k}=2n^2(1-\zeta_{n,k})$ and $r_{n,k}(w)=2/(s_{n,k}+w)$. Since $s_{n,1}\to\sigma_1$ and $W$ is separated from the excluded ray, $|s+w|/(1+s)$ has a positive lower bound for $s\ge s_{n,1}$, $w\in W$. Thus $|r_{n,k}(w)|\le C_Wr_{n,k}$ and $|\Lambda_k(w)|\le C_W\Lambda_k$, and the finite shifted resolvents are defined for all large $n$. A finite spectral cutoff and \cref{lem:spec} give uniform $\ell^1$ convergence of $r_n(w)$ to $\Lambda(w)$. For real $w\in W$, $s_{n,1}+w>0$ for all large $n$, so $z_n(w)>\zeta_{n,1}$.

Now apply the proof in \cref{app:gram} with this common majorant. The contribution $k>K_n$ is at most $4/n$, and the fixed-wave limits are unchanged. It gives the diagonal bound and convergence for any moving $i_n$. Bounded convergence, using both $|a^2-b^2|$ and $||a|^2-|b|^2|\le|a-b|(|a|+|b|)$, gives the two energy formulas uniformly on $W$. Persymmetric uniqueness gives $J_n^*(\alpha,\beta)=-D J_n^*(\beta,\alpha)D$, $D=\diag((-1)^i)$, by reflection of the spectral measure. Taking resolvent diagonals proves the last assertion.
\end{proof}

At $U$, the shifted profile is
\[
 \phi_{1/2}(t;w)=\frac{2\sinh(\sqrt w\,t)\sinh(\sqrt w(1-t))}{\sqrt w\sinh\sqrt w},
\]
with value $2t(1-t)$ at $w=0$ by continuity.
For real $w>-\sigma_1$, the shifted spectrum is positive. Logarithmic differentiation of the Bessel product gives, with $I_\alpha$ the modified Bessel function of the first kind,
\[
 M(w)=\sum_k\Lambda_k(w)=\frac{I_{\alpha+1}(\sqrt w)}{\sqrt w I_\alpha(\sqrt w)},\qquad
 T_2(w)=\sum_k\Lambda_k(w)^2=-2M'(w),
\]
with analytic continuation through zero and to negative $w$ in the stated interval. Define $e(\alpha,\beta;w)$ by the same $\liminf$ as $e$, using $d_i(J;z_n(w))$.
\begin{corollary}\label{cor:realw}
For real $w>-\sigma_1$,
\[
 \ell(\Lambda(w))\le e(\alpha,\beta;w)\le\int_0^1\phi_\alpha(t;w)^2\dd t.
\]
In particular the lower side includes $\max\{\frac{2\sqrt3}{5}(5T_2(w)/2)^{3/4},\mathsf{Kel}(\Lambda_1(w),M(w))\}$.
\end{corollary}
\begin{proof}
The upper bound was just proved. For the lower one set $c_n'=(1-z_n)_+/2$, $\theta_n=z_n+c_n'\to1$. The affine fibre map $J'=(J+c_n'I)/\theta_n$ has spectrum in $(-1,1)$, because $z_n>\zeta_{n,1}$ and $2c_n'\ge1-z_n$. It satisfies $\E(J')=\theta_n^2\sum_i d_i(J;z_n)^2$ and has scaled resolvent spectrum $\theta_nr_n(w)\to\Lambda(w)$ in $\ell^1$. Apply \cref{thm:transfer-main} and let $n\to\infty$.
\end{proof}

\section{Packing, variational problems and calibrations}

We first establish spectral stability and compactness for densities, then apply the finite packing and solve the constrained problems. Throughout this appendix the variational classes consist of nonnegative $m\in L^1(\R,\dd x)$; the sharp free inequality in \cref{app:moment} also applies to the full $L^{3/2}(\R,\dd x)$ class. The functionals $\fm,\fr,E$ are those of \eqref{eq:density-functionals}, with extended values allowed.

\subsection{Variational foundation}\label{app:foundation}

\paragraph{One linear operator model.}
Identify $H^1(\R)$ functions with their continuous representatives. Equip $\mathcal H=H^1(\R)$ with
\[
 \ip{u}{w}_{\mathcal H}=\int_\R\bigl(u'(x)w'(x)+u(x)w(x)\bigr)\dd x.
\]
Set $v_x(y)=2^{-1/2}e^{-|y-x|}$ and, for $a,b\in\mathcal H$, $(a\otimes b)u=\ip{u}{b}_{\mathcal H}a$. For real $\delta\in L^1(\R,\dd x)$ define the trace-norm convergent Bochner integral
\[
 \cK_\delta=\int_\R\delta(x)\,(v_x\otimes v_x)\dd x.
\]
For bounded $B$ on $\mathcal H$, integration by parts gives
\begin{equation}\label{eq:features}
 \begin{gathered}
 \ip u{v_x}_{\mathcal H}=\sqrt2u(x),\quad \nrm{v_x}_{\mathcal H}=1,\quad
 |u(x)|\le2^{-1/2}\nrm u_{\mathcal H},\\
 \ip{v_x}{v_y}_{\mathcal H}=e^{-|x-y|},\quad \nrm{v_x-v_y}_{\mathcal H}^2\le2|x-y|,\\
 \ip{\cK_\delta u}w_{\mathcal H}=2\int_\R\delta(x)u(x)w(x)\dd x,\\
 \nrm{\cK_\delta}\le\nrm{\cK_\delta}_1\le\nrm\delta_1,\quad \Tr\cK_\delta=\int_\R\delta(x)\dd x,\\
 \Tr(B\cK_\delta)=\int_\R\delta(x)\ip{Bv_x}{v_x}_{\mathcal H}\dd x.
 \end{gathered}
\end{equation}
For $m\ge0$, $\cK_m=S^*S$ and $\sK_m=SS^*$ with $S:\mathcal H\to L^2(\R,\dd x)$, $Su=\sqrt{2m}\,u$. Thus both are positive trace class with the same nonzero spectra, including multiplicities and eigenspaces under $u\mapsto\sqrt m\,u$. Min--max and the trace-norm eigenvalue inequality \cite{Bha,Sim} imply
\[
 \sum_j|\lambda_j(\sK_m)-\lambda_j(\sK_{m'})|
 \le\nrm{\cK_m-\cK_{m'}}_1\le\nrm{m-m'}_1.
\]
The eigenvalue inequality also follows by comparing positive $X,Y$ with $X+(X-Y)_-=Y+(X-Y)_+$, where $A_\pm=(|A|\pm A)/2$ are the positive and negative spectral parts. Set $G_\infty(m)=\lambda_1(\sK_m)$. Expanding the Bochner products gives, for all integers $p\ge2$,
\begin{equation}\label{eq:Gp}
 \begin{aligned}
 G_p(m)&=\Tr\cK_m^p
 =\int_{\R^p}\prod_{i=1}^p\bigl[m(x_i)e^{-|x_i-x_{i+1}|}\bigr]\dd x_1\cdots\dd x_p,
 \qquad x_{p+1}=x_1,\\
 G_p(m)&\le G_\infty(m)^{p-1}\fm(m).
 \end{aligned}
\end{equation}
The integrable majorant $\prod_i m(x_i)$ justifies Fubini for every $L^1$ density.

The ground-level formula is
\begin{equation}\label{eq:app-BS}
 \lambda_1(\sK_m)=\sup_{0\ne u\in\mathcal H}\frac{2\int_\R m(x)u(x)^2\dd x}{\int_\R\bigl(u'(x)^2+u(x)^2\bigr)\dd x},\qquad
 -u''+u=\frac2\lambda mu
 \ \Longleftrightarrow\ \cK_mu=\lambda u
 \ \Longleftrightarrow\ u=\lambda^{-1}e^{-|\cdot|}*(mu).
\end{equation}
The quotient ranges over $H^1(\R)\setminus\{0\}$, and the equivalences hold for $\lambda>0$, $u\in H^1(\R)$, with the differential equation understood weakly. These are the Birman--Schwinger principle in our normalization. If $m\ne0$, the largest level is positive and simple: replacing a Rayleigh maximizer by its absolute value preserves the quotient, and the strictly positive convolution makes it everywhere positive. Two independent ground states would have a linear combination vanishing at a point, a contradiction. The weak solutions are $C^2$ when $m$ is continuous.

\paragraph{First variations with quantitative remainders.}
Fix a nonzero density $m$, let $u>0$ be its ground state in \eqref{eq:app-BS}, and abbreviate $\lambda_j=\lambda_j(\sK_m)$. Let $\delta\in L^1(\R,\dd x)$ with $m+\delta\ge0$, and write $W_2=e^{-2|\cdot|}*m$ and $G_2^{\rm bil}(h,k)=\iint_{\R^2}h(x)e^{-2|x-y|}k(y)\dd x\dd y$. Bilinearity and perturbation of the isolated simple top eigenvalue give
\begin{equation}\label{eq:variations}
 \begin{aligned}
 G_2(m+\delta)&=G_2(m)+2\int_\R W_2(x)\delta(x)\dd x+G_2^{\rm bil}(\delta,\delta),
 &|G_2^{\rm bil}(\delta,\delta)|&\le\nrm\delta_1^2,\\
 \lambda_1(\sK_{m+\delta})&=\lambda_1(\sK_m)+\int_\R g_\infty(x)\delta(x)\dd x+R,
 &|R|&\le2\nrm\delta_1^2/\epsilon_*,
 \end{aligned}
\end{equation}
where $g_\infty(x)=\lambda_1u(x)^2/\bigl(\int_\R m(y)u(y)^2\dd y\bigr)$, $\epsilon_*=\lambda_1-\lambda_2>0$ and $\nrm\delta_1\le\epsilon_*/4$. For the second line use $\nrm{\cK_\delta}\le\nrm\delta_1$; separating the top eigendirection in the Rayleigh quotient gives the stated quadratic remainder (equivalently, simple-eigenvalue perturbation \cite{Kat}).

For Fredholm let $f_s(x)=-\log(1-sx)-sx$ for $x<1/s$ and $+\infty$ for $x\ge1/s$. If $s\lambda_1\le1-\epsilon$, the coercive form $-\partial^2+1-2sm$ has continuous positive Green kernel
\[
 G(x,y)=\tfrac12\ip{(I-s\cK_m)^{-1}v_x}{v_y}_{\mathcal H},\qquad
 \tfrac12\le g(x):=G(x,x)\le[2(1-s\lambda_1)]^{-1}.
\]
Moreover $2G(x,y)\ge e^{-|x-y|}$ by the positive Neumann expansion. The gradient of $\hD_s$ is $g_s=s(2g-1)$, and for $\nrm\delta_1\le\epsilon/(2s)$,
\begin{equation}\label{eq:fred-variation}
 |\hD_s(m+\delta)-\hD_s(m)-\int_\R g_s(x)\delta(x)\dd x|
 \le2s^2\nrm\delta_1^2/\epsilon^2.
\end{equation}
On two densities below this same pole margin, $|\hD_s(m)-\hD_s(m')|\le(s/\epsilon)\nrm{m-m'}_1$. The trace derivative follows from $f_s'(X)=s(I-sX)^{-1}-sI$ and the resolvent identity. The trace remainder is between zero and $\Tr[(f_s'(X+Y)-f_s'(X))Y]$ by Klein's inequality. For compact self-adjoint $X,Z$ on a separable space, with spectra in a compact interval $I$, convex $f\in C^1(I)$ and trace-class $f(X),f(Z),Z-X$, its infinite-dimensional version follows in orthonormal eigenbases $Xe_j=x_je_j$, $Zh_k=z_kh_k$ by summing the nonnegative quantities
\[
 \sum_k|\ip{e_j}{h_k}|^2[f(z_k)-f(x_j)-f'(x_j)(z_k-x_j)].
\]
Bounded $f'$ and the trace-class assumptions justify the trace. This proves the inequality and \eqref{eq:fred-variation}; the spatial determinant representation will be proved in \cref{app:fredholm}.

\paragraph{Rearrangement and compactness.}
Symmetric decreasing rearrangement preserves $\fm$, $\fr$ and $E$ from \eqref{eq:density-functionals}, and increases every $G_p$, $G_\infty$ and $\hD_s$. Apply Brascamp--Lieb--Luttinger \cite{BLL} to the cyclic integral \eqref{eq:Gp}; the top level follows from its moment limit or the Rayleigh rearrangement, and the determinant from its positive series, with the pole case treated separately. When $m$ is even and nonincreasing, $\cK_m u=e^{-|\cdot|}*(mu)$ preserves that cone. The normalized powers $\lambda_1^{-j}\cK_m^jv_0$ converge in $H^1(\R)$ to a positive multiple of the simple ground state; the coefficient is positive since $\langle v_0,u\rangle_{\mathcal H}=\sqrt2u(0)>0$. The ground state is therefore even and nonincreasing.

For $a>0$ and integrable $h\ge0$, $k_a*h$, $k_a=e^{-a|\cdot|}$, is $C^1$, tends to zero together with its first derivative, and satisfies $-(k_a*h)''+a^2(k_a*h)=2ah$ distributionally. For nonzero even decreasing $h$, its derivative is strictly negative on $x>0$: the two half-convolutions differ by $\int_0^\infty e^{-av}[h(x-v)-h(x+v)]\dd v>0$. Equality would make $h$ periodic with period $2x$, impossible for a nonzero integrable function. Applying this also to $h=mu$ in \eqref{eq:app-BS} proves $u'<0$ and $W_2'<0$ on $x>0$. Thus both moment gradients in \eqref{eq:variations} are strictly decreasing and continuous.

For an even decreasing density with mass at most $M'$, energy at most $E_1$ and, when assumed, root mass at most $r_1$,
\[
 m(x)\le\min\left(\frac{M'}{2|x|},\left(\frac{E_1}{2|x|}\right)^{2/3},
                         \frac{r_1^2}{4x^2}\right).
\]
The third bound applies when the root cap is imposed. For $m_R=m\ind_{[-R,R]}$ the spectral tails satisfy
\begin{equation}\label{eq:tail-foundation}
 |G_2(m)-G_2(m_R)|\le3M'^2/R,\quad
 0\le G_\infty(m)-G_\infty(m_R)\le M'/R,
\end{equation}
and $0\le\hD_s(m)-\hD_s(m_R)\le s^2M'^2/(\epsilon R)$ below pole margin $\epsilon$. Indeed the tail feature operator has norm at most $2\nrm{m-m_R}_\infty\le M'/R$ and trace at most $M'$; expand $\Tr(X+Y)^2$, apply min--max, and use Klein with $f_s'(X+Y)\le(s^2/\epsilon)(X+Y)$.

\begin{lemma}[Spectral compactness]\label{lem:spectral-compactness}
Let $m_j$ be even nonincreasing densities with $\int_\R m_j(x)\dd x\le M'$ and $E(m_j)\le E_1$. A subsequence converges locally in $L^1$ to an even nonincreasing density $m$, and
\[
 \begin{gathered}
 E(m)\le\liminf_jE(m_j),\qquad
 \fr(m)\le\liminf_j\fr(m_j),\\
 \fm(m)\le M',\qquad G_i(m_j)\to G_i(m)\quad(i=2,\infty).
 \end{gathered}
\]
For each fixed $s>0$, $\limsup_j\hD_s(m_j)\le\hD_s(m)$. If the root masses are uniformly bounded, convergence is global in $L^1$ and $\int_\R m_j(x)\dd x\to\int_\R m(x)\dd x$.
\end{lemma}
\begin{proof}
Helly's selection theorem on intervals away from zero and a diagonal subsequence give pointwise convergence at almost every $x\ne0$. Fatou gives the displayed lower semicontinuity and mass cap. The energy bound gives the locally integrable majorant $C|x|^{-2/3}$, hence local $L^1$ convergence. Apply \eqref{eq:tail-foundation} before letting $R\to\infty$ to obtain convergence of $G_2,G_\infty$. If $sG_\infty(m)\ge1$, the determinant statement follows from its extended value. Otherwise top-level convergence gives a common pole margin; local $L^1$ continuity of the determinant and its tail bound then give convergence. A root-mass bound adds the integrable $Cx^{-2}$ tail, giving global $L^1$ convergence and convergence of the masses.
\end{proof}

To state the full variants, for $S\subseteq\{2,\infty\}$ and thresholds $T'=(T_i')_{i\in S}$ put
\[
 \begin{aligned}
 \cL_S^{\rm r}(M',T')&=\inf\{E(m):\fm=M',\ \fr\le1,\ G_i\ge T_i'\ (i\in S)\},\\
 \cL_S^{\rm rf}(M',T')&=\inf\{E(m):\fm\le M',\ G_i\ge T_i'\ (i\in S)\}.
 \end{aligned}
\]
Define $\cB_s^{\rm r}(M',d')$ by the first scalar constraints and $\hD_s\ge d'$, and $\cB_s(M',d')$ with $\fm\le M'$ as its scalar constraint. An argument $\Lambda$ means $M'=\sum_k\Lambda_k$, $T_2'=\sum_k\Lambda_k^2$ and $T_\infty'=\Lambda_1$. At these data the short names are $L_2=\cL_{\{2\}}^{\rm rf}$, $L_\infty=\cL_{\{\infty\}}^{\rm rf}$ and $B_s=\cB_s$. \Cref{lem:spectral-compactness} proves lower semicontinuity of $\cL_S^{\rm r}$ and $\cB_s^{\rm r}$ under simultaneous convergence of the mass and thresholds. Rearrange approximate minimizers, apply global $L^1$ compactness, and pass each constraint in its established direction.

\begin{proposition}[Existence]
The root-free moment problem is finite and attained at an even decreasing nonzero density whenever $S\ne\emptyset$, $0<T'_2<M'^2$ and $0<T'_\infty<M'$ for the constraints included. The same holds for $\cB_s(M',d')$ if $0<d'<f_s(M')$, with $f_s(M')=+\infty$ for $sM'\ge1$. These conditions hold for every edge spectrum used here.
\end{proposition}
\begin{proof}
The narrow box $(M'/\varepsilon)\ind_{[-\varepsilon/2,\varepsilon/2]}$ has feature operator tending in trace norm to $M'v_0\otimes v_0$. It satisfies the strict thresholds for small $\varepsilon$ and has finite energy. Apply the compactness just proved to a rearranged minimizing sequence. For edge data, $T_2<M^2$, $\Lambda_1<M$ and $\hd_s<f_s(M)$; the last inequality follows from strict increase of $f_s(x)/x$, or is immediate if $sM\ge1$.
\end{proof}

\paragraph{Normal multipliers and the support.}
For $m_t(x)=tm(tx)$, $t\ge1$, mass is preserved, $\fr(m_t)=t^{-1/2}\fr(m)$ and $E(m_t)=t^{1/2}E(m)$. Also $\nrm{m_t-m}_1\le2M'(t-1)$ for even decreasing $m$. Each of $G_2,G_\infty$ and the subcritical determinant has a strictly positive linear gain under this concentration. For moments, rescaling \eqref{eq:Gp} differentiates $e^{-|x-y|/t}$; its positive distance factor has positive integral. For the top level, unitary rescaling gives the kernel $\sqrt{m(x)}e^{-|x-y|/t}\sqrt{m(y)}$. At $t=1$ its derivative is $\sqrt{m(x)}|x-y|e^{-|x-y|}\sqrt{m(y)}$. Its quadratic form on the positive normalized top eigenfunction is strictly positive; testing that function in Rayleigh's principle gives the linear gain. For Fredholm keep the second-moment term of its positive series; the top-level Lipschitz bound retains the pole margin.

\begin{lemma}[Common multiplier rule]\label{lem:app-kkt}
Let $m_*\ne0$ minimize $E$ under $\fm\le M'$ and $\Phi\ge T$, with $m_*$ even decreasing and of finite energy. Suppose $\Phi(m_*)<\infty$, has bounded gradient $g$ and an $O(\nrm\delta_1^2)$ local remainder, and has the positive concentration gain just described. Then $\Phi(m_*)=T$ and
\[
 \sqrt{m_*}=\tfrac23(\mu g-\lambda_{\fm})_+,
 \quad\mu>0,\quad\lambda_{\fm}\ge0,\quad
 \lambda_{\fm}(M'-\fm(m_*))=0.
\]
\end{lemma}
\begin{proof}
If the constraint were slack, decreasing the amplitude slightly would lower $E$ and remain admissible. On the cone of directions $d=\theta(\varphi-m_*)$, $\theta\ge0$, $\varphi\ge0$ of finite mass and energy, the derivatives are $E'd=(3/2)\int_\R\sqrt{m_*(x)}d(x)\dd x$, $\fm'd=\int_\R d(x)\dd x$ and $\Phi'd=\int_\R g(x)d(x)\dd x$. Minimality excludes $E'd<0$, $\Phi'd>0$ with nonpositive mass derivative at an active cap. Separation of these finite-dimensional derivative images gives nonnegative $a,\lambda_{\fm},\mu$, with $(a,\lambda_{\fm},\mu)\ne0$ and
\[
 aE'd+\lambda_{\fm}\fm'd-\mu\Phi'd\ge0
\]
on the cone; omit the mass coordinate when slack. For $d=m_t-m_*$, the concentration gain and quadratic remainder make $\Phi'd>0$ and $\fm'd=0$. If $a=0$, this forces $\mu=0$, and $d=\pm m_*$ forces $\lambda_{\fm}=0$, a contradiction. Normalize $a=1$. Adding a bounded indicator and removing the density on a bounded measurable set show that $\Psi=(3/2)\sqrt{m_*}-\mu g+\lambda_{\fm}$ is nonnegative and $\Psi m_*=0$. This is the clipping formula. If $\mu=0$, it would give $m_*=0$, again impossible.
\end{proof}

For $\Phi=G_i$, $i=2,\infty$, all hypotheses are \eqref{eq:variations} and concentration. Thus the constraint is active, the clipped coefficient is continuous, even and strictly decreasing, and the positive set is one centered interval. For Fredholm the constraint is active and the clipping formula holds with $g_s=s(2g-1)$ at a subcritical minimizer. The remainder \eqref{eq:fred-variation} verifies its hypothesis; finite-energy admissible approximations are obtained by concentrating slightly, then truncating in height and space. Exclusion of a pole minimizer for our calibration is proved below.

\subsection{Packing and transfer}\label{app:packing}

Part~II, Theorem~2.3, constructs the Green coordinates and disjoint intervals, including reducible blocks. We apply that construction with the present normalization and then transfer general limiting spectra.

\begin{proof}[Proof of \Cref{prop:packing}]\refstepcounter{prooflink}\label{proof:prop:packing}
For irreducible $J$, the centering lift of Part~II gives a zero-diagonal chain with inverse $Q$, points $x_k$, and $Q_{2i-1,2j-1}=2R_{ij}$, $R=(I-J)^{-1}$. Its exponential Green coordinates give $Q_{k\ell}=\sqrt{Q_{kk}Q_{\ell\ell}}e^{-|x_k-x_\ell|}$. If $\varrho_k=\coth(x_{k+1}-x_k)$, with endpoint value $\varrho=1$, then $d_i=(\varrho_{2i-2}+\varrho_{2i-1})/4$. The intervals extending to the left and right of $y_i=x_{2i-1}$ by $4\varrho_{2i-1}/(\varrho_{2i-2}+\varrho_{2i-1})^2$ and $4\varrho_{2i-2}/(\varrho_{2i-2}+\varrho_{2i-1})^2$ have length $1/d_i$ and disjoint interiors. For $a,b>0$, $4a/(a+b)^2\le1/b$. Each extension is therefore at most the corresponding $\tanh l\le l$, so neighbouring intervals have disjoint interiors.

The spectral meaning is especially simple. For $X=\sum_i(d_i/n^2)v_{y_i}\otimes v_{y_i}$, its finite Gram matrix is
\begin{equation}
 \Gamma_{ij}=\frac{\sqrt{d_id_j}}{n^2}e^{-|y_i-y_j|}=R_{ij}/n^2.
\end{equation}
It therefore has the same nonzero spectrum. Replacing each point mass by its interval of constant height $n^{-2}d_i^2$ preserves the three integrals in \eqref{eq:packing-main}. Moreover $\nrm{v_x\otimes v_x-v_y\otimes v_y}_1\le2\sqrt{2|x-y|}$ and $\int_{I_i}|x-y_i|\dd x\le1/(2d_i^2)$. Cauchy--Schwarz against the mass $M_J$ gives
\begin{equation}
 \nrm{\cK_{m_J}-X}_1\le2\sqrt{M_J/n}.
\end{equation}
For reducible $J$, perform the same construction in each irreducible block with the global scale $n^{-2}$. Separate the translated blocks by distance $\tau$. The cross-block Gram error is at most $e^{-\tau}nM_J$ in trace norm. Choose $\tau$ to make it at most $2\sqrt{M_J/n}$. The eigenvalue inequality in \eqref{eq:features} gives the uniform constant four in \eqref{eq:packing-main}. \end{proof}

\begin{theorem}[Full transfer]
Under the hypotheses of \cref{thm:transfer-main}, for every $S\subseteq\{2,\infty\}$ and $0<s<1/\Lambda_1$,
\[
 \liminf_n\inf_{\cJ(p_n)}n^{-3}\E
 \ge\cL_S^{\rm r}(\Lambda)\ge\cL_S^{\rm rf}(\Lambda),\qquad
 \liminf_n\inf_{\cJ(p_n)}n^{-3}\E
 \ge\cB_s^{\rm r}(M,\hd_s)\ge B_s(\Lambda).
\]
\end{theorem}
\begin{proof}
Take a subsequence and realizations whose normalized energies converge to the finite lower limit; the case of an infinite lower limit satisfies the claim immediately. Put $M_n=\sum_k r_{n,k}$ and $\delta_n=4\sqrt{M_n/n}+\nrm{r_n-\Lambda}_{\ell^1}\to0$. Every packed competitor satisfies
\begin{equation}
 \nrm{\lambda(\sK_{m_J})-\Lambda}_{\ell^1}\le\delta_n.
\end{equation}
Thus the top-level error is at most $\delta_n$, and the second-moment error at most $(M_n+M)\delta_n$. For fixed $s$, both spectra eventually stay below one common pole margin $\epsilon$, so the determinant error is at most $s\delta_n/\epsilon$. The packed density is admissible for the root-constrained problems at mass $M_n$, thresholds reduced by these errors, and has exactly the finite normalized energy. Lower semicontinuity in all data, proved in \cref{app:foundation}, gives the two first inequalities. The last inequalities are inclusions of admissible classes.
\end{proof}

\subsection{Moment and ground-state calibrations}\label{app:moment}

\paragraph{The classical free constant.}
The sharp $H^1(\R)\to L^3(\R,\dd x)$ inequality \cite{DELL14}, after scaling the derivative term, reads
\[
 \nrm v_3^2\le B\int_\R\bigl(v'(x)^2+4v(x)^2\bigr)\dd x,\qquad
 B=\frac5{32}(16/15)^{2/3}.
\]
Its nonnegative nonzero extremizers are $a\sech^2(x-x_0)$. This fixes the normalization by the integrals of $\sech^4$ and $\sech^6$. For $W=e^{-2|\cdot|}*m$, $\int_\R\bigl(W'(x)^2+4W(x)^2\bigr)\dd x=4\int_\R m(x)W(x)\dd x$. H\"older and this inequality give $G_2\le4B\nrm m_{3/2}^2=C_{\rm GN}E^{4/3}$. Equality in H\"older requires $m=cW^2$ for some $c>0$, yielding precisely the equality cases in \eqref{eq:sharp-main}; conversely $e^{-2|\cdot|}*\sech^4=\tfrac23\sech^2$. Approximation of nonnegative $L^{3/2}$ data in the dual Sobolev space, or monotone truncation followed by this bound, proves the statement even with infinite mass. In particular $G_2(\sech^4)=32/45$, $E(\sech^4)=16/15$ and $\int_\R\sech^4x\dd x=4/3$. At the $U$ data, the full-spectrum calculation in \cref{sec:second} and the sharp bound imply $\cL_S^{\rm r}=\cL_S^{\rm rf}=2/15$ for every $S\ni2$, as well as the value $2/15$ under the full spectrum and both scalar normalizations.

\paragraph{A common quadrature estimate.}
The following estimate also supplies the Fredholm monotonicity proof. For $0<\tau\le1$, $0\le\chi<1$ put
\begin{equation}\label{eq:common-weight}
 \begin{gathered}
 L=1-\tau(1-v),\quad N_\chi=1-\chi\tau(1-v),\quad
 P=L^2-v^3N_\chi=(1-v)r,\\
 q=v^2+(1-\tau^2)v+(1-\tau)^2,\quad r=q+\chi\tau v^3,
 \quad I[g]=\int_0^1g(v)P(v)^{-1/2}\dd v.
 \end{gathered}
\end{equation}
For $k\ge2$ write $\langle h\rangle_k=I[v^kh]/I[v^k]$. Then
\[
 v^2\le L^2\le q\le r\le1+v+v^2\le3,\qquad
 \langle v\rangle_2:=I[v^3]/I[v^2]\ge2/3.
\]
The first chain is immediate by expansion and $\partial_\tau r\le0$. For the mean, split $v^2(v-2/3)$ at $2/3$, bounding the negative part with $P^{-1/2}\le v^{-1}(1-v)^{-1/2}$ and the positive part with $P^{-1/2}\ge[3(1-v)]^{-1/2}$. The resulting elementary integrals satisfy
\[
 \int_0^{2/3}\frac{v(2/3-v)}{\sqrt{1-v}}\dd v<8/81
 <368/2835=\frac1{\sqrt3}\int_{2/3}^1\frac{v^2(v-2/3)}{\sqrt{1-v}}\dd v.
\]
For $0<v<1$, also $0<\xi=(2L-\chi v^3)/(2r)\le1/v$, $0\le v^3/r\le v$ and $\xi'<0$, with
\begin{equation}\label{eq:weight-derivatives}
 \partial_\tau P^{-1/2}=\xi P^{-1/2},\qquad
 \partial_\chi P^{-1/2}=-(\tau/2)(v^3/r)P^{-1/2}.
\end{equation}
For the strict sign, $-2r^2\xi'$ is the positive polynomial
\[
 \chi v^2[v^2+2(1+\tau^2)v+3(1-\tau^2)]
 +2(1-\tau)^2+4(1-\tau)v+2\tau v^2.
\]
These formulas supply integrable majorants for all subsequent differentiations, including the physical one-sided limits.

\begin{theorem}[Mass-constrained second moment]\label{thm:elliptic}
Let
\[
 Q_\tau=(1-v)[v^2+(1-\tau^2)v+(1-\tau)^2],\quad
 i_k(\tau)=\int_0^1\frac{v^k\dd v}{\sqrt{Q_\tau(v)}},\quad
 \kappa(\tau)=\int_0^1\frac{v^2(1-\tau+\tau v)\dd v}{\sqrt{Q_\tau(v)}}.
\]
The function $\mathsf r=2\kappa/(3i_2^2)$ is continuous and strictly decreasing on $(0,1]$, from $1$ to $2/5$. If $5T_2>2M^2$, its unique root $\mathsf r(\tau_\Lambda)=T_2/M^2$ lies in $(0,1)$ and
\[
 L_2=M^{3/2}\Xi(\tau_\Lambda)>\ell_Y,\qquad
 \Xi(\tau)=\frac{i_3(\tau)}{\sqrt\tau\,i_2(\tau)^{3/2}}.
\]
For Jacobi data this is $-1<\alpha<1/2$ and $\mathsf r(\tau_\alpha)=1/(\alpha+2)$. In the remaining regime $L_2=\ell_Y$, attained by the free extremizer described in \cref{sec:second}.
\end{theorem}
\begin{proof}
Start from the optimality equation \eqref{eq:moment-mechanism}. With the multiplier $\mu>0$ from \cref{lem:app-kkt}, set $A=27/(16\mu)$ and $Ay=2\mu W_2-\lambda_{\fm}$, $\sigma=-\lambda_{\fm}/A\le0$. Then
\begin{equation}
 -y''+4y=4\sigma+6y_+^2,\qquad
 y'^2=4[(y-\sigma)^2-y_+^3],\quad
 m=\tfrac49A^2y_+^2,\quad W_2=\tfrac{8A^2}{27}(y-\sigma).
\end{equation}
The first integral is \eqref{eq:moment-mechanism} in these coordinates. The even solution decreases strictly to $\sigma$. Its positive maximum $y_0$ satisfies $y_0^{3/2}-y_0=|\sigma|$, giving a unique $y_0\ge1$. The polynomial $P_\sigma(y)=(y-\sigma)^2-y_+^3$ is positive on $(\sigma,y_0)$: for $y>0$ compare the strictly decreasing $y^{-1/2}+|\sigma|y^{-3/2}$ with one. Hence one change of variable gives, for all nonnegative $F$ vanishing on $y\le0$,
\[
 \int_\R F(y(x))\dd x=\int_0^{y_0}\frac{F(y)}{\sqrt{P_\sigma(y)}}\dd y.
\]
Writing $I_k=\int_0^{y_0}y^k/\sqrt{P_\sigma(y)}\dd y$ and $J=\int_0^{y_0}y^2(y-\sigma)/\sqrt{P_\sigma(y)}\dd y$ yields
\[
 \fm=\tfrac49A^2I_2,\quad E=\tfrac8{27}A^3I_3,\quad
 G_2=\tfrac{32}{243}A^4J,\quad\fr=\tfrac23AI_1.
\]
If $\sigma=0$, the normalized Cauchy problem gives $y=\sech^2x$, whose mass exceeds $M$ in the present regime. Thus $\lambda_{\fm}>0$ and $\fm=M$. Set $y_0=\tau^{-2}$, $\sigma=\tau^{-2}-\tau^{-3}$ and $y=y_0v$. Then $P_\sigma=y_0^3Q_\tau$, so $T_2/M^2=\mathsf r(\tau)$ and $E=M^{3/2}\Xi(\tau)$.

The common weight bounds give continuity up to both ends. At $\tau=1$, $i_2=4/3$, $i_3=\kappa=16/15$; at $\tau=0$, $i_2=\kappa=2/3$. For strict monotonicity put $\phi=L/q$ and take means with weight $v^2Q_\tau^{-1/2}$. Differentiation justified by \eqref{eq:weight-derivatives} gives
\[
 \frac{\kappa'i_2-2\kappa i_2'}{i_2^2}
 =-\langle1-v\rangle_2+\langle L\phi\rangle_2
                   -2\langle L\rangle_2\langle\phi\rangle_2<0.
\]
Indeed $L\le1$ and $\langle L\rangle_2\ge\langle v\rangle_2\ge2/3$ bound the right side by $-\langle1-v\rangle_2-\tfrac13\langle\phi\rangle_2<0$. This proves the unique root and exact value, since a minimizer exists. The free extremizer has mass exceeding $M$, so the equality classification in \eqref{eq:sharp-main} gives strictness.
\end{proof}

\begin{theorem}[Keller with a mass cap]\label{thm:keller}
For $c\ge1$ define
\[
 j_k(c)=\int_0^1\frac{(1-\sigma^2)^k\dd\sigma}
 {\sqrt{(1-\sigma^2/c)[c(3-3\sigma^2+\sigma^4)-1]}},\quad k=1,2,3.
\]
For $0<\lambda<M$, $\mathsf{Kel}(\lambda,M)=\inf\{E(m):m\ge0,\ \fm(m)\le M,\ \lambda_1(\sK_m)\ge\lambda\}$ is
\[
 \begin{cases}
 \displaystyle\frac\pi4(3\lambda/2)^{3/2},&M\ge3\lambda/2,\\[1mm]
 \displaystyle2(3\lambda/2)^{3/2}c\,j_3(c),&\lambda<M<3\lambda/2,
 \end{cases}
\]
where the latter $c>1$ is the unique root $3\lambda\sqrt c\,j_2(c)=M$. Indeed $\sqrt c\,j_2(c)$ decreases continuously from $1/2$ to $1/3$. Thus $L_\infty=\mathsf{Kel}(\Lambda_1,M)$; the constrained branch for Jacobi data occurs exactly when $j_{\alpha,1}^2<6(\alpha+1)$.
\end{theorem}
\begin{proof}
The free value is Keller's inequality \cite{Kel,LT}, with density $(3\lambda/2)\sech^2(2x)$, mass $3\lambda/2$. To classify the cap, normalize the minimizer's positive ground state by $u(0)=1$. The common rule gives $m=a^2(u^2-v_e)_+^2$, $0\le v_e<1$. Its ODE and decay give
\[
 u'^2=u^2-\frac{2a^2}{3\lambda}(u^2-v_e)_+^3,
 \quad c=(1-v_e)^{-1},\quad a^2=(3\lambda/2)c^3.
\]
On the positive half-support substitute $u^2=1-\sigma^2/c$. The first integral gives
\[
 \dd x=\frac{\dd\sigma}{\sqrt c\sqrt{(1-\sigma^2/c)[c(3-3\sigma^2+\sigma^4)-1]}},
\]
so $\fm=3\lambda\sqrt c\,j_2$, $E=2(3\lambda/2)^{3/2}c\,j_3$ and $\fr=2(3\lambda/2)^{1/2}j_1$. For $\varepsilon=1/c$, the integrand defining $\sqrt c\,j_k$ is strictly increasing in $\varepsilon$ and dominated by its bounded value at $\varepsilon=1$. At $c=1$, $j_1=\pi/4$, $j_2=1/2$, $j_3=\pi/8$; at infinity the $j_2$ limit is $1/3$, using the derivative of $\sigma\sqrt{3-3\sigma^2+\sigma^4}$. Thus the mass curve is strictly decreasing. Complementarity gives $c=1$ in the free regime and the unique mass-active $c>1$ otherwise. All values are attained by existence; the algebraic substitution of \eqref{eq:edgedata} proves the final criterion.
\end{proof}

\subsection{Fredholm calibration}\label{app:fredholm}

We derive the spatial determinant identity, classify the minimizing profiles, and construct their positive Green kernels. The parameter monotonicity then selects the admissible branches.

\paragraph{The Green diagonal and determinant.}
For continuous bounded $m$ below the pole, the coercive Green kernel from \cref{app:foundation} factors as
\[
 G(x,y)=\varphi_-(\min(x,y))\varphi_+(\max(x,y))/W,
 \quad W=\varphi_+\varphi_-'-\varphi_+'\varphi_->0,
\]
where $\varphi_\pm>0$ span the solutions square integrable at the respective infinities for $-\varphi''+(1-2sm)\varphi=0$. These spaces are one-dimensional: an $L^2$ solution has $L^2$ derivative by a cutoff energy estimate, so a nonzero constant Wronskian between two such solutions would be integrable. The positive Green columns construct the two solutions; coercivity excludes a common whole-line $H^1$ solution and fixes their independence. The derivative jump of $G$ fixes the displayed normalization. Direct differentiation of $g=G(x,x)$ gives
\begin{equation}\label{eq:GD}
 2gg''-g'^2-4(1-2sm)g^2+1=0.
\end{equation}
Thus $g\in C^2$. The variation formula for $L^1$ densities is \eqref{eq:fred-variation}; the continuous profiles constructed below satisfy this classical ODE representation.

For coupling $0\le\sigma\le s$, resolvent differentiation gives $\partial_\sigma g_\sigma(x)=2\int_\R m(z)G_\sigma(x,z)^2\dd z$. The ratio $\varphi_+/\varphi_-$ decreases to zero at $+\infty$ and its derivative is $-W/\varphi_-^2$; otherwise both solutions would be square integrable there. For $x\ge z$, $G_\sigma(x,z)/g_\sigma(x)=\varphi_-(z)/\varphi_-(x)$, so integrating this ratio derivative gives
$\int_z^\infty G_\sigma(x,z)^2/g_\sigma(x)^2\dd x=g_\sigma(z)$.
The reflected argument gives the same value on $(-\infty,z)$, hence
\[
 \int_\R G_\sigma(x,z)^2/g_\sigma(x)^2\dd x=2g_\sigma(z).
\]
Integrate $\partial_\sigma(1-1/(2g_\sigma))$ and use Tonelli, all terms being nonnegative. The result is
\[
 \begin{aligned}
 \int_0^s\int_\R 2m(x)g_\sigma(x)\dd x\dd\sigma
 &=\int_0^s\Tr[(I-\sigma\cK_m)^{-1}\cK_m]\dd\sigma\\
 &=-\log\det(I-s\sK_m).
 \end{aligned}
\]
Subtracting $s\fm(m)$ gives
\begin{equation}\label{eq:Dg}
 \hD_s(m)=\int_\R\left(1-\frac1{2g(x)}-sm(x)\right)\dd x.
\end{equation}
In particular both integral representations are convergent. For Jacobi data the target follows from the Bessel product:
\begin{equation}
 \prod_{k\ge1}(1-z^2/j_{\alpha,k}^2)=\Gamma(\alpha+1)(2/z)^\alpha J_\alpha(z),\qquad
 \hd_s(\alpha)=-\log[\Gamma(\alpha+1)(2/z)^\alpha J_\alpha(z)]-sM,
 \quad z=\sqrt{2s}.
\end{equation}

\begin{theorem}[Exact Fredholm value]\label{thm:fred}
Use $L,N_\chi,P,I$ from \eqref{eq:common-weight} and define the dimensionless state functions
\[
 \fm(\tau,\chi)=\tfrac32\chi I[v^2],\quad
 \fD(\tau,\chi)=\log\frac{1-\chi\tau}{1-\chi}
          +\tau\chi I[L/N_\chi]-\fm(\tau,\chi),\quad
 \fE(\tau,\chi)=(3\chi)^{3/2}\tau^{-1/2}I[v^3].
\]
For any edge spectrum and $0<s<1/\Lambda_1$, the mass-active system
\[
 \fm(\tau,\chi)=sM,\qquad \fD(\tau,\chi)=\hd_s(\Lambda),\qquad 0<\tau,\chi<1
\]
has at most one solution. The free equation $\fD(1,\chi)=\hd_s(\Lambda)$ has at most one solution $0<\chi<1$. At least one of the following candidates is admissible, and $B_s$ is their minimum:
\[
 (2s)^{-3/2}\fE(\tau,\chi)\quad\text{for the mass-active solution};\qquad
 (2s)^{-3/2}\fE(1,\chi)\quad\text{if }\fm(1,\chi)\le sM.
\]
\end{theorem}

\paragraph{Reduction and shape.}
An even decreasing minimizer exists. To prove $s\lambda_1<1$, suppose $s\lambda_1\ge1$. Then amplitudes $a\uparrow[s\lambda_1(\sK_m)]^{-1}\le1$ have determinant tending to infinity and, for $a<1$ close enough, are admissible with lower energy. The common multiplier rule therefore applies. Put
\begin{equation}\label{eq:app-mstar}
 g_e=\tfrac12+\lambda_{\fm}/(2\mu s),\quad
 \widetilde a=\tfrac43\mu s\sqrt{2s},\qquad
 m=\frac{\widetilde a^2}{2s}(g-g_e)_+^2.
\end{equation}
This makes $m$ continuous and bounded. Moreover $g\to1/2$ by the bound $2g-1\le s(1-s\lambda_1)^{-1}(e^{-2|\cdot|}*m)$. Integrating \eqref{eq:GD} after dividing by $g^2$ gives
\begin{equation}\label{eq:Fg}
 g'^2=F(g),\qquad F(\gamma)=(2\gamma-1)^2
                  -\tfrac43\widetilde a^2\gamma(\gamma-g_e)_+^3.
\end{equation}
The integration constant is fixed by $g\to1/2$: the identity first makes $g'^2$ converge, and a nonzero limit would contradict convergence of $g$, so $g'\to0$.

There is exactly one zero $g_0>g_e$ of $F$, with $F'(g_0)<0$, $F>0$ on $(1/2,g_0)$ and $F<0$ on $(g_0,\infty)$. Indeed $R(\gamma)=\gamma(\gamma-g_e)^3/(2\gamma-1)^2$ increases from zero to infinity, since
$(\log R)'=1/\gamma+3/(\gamma-g_e)-2/(\gamma-1/2)>0$.
The even solution has $g(0)=g_0$. The other possibility $g(0)=1/2$ would give the zero density by ODE uniqueness. Any later zero of $g'$ would similarly force a constant or a periodic solution, so $g'<0$ on $(0,\infty)$. Writing $a_e$ for the support half-length, it is finite exactly when $g_e>1/2$; outside the support, $g-1/2=(g_e-1/2)e^{-2(|x|-a_e)}$.

\paragraph{Construction and all state integrals.}
Every physical pair $0<\tau\le1$, $0<\chi<1$ determines
\begin{equation}
 g_0=\frac1{2(1-\chi)},\quad
 \delta=\frac{\tau\chi}{2(1-\chi)},\quad g_e=g_0-\delta,\quad
 \widetilde a=\delta^{-1}\sqrt{3\chi/\tau}.
\end{equation}
Conversely these parameters are recovered from $g_0>g_e\ge1/2$ by $\chi=(2g_0-1)/(2g_0)$ and $\tau=2(g_0-g_e)/(2g_0-1)$. With $g=g_e+\delta v$, \eqref{eq:Fg} becomes $F(g)=(2g_0-1)^2P(v)$. Define the even $g$ by $|x|=\int_{g(x)}^{g_0}F(r)^{-1/2}\dd r$, and $m$ by \eqref{eq:app-mstar}. Its simple central turn and double limiting zero give a $C^2$ solution on $\R$. On the support,
\begin{equation}
 x=\frac\tau2\int_{v(x)}^1P(v)^{-1/2}\dd v,\qquad
 m(x)=\frac{3\chi}{2s\tau}v(x)^2\quad(x\ge0).
\end{equation}
The support half-length is $a_e=(\tau/2)I[1]$, finite for $\tau<1$ and infinite for $\tau=1$. In the latter case the first integral gives exponential decay; hence all constructed densities are integrable, continuous and bounded.

The constructed profiles satisfy $s\lambda_1<1$. Define $\varphi_\pm=\sqrt g\exp(\mp\int_0^x\frac{\dd r}{2g(r)})$. Equation~\eqref{eq:GD} gives $\varphi_\pm''=(1-2sm)\varphi_\pm$ and Wronskian one. Their positive factorization gives a nonnegative Schr\"odinger form, hence $s\lambda_1\le1$. Equality would give an $H^1$ zero solution; square integrability at each infinity forces incompatible multiples of the independent $\varphi_\pm$, because $1/2\le g\le g_0$. Thus the inequality is strict and $g$ is the Green diagonal.

One substitution, $\dd x=(\tau/2)P^{-1/2}|\dd v|$ on each half, gives
\[
 s\fm(m)=\fm(\tau,\chi),\quad
 (2s)^{3/2}E(m)=\fE(\tau,\chi),\quad
 \sqrt{2s}\,\fr(m)=\sqrt{3\chi\tau}\,I[v],\quad
 \hD_s(m)=\fD(\tau,\chi).
\]
For the determinant use \eqref{eq:Dg}: the two exterior tails contribute $\log(2g_e)$, and the core contributes $\tau\chi I[L/N_\chi]$, before subtracting $s\fm$. The multiplier is $\mu=3\widetilde a/(4s\sqrt{2s})$, $\lambda_{\fm}=\mu s(2g_e-1)$, so its clipped coefficient is nonpositive outside the support. Every even decreasing minimizer is one of these profiles by the inverse first integral, and $\tau=1$ exactly when its mass multiplier is zero.

\paragraph{Parameter monotonicity and the envelope.}
Put $\Psi=I[v^3]/(\sqrt\tau I[v^2]^{3/2})$, so $\fE=(2\fm)^{3/2}\Psi$. The common bounds justify $C^1$ differentiation, including one-sided $\tau=1$, and give
\[
 \fm_\tau>0,\quad\fm_\chi>0,\quad
 (\log\Psi)_\tau=-\frac1{2\tau}+\langle\xi\rangle_3-\tfrac32\langle\xi\rangle_2<0,\quad
 (\log\Psi)_\chi=-\tfrac\tau2(\langle v^3/r\rangle_3-\tfrac32\langle v^3/r\rangle_2)\ge0.
\]
For the first shape sign, decreasing $\xi$ gives $\langle\xi\rangle_3\le\langle\xi\rangle_2$ by Chebyshev. For the second, write $h=v^3/r\ge0$; then $\langle h\rangle_3=\langle vh\rangle_2/\langle v\rangle_2\le(3/2)\langle h\rangle_2$. The mass derivatives follow directly from \eqref{eq:weight-derivatives}; for $\fm_\chi$ the integrand factor is $1-(\chi\tau/2)(v^3/r)>0$.

Along a constant-mass curve $\vartheta(\theta)=(\tau(\theta),\chi(\theta))$, or the free branch $\tau=1$, the exact identity is
\begin{equation}\label{eq:fred-envelope}
 \frac{\dd\fE}{\dd\theta}=\tfrac32\widetilde a\frac{\dd\fD}{\dd\theta}.
\end{equation}
We justify it despite the moving support. Let $m_\vartheta$ be the constructed density for $\vartheta=(\tau,\chi)$. For $0<w\le1$ put $X(w;\vartheta)=(\tau/2)\int_w^1P(v;\vartheta)^{-1/2}\dd v$; the common estimates give
\[
 |X(w;\vartheta')-X(w;\vartheta)|\le Cw^{-1}|\vartheta'-\vartheta|.
\]
Layer cake, integrated against $2w\dd w$, proves local $L^1$ Lipschitz dependence of $m_\vartheta$. Also $-X_w\ge\tau/(2\sqrt3)$ gives uniform convergence of $\sqrt m$, and $a_e$ is continuous for $\tau<1$. Consequently, for $\Delta=m_{\vartheta'}-m_\vartheta$ along a $C^1$ curve,
\[
 E(m_{\vartheta'})-E(m_\vartheta)
 =\mu[\hD_s(m_{\vartheta'})-\hD_s(m_\vartheta)]
       -\lambda_{\fm}\int_\R\Delta(x)\dd x+o(|\theta'-\theta|).
\]
The energy remainder is bounded by $(3/2)\nrm\Delta_1\nrm{\sqrt{m_{\vartheta'}}-\sqrt{m_\vartheta}}_\infty$, and the determinant remainder is quadratic by \eqref{eq:fred-variation}. The stationarity coefficient vanishes on the old support. Outside it this coefficient is Lipschitz and zero at the boundary, so the newly occupied layer contributes at most $C|a_e(\vartheta')-a_e(\vartheta)|\nrm\Delta_1=o(|\theta'-\theta|)$. On the free branch stationarity holds on the entire line. The mass term vanishes on the designated curves. Rescaling gives \eqref{eq:fred-envelope}, including differentiability of $\fD$.

\begin{proof}[Completion of \cref{thm:fred}]
On a fixed mass graph $\fm=sM$, the positive mass derivatives make $\chi=k(\tau)$ a $C^1$ decreasing function on an interval. The domain is an interval because $\sup_\chi\fm(\tau,\chi)$ is nondecreasing and $\fm(\tau,0)=0$. The signs above imply $\dd\fE/\dd\tau<0$, hence $\dd\fD/\dd\tau<0$ by the envelope. This proves uniqueness of the interior candidate. On the free branch $\fE(1,\chi)$ increases strictly, so does $\fD(1,\chi)$, proving its uniqueness.

Every listed candidate is admissible by the profile and state-integral construction. Conversely, an even decreasing minimizer exists and is subcritical and of that form. Complementarity places it either on the mass-active graph or the admissible free branch, with determinant exactly $\hd_s$. Thus at least one candidate occurs and the least energy among them is the infimum, as claimed.
\end{proof}

\subsection{Root-mass calibration}\label{app:root}

We solve the inner problem in \eqref{eq:root-main} by a global pointwise minimum principle.

\begin{theorem}[Explicit root-mass bound]\label{thm:dual}
Let $\eta\ge-34/5$, $a=\eta/2$, $y>0$, and for real $n\ge1$ put
\[
 p(n)=3n^2+an+1,\quad q(n)=n^2+an+3,\quad
 \psi(n)=\frac{n^3(n^2-1)}{p(n)^2},\qquad n_u=\psi^{-1}(y(1-u^2)).
\]
The inverse exists. Define
\begin{equation}
 \begin{aligned}
 \hat r&=2y\int_0^1\frac{p(n_u)^2}{n_u^2q(n_u)}\dd u,&
 \hat m&=2y\int_0^1\frac{p(n_u)^2}{n_u q(n_u)}\dd u,\\
 \hat E&=2y\int_0^1\frac{p(n_u)^2}{q(n_u)}\dd u,&
 \hat G&=4y^2\int_0^1\frac{p(n_u)^3}{n_u^2q(n_u)}\dd u.
 \end{aligned}
\end{equation}
The function $y\mapsto\hat G_\eta(y)$ is $C^1$, strictly increasing from zero to infinity. For $b>0$ let $y$ be the unique root $\hat G_\eta(y)=T_2/b^4$. Then
\[
 \cL_{\{2\}}^{\rm r}(\Lambda)\ge\ell_H(\Lambda;b,\eta)
 =b^3\hat E+\frac{\eta b}{4}(b^2\hat m-M)+b^2(b\hat r-1).
\]
\end{theorem}

We prove the inner minimization and change of variable for $\eta>-4\sqrt3$, and the monotonicity certificate for $\eta\ge-34/5$. In the auxiliary range $p\ge4+a>0$, $p'>0$, $q>0$, and
\begin{equation}\label{eq:psi-derivative}
 \psi'(n)=\frac{n^2(3n^2-1)q(n)}{p(n)^3}>0.
\end{equation}
Thus $\psi$ is a smooth bijection $[1,\infty)\to[0,\infty)$, with $\psi(n)\sim n/9$. Positivity of $q$ follows either from monotonicity for $a\ge-2$ or from completing its square for $-2\sqrt3<a<-2$.

\paragraph{Existence and activity of the inner minimum.}
Set $\lambda_{\fm}=\eta b/4$ and $f(z)=z^3+\lambda_{\fm}z^2+b^2z$. The inner functional is $\mathcal F(m)=\int_\R f(\sqrt{m(x)})\dd x$, and its infimum at $G_2\ge T>0$ is $V_{b,\eta}(T)$. The inequality $bz^2\le(z^3+b^2z)/2$ gives
\[
 f(z)\ge c_\eta(z^3+b^2z),\qquad
 c_\eta=1-\max(0,-\eta)/8>0.
\]
Narrow boxes show finiteness. For a minimizing sequence coercivity bounds $E$ and $\fr$, hence $\fm\le(\fr E)^{1/2}$. Rearrange and apply the compactness in \cref{app:foundation}; Fatou for the nonnegative continuous $f$ gives a minimizer. The constraint is active for every minimizer: if $G_2>T$, remove a positive-measure subset of the positive density with sufficiently small mass. Continuity of $G_2$ keeps admissibility, while $f>0$ makes the cost strictly smaller.

\paragraph{A global pointwise principle.}
Let $m_*$ be an even decreasing minimizer, $z_* =\sqrt{m_*}$, and $W_2=e^{-2|\cdot|}*m_*$. Put
\[
 K_0=\tfrac12\int_\R f'(z_*(x))z_*(x)\dd x,\qquad \mu=K_0/(2T).
\]
Then $z_*(x)$ globally minimizes $f(z)-2\mu W_2(x)z^2$ over $z\ge0$ for almost every $x$. Here is the constrained needle argument. At a common Lebesgue point $x_0$ of $z_*^j$, $j=1,2,3$, replace $m_*$ by a constant $\zeta\ge0$ on $B_\varepsilon=[x_0,x_0+\varepsilon]$, and multiply it by $1+\theta$ elsewhere. Bilinearity gives
\[
 G_2(m_{\varepsilon,\theta})=(1+\theta)^2T
 +2\varepsilon W_2(x_0)(\zeta-m_*(x_0))+O(\varepsilon|\theta|)+o(\varepsilon).
\]
Its derivative in $\theta$ is $2T+O(\varepsilon)+O(|\theta|)$, so exact equality to $T$ is obtained with
\[
 \theta=-\frac\varepsilon T W_2(x_0)(\zeta-m_*(x_0))+o(\varepsilon).
\]
The objective change is $\theta K_0+\varepsilon[f(\sqrt\zeta)-f(z_*(x_0))]+O(\theta^2)+O(\varepsilon|\theta|)+o(\varepsilon)$. Minimality, division by $\varepsilon$ and this correction give the claimed global inequality for every $\zeta$. One common full-measure set suffices, or first use rational $\zeta$ and continuity. Finally $\mu\le0$ would make this pointwise objective strictly positive for every $z>0$, since $f(z)>0$, forcing $m_*=0$. Therefore $\mu>0$.

\paragraph{The selected profile and its physical integral.}
Write $c=-\lambda_{\fm}+2\mu W_2$. This continuous even function decreases strictly on $x>0$ and tends to $-\eta b/4<2b$. The cubic $z^3-cz^2+b^2z$ has unique global minimizer zero when $c<2b$, and unique global minimizer $z_+(c)=(c+\sqrt{c^2-3b^2})/3$ when $c>2b$. At equality both zero and $b$ minimize it. Thus there is one finite half-support $x_*$, $c(x_*)=2b$, and $\sqrt{m_*}=b\hat n$ on the core, where $\hat n$ decreases continuously from $n_0>1$ to $1$; the density is zero outside. This proves the jump and chooses the larger root globally.

Put $\hat c=c/b$, $\varkappa=8\mu b$ and $y=1/\varkappa$. On the core,
\[
 \hat c=\frac{3\hat n^2+1}{2\hat n},\qquad
 \hat c''=4\hat c+\eta-\varkappa\hat n^2.
\]
Outside, $\hat c=-\eta/4+(2+\eta/4)e^{-2(x-x_*)}$ for $x\ge x_*$, hence $\hat c'(x_*)=-(4+a)=-p(1)$. Multiplying the core equation by $2\hat c'$ and using $(3\hat n^2-1)\hat n'=2\hat n^2\hat c'$ gives the zero first integral
\[
 \hat c'^2-(2\hat c+\eta/2)^2+\varkappa(\hat n^3-\hat n)=0,
 \qquad \hat n^2\hat c'^2=R_y(\hat n),
 \quad R_y(n)=p(n)^2-\frac{n^3(n^2-1)}y.
\]
The central turn is $n_0=\psi^{-1}(y)$; it is simple, since $R_y'(n_0)=-p(n_0)^2\psi'(n_0)/y<0$. Thus, for every continuous $F$ on $[1,n_0]$,
\begin{equation}
 \int_{-x_*}^{x_*}F(\hat n(x))\dd x
 =\langle F\rangle_y:=\int_1^{n_0}\frac{F(n)(3n^2-1)}{n\sqrt{R_y(n)}}\dd n
 =2y\int_0^1\frac{F(n_u)p(n_u)^2}{n_u^3q(n_u)}\dd u.
\end{equation}
The last substitution is $\psi(n)=y(1-u^2)$, with Jacobian from \eqref{eq:psi-derivative}; it removes the integrable turning singularity. Consequently
\[
 \begin{gathered}
 \fr(m_*)=b\hat r,\qquad\fm(m_*)=b^2\hat m,\qquad E(m_*)=b^3\hat E,\\
 G_2(m_*)=2b^4y\langle np(n)\rangle_y=b^4\hat G,\qquad
 \hat G=2y(3\hat E+a\hat m+\hat r).
 \end{gathered}
\]
For the fourth identity use $W_2=2b^2yp(\hat n)/\hat n$ on the core. Every even decreasing inner minimizer has these state integrals throughout the auxiliary range.

\paragraph{A compact monotonicity certificate.}
Let $\mathcal Q(n)=n^4(n^2-1)^2/[p(n)q(n)]$. For $t=n-1\ge0$ the numerator of its logarithmic derivative is
\[
 \begin{aligned}
 N_a(t)={}&12t^6+4(5a+18)t^5+2(3a^2+50a+120)t^4\\
 &+8(a+6)(3a+10)t^3+34(a+4)^2t^2+20(a+4)^2t+4(a+4)^2,
 \end{aligned}
\]
namely $N_a=n(n^2-1)pq(\log\mathcal Q)'$. Its derivative in $a$ has nonnegative coefficients for $a\ge-4$. Hence for $a\ge-17/5$, $N_a\ge\Pi_\sharp/25$, where
\[
 \Pi_\sharp(t)=300t^6+100t^5-766t^4-104t^3+306t^2+180t+36.
\]
For $z=(1,t,t^2,t^3)^\top$, $w=(1,t,t^2)^\top$ the exact identity
\[
 2\Pi_\sharp=z^\top S_0z+t\,w^\top S_1w,\quad
 S_0=\begin{pmatrix}72&81&4&-114\\81&624&-14&-538\\4&-14&116&-89\\-114&-538&-89&600\end{pmatrix},\quad
 S_1=\begin{pmatrix}198&-10&-144\\-10&336&-286\\-144&-286&378\end{pmatrix}
\]
proves positivity. Sylvester's criterion gives positive leading minors $72,38367,4417404,101099401$ and $198,66428,1123200$. Thus $\mathcal Q'>0$ on $(1,\infty)$.

Now $\hat G(y)=4\int_0^1\mathcal Q(n_u)/(1-u^2)^2\dd u$. Since $\psi'>0$, its derivative integrand is positive. Differentiation is legitimate at $u=1$: $\mathcal Q=(n-1)^2k_2(n)$ and $\mathcal Q'=(n-1)k_1(n)$ with smooth $k_i$, while $(n_u-1)/(1-u^2)$ extends continuously as $y/\psi'(1)$. This supplies continuous majorants on compact $y$ intervals and proves $C^1$ strict increase. Finally $p^3/(n^2q)$ is positive with a positive global minimum and tends to infinity; its boundedness on $[1,\psi^{-1}(y)]$ as $y\downarrow0$ gives $\hat G(0+)=0$, and $\hat G(y)\ge Cy^2$ gives divergence at infinity.

\begin{proof}[Completion of \cref{thm:dual}]
The inner minimizer exists and has active $G_2=T_2$. The physical formulas therefore select the unique $y$ in the theorem and give $V_{b,\eta}(T_2)=b^3\hat E+(\eta b/4)b^2\hat m+b^3\hat r$. Substituting this in the weak-duality formula \eqref{eq:root-main} proves the stated bound.
\end{proof}

\section{Analytic dependence and evaluation}\label{app:analytic}

\subsection{Analytic dependence}

We construct a common complex neighborhood on which all Bessel levels, the norming constants and the GL transformation depend holomorphically on the order. This proves the analytic input of \cref{thm:C}.

\begin{lemma}
For each real $\alpha_0>-1$ there is a disk $\mathbb D_r(\alpha_0)=\{\alpha:|\alpha-\alpha_0|<r\}$, $r<\alpha_0+1$, on which all $\sigma_k$ continue holomorphically, agree with $j_{\alpha,k}^2$ on the real slice, and satisfy uniformly
\[
 |\sigma_k-\tau_k|\le Ck^{-2},\qquad
 |\sigma_\ell-\sigma_k|\ge c|\ell^2-k^2|,\qquad |\sigma_k|\ge ck^2.
\]
\end{lemma}
\begin{proof}
On a closed order disk in $\Re\alpha>-1$, complex-order Hankel remainders \cite[10.17.13--15]{DLMF} give, uniformly in $|\arg z|\le\pi/4$,
\[
 (\pi z/2)^{1/2}J_\alpha(z)
 =P\cos\Theta-Q\sin\Theta,\quad
 P=1+O(z^{-2}),\quad Q=(4\alpha^2-1)/(8z)+O(z^{-3}),
\]
where $\Theta=z-\alpha\pi/2-\pi/4$. For $z=\varpi_k+\zeta$, $\varpi_k=\pi(k-1+\gamma)$, the zero equation is $P\sin\zeta+Q\cos\zeta=0$. Rouch\'e on $|\zeta|=1/2$ gives one simple zero for all sufficiently large $k$, on the same disk, and
\[
 z_k=\varpi_k-\frac{4\alpha^2-1}{8\varpi_k}+O(k^{-3}),\qquad
 z_k^2=\tau_k+O(k^{-2}).
\]
Uniqueness glues the local holomorphic zeros. McMahon identifies their numbering on the real slice. The finitely many remaining simple positive zeros continue by the implicit function theorem after one common shrinking. The reference gap $|\tau_\ell-\tau_k|\ge(\pi^2/2)|\ell^2-k^2|$ handles all pairs with a large index; shrinking again handles finitely many remaining pairs and keeps them away from zero.
\end{proof}

On a smaller disk continue $c_k^0=1/\fn_k$ and $c_k=\rho_k/\fn_k$ by their bilinear normalizations. The whole construction is holomorphic, with uniform bounds:
$\rho_k-1=O(k^{-2})$, $\sup|F_k|\le Ck^{-2}$, $\sup|K|\le C$, and $|u(t,\sigma_k)|\le Ck^{-\Re\gamma}$. Indeed the product proof in \cref{app:closure} uses the uniform gap and displacement estimates, which also hold for complex orders. Along the segment from $\tau_k$ to $\sigma_k$, $\sqrt{z-c_0}=\varpi_k+O(k^{-3})$ stays in one bounded imaginary strip. The complex reference estimates of \cref{app:reference} therefore give normally convergent $F$ in $C([0,1/2]^2)$.

At $\alpha_0$ every section is coercive by the same $\kappa_*$. Uniform kernel continuity preserves invertibility on one smaller disk, with inverse norm at most $2/\kappa_*$. On $L^2((0,1/2),\dd s)$ write $\cF=\cF_{\alpha,1/2}$, let $P_t$ be multiplication by $\ind_{(0,t)}$, and use
\[
 k_t=-(I+P_t\cF P_t)^{-1}P_tF(t,\cdot),\qquad
 K(t,s)=-F(t,s)-\int_0^{1/2}k_t(r)F(r,s)\dd r.
\]
Holomorphic inversion and this bounded bilinear pairing prove holomorphy of $K$ and $u$. The continued bilinear normalizations are the nonzero holomorphic functions $\fn_k$ of \cref{app:reference}. Thus the series $\sum_k(2/\sigma_k)(\rho_k/\fn_k)u(t,\sigma_k)^2$ is holomorphic with a uniform $Ck^{-2}$ majorant. It follows that $c_* =2\int_0^{1/2}\phi_\alpha(t)^2\dd t$ is holomorphic. On the real slice uniqueness identifies it with the original construction, proving real analyticity throughout $(-1,\infty)$.

\subsection{Numerical evaluation}\label{app:table}

Table~\ref{tab:main} illustrates the formulas at the ten stated indices. Values are displayed to six significant digits, with two digits for the relative gap.

\emph{Canonical cost.} Bessel zeros are obtained by bracketed root solves. The normalized reference waves use the Gegenbauer recurrence in \eqref{eq:reference-norm}. Gauss--Jacobi quadrature of \eqref{eq:app-fresc} evaluates the small difference $f(t,\sigma_k)-f(t,\tau_k)$; the difference of cosines is computed as a product of sines. With $R_k^0=f(\cdot,\tau_k)/\sqrt{\fn_k}$ and $D_k=\sqrt{c_k}f(\cdot,\sigma_k)-R_k^0$, stack the waves as columns of $R^0,D$ and form the discretized $F$ as $D(R^0)^\top+R^0D^\top+DD^\top$. This preserves the cancellation in its spectral terms. Each GL section is solved by positive-definite Nystr\"om discretization with $s=tv^2$, and the outer energy integral uses $t=y^2/2$. Subtract the reference Green diagonal before summing the profile:
\[
 \phi_\alpha(t)=2\mathfrak G_0(t,t)
 +\sum_k\left[\frac{2c_k}{\sigma_k}u(t,\sigma_k)^2
                    -\frac{2}{\tau_k\fn_k}f(t,\tau_k)^2\right].
\]
Here $\tau_1>0$ at every tabulated index. The reference Green kernel is evaluated from \cref{app:geometry} with 28-digit special functions. The triples (spectral cutoff, section and outer nodes, oscillatory nodes) are $(200,64,256)$, $(320,96,400)$ and $(480,128,576)$. The $U$ and sine values agree with their explicit reference values within $2\cdot10^{-13}$. The largest absolute change between the last two discretizations is $6.4\cdot10^{-9}$.

\emph{Lower values.} The second-moment and Keller roots are bracketed in their monotone shape parameters. For Fredholm, the free branch is tested first; when its mass exceeds the cap, solve the mass equation for $\chi$ and then the strictly monotone determinant equation for $\tau$. Search $s\Lambda_1$ on a grid reaching $0.9999$ and refine interior maxima. The small-$s$ target is evaluated with 40-digit Bessel data; on the free branch its quadrature is rearranged to factor out $\chi^2$, preserving the second-order cancellation. For the root bound, set $b=(T_2/\hat G_\eta(y))^{1/4}$ and search in $(\log y,\log(\eta+7))$, with $-32\le\log y\le12$ and $-34/5\le\eta\le2000$. Monotone safeguarded inversion gives $n_u$. Multiple initial points are followed by local optimization.

All lower quadratures use 512 nodes and are repeated with 1024 nodes at the selected parameters. The largest absolute changes are $5.2\cdot10^{-14}$ for $L_2$, $1.2\cdot10^{-13}$ for Keller, $2.7\cdot10^{-12}$ for Fredholm and $1.8\cdot10^{-8}$ for the root bound. The columns $\widehat B$ and $\widehat H$ are the resulting evaluations of $B_s$ and $\ell_H$; each admissible parameter gives an analytic lower bound. The numerical gaps compare these evaluations with the computed canonical cost. The certified Chebyshev value from Part~II appears separately in \cref{sec:comparison}.

\end{document}